\documentclass[12pt]{amsart}

\usepackage{enumitem}
\usepackage{psfrag}
\usepackage{graphicx}
\usepackage{pinlabel}
\usepackage{graphicx}
\usepackage{amsmath}
\usepackage{amssymb}
\usepackage{amscd}
\usepackage{autobreak}
\usepackage{picinpar}
\usepackage{times}
\usepackage{pb-diagram}
\usepackage{graphicx}
\usepackage{wrapfig}
\usepackage{xspace}
\usepackage{hyperref}
\usepackage{url}
\usepackage{caption}
\usepackage{subcaption}
\usepackage{fancyhdr}
\usepackage{overpic}
\usepackage{color}
\usepackage[square,comma,sort&compress,numbers]{natbib}
\usepackage{latexsym}
\usepackage{mathrsfs}
\usepackage{amsfonts}
\usepackage{hyperref}
\usepackage{amsthm}
\usepackage{appendix}
\usepackage{pdfsync}

\theoremstyle{definition}
\newtheorem{theorem}{Theorem}[section]
\newtheorem{lemma}[theorem]{Lemma}
\newtheorem{proposition}[theorem]{Proposition}
\newtheorem{corollary}[theorem]{Corollary}
\newtheorem{definition}[theorem]{Definition}

\newtheorem{remark}[theorem]{Remark}
\newtheorem*{theorem*}{Theorem}

\def\qed{\hfill{Q.E.D.}\smallskip}

\begin{document}

\title{\bf Discrete conformal structures on surfaces with boundary (II)---Rigidity and Existence}
\author{Xu Xu, Chao Zheng}

\date{\today}

\address{School of Mathematics and Statistics, Wuhan University, Wuhan, 430072, P.R.China}
 \email{xuxu2@whu.edu.cn}

\address{International Center for Mathematical Research (BICMR), Beijing (Peking) University, Beijing, 100871, P.R. China}
\email{czheng@whu.edu.cn}



\begin{abstract}
In \cite{X-Z DCS1}, we introduced discrete conformal structures on surfaces with boundary via an axiomatic framework, and provided a classification of such discrete conformal structures.
The present work focuses on the rigidity and existence of these discrete conformal structures on surfaces with boundary. 
As a direct consequence, the results by Guo-Luo in \cite{GL2} and Guo in \cite{Guo}, which deal with the rigidity and existence of discrete conformal structures on surfaces with boundary, are extended to a very general context.
\end{abstract}

\maketitle

\tableofcontents

\section{Introduction}\label{section 1}

\subsection{Basic definition and previous results}

Let $\widetilde{\mathcal{T}}=(V,\widetilde{E},\widetilde{F})$ be a triangulation of a closed surface $\widetilde{\Sigma}$, 
where $V$, $\widetilde{E}$, and $\widetilde{F}$ denote the sets of vertices, unoriented edges, and faces, respectively. 
Here $V$ is a finite subset of $\widetilde{\Sigma}$ with $|V|=N$. 
Let $N(V)$ be a small open regular neighborhood of the union of all vertices, chosen such that the neighborhoods around distinct vertices are disjoint. 
Then $\Sigma=\widetilde{\Sigma}\setminus N(V)$ is a compact surface with $N$ boundary components. 
The intersection $\mathcal{T}=\widetilde{\mathcal{T}}\cap \Sigma$ is referred to as an ideal triangulation of $\Sigma$. 
Let $E=\widetilde{E}\cap \Sigma$ and $F=\widetilde{F}\cap \Sigma$ denote the sets of unoriented ideal edges and ideal faces of $\Sigma$ with respect to $\mathcal{T}$, respectively. 
A boundary arc is defined as a connected component of the intersection between an ideal face and the boundary $\partial\Sigma$. 
For simplicity, we index the boundary components of $\Sigma$ by the set $B = \{1, 2, \dots, N\}$, where each element $i \in B$ corresponds to a unique boundary component. 
Let $E_+$ denote the set of oriented ideal edges. 
For any two adjacent boundary components $i, j \in B$, the unoriented ideal edge connecting them is denoted by $\{ij\} \in E$, with the corresponding oriented ideal edge given by $(i, j)$. 
An ideal face $\{ijk\} \in F$ is adjacent to the boundary components $i, j, k \in B$.
Specifically, such an ideal face is a hexagon, which corresponds to the triangle $v_iv_jv_k\in \widetilde{F}$ in $\widetilde{\mathcal{T}}$. 
The sets of real-valued functions on $B$, $E$, and $E_+$ are denoted by $\mathbb{R}^N$, $\mathbb{R}^E$, and $\mathbb{R}^{E_+}$, respectively.

A edge length function associated to $\mathcal{T}$ is a vector $l \in \mathbb{R}^E_{>0}$ that assigns to each ideal edge $\{ij\} \in E$ a positive number $l_{ij} = l_{ji}$. 
For any ideal face $\{ijk\} \in F$, there exists a unique right-angled hyperbolic hexagon (up to isometry) whose three non-adjacent edges have lengths $l_{ij}, l_{ik}, l_{jk}$ (see \cite[Theorem 3.5.14]{Ratcliffe}). 
By gluing all such right-angled hyperbolic hexagons isometrically along the ideal edges in pairs,
one can construct a hyperbolic surface with totally geodesic boundary from the ideal triangulation $\mathcal{T}$. 
Conversely, any ideally triangulated hyperbolic surface $(\Sigma, \mathcal{T})$ with totally geodesic boundary induces a function $l \in \mathbb{R}^E_{>0}$, where $l_{ij}$ is the length of the shortest geodesic connecting the boundary components $i, j \in B$. Such an edge length function $l \in \mathbb{R}^E_{>0}$ is referred to as a \textit{discrete hyperbolic metric} on $(\Sigma, \mathcal{T})$.
The length $K_i$ of the boundary component $i \in B$ is called the \textit{generalized combinatorial curvature} of the discrete hyperbolic metric $l$ at $i$. Specifically, the generalized combinatorial curvature $K_i$ is defined by
\begin{equation}\label{Eq: K}
K_i = \sum_{\{ijk\} \in F} \theta^{jk}_i,
\end{equation}
where the summation is taken over all right-angled hyperbolic hexagons adjacent to $i$, and the generalized angle $\theta^{jk}_i$ denotes the length of the boundary arc of the right-angled hyperbolic hexagon $\{ijk\}$ at $i$.

Motivated by the work of Glickenstein \cite{Glickenstein} and that of Glickenstein-Thomas \cite{GT} on discrete conformal structures for triangulated closed surfaces, 
we introduced in \cite{X-Z DCS1} the following two definitions regarding partial edge lengths and discrete conformal structures on ideally triangulated surfaces with boundary.

\begin{definition}[\cite{X-Z DCS1}, Definition 1.1]\label{Def: partial edge length}
Let $(\Sigma,\mathcal{T})$ be an ideally triangulated surface with boundary.
An assignment of partial edge lengths is a map $d\in \mathbb{R}^{E_{+}}$ satisfying $l_{ij}=d_{ij}+d_{ji}>0$ for every edge $\{ij\}\in E$, and
\begin{equation}\label{Eq: compatible condition}
\sinh d_{ij} \sinh d_{jk} \sinh d_{ki}=
\sinh d_{ji} \sinh d_{kj} \sinh d_{ik}
\end{equation}
for every ideal face $\{ijk\}\in F$.
\end{definition}

\begin{definition}[\cite{X-Z DCS1}, Definition 1.2]\label{Def: DCS}
Let $(\Sigma,\mathcal{T})$ be an ideally triangulated surface with boundary.
A discrete conformal structure $d=d(f)$ on $(\Sigma,\mathcal{T})$ is a smooth map that sends a function $f\in \mathbb{R}^N$ defined on the boundary components $B$ to a partial edge length function $d\in \mathbb{R}^{E_+}$, satisfying
\begin{equation}\label{Eq: variation 1}
\frac{\partial l_{ij}}{\partial f_i}
=\coth d_{ij}
\end{equation}
for each $(i,j)\in E_{+}$, and
\begin{equation}\label{Eq: variation 2}
\frac{\partial d_{ij}}{\partial f_k}=0
\end{equation}
if $k\neq i$ and $k\neq j$.
The function $f\in \mathbb{R}^N$ is called a discrete conformal factor.
\end{definition}

\begin{remark}\label{Rmk: 5}
In general, the partial edge lengths $d \in \mathbb{R}^{E_+}$ may take negative values and do not satisfy the symmetry condition, i.e., $d_{ij} \neq d_{ji}$. 
Let $E_{ij}$ denote the hyperbolic geodesic in the hyperbolic plane $\mathbb{H}^2$ that extends the ideal edge $\{ij\} \in E$. 
The partial edge length $d_{ij}$ represents the signed distance from a point $c_{ij} \in E_{ij}$ to the boundary component $i \in B$, while $d_{ji}$ represents the signed distance from $c_{ij}$ to $j \in B$. 
This point $c_{ij}$ is referred to as the \emph{edge center} of the ideal edge $\{ij\}$. 
Note that the signed distance $d_{ij}$ from $c_{ij}$ to $i$ is positive if $c_{ij}$ lies on the same side of $i$ as $j$ along the hyperbolic geodesic $E_{ij}$, and negative otherwise.
Please refer to Figure \ref{Figure 1}.
\begin{figure}[!ht]
  \centering
  \includegraphics[scale=0.9]{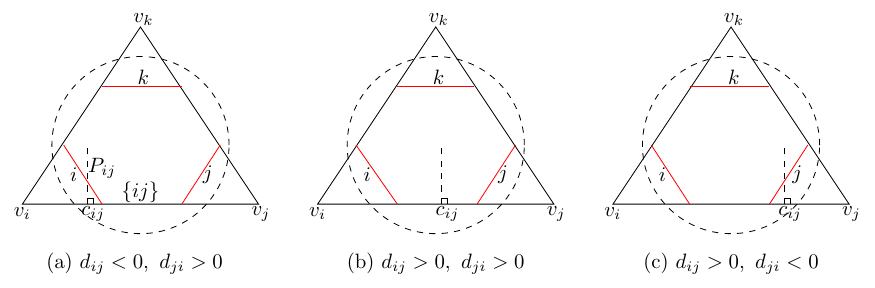}
  \caption{Poincar\'{e} dual to the ideal edge $\{ij\}\in E$ with the position determined by $d_{ij}$ and $d_{ji}$ in the Klein model (dotted). The red lines represent the boundary arcs of the ideal face $\{ijk\}$, and the black triangle $v_iv_jv_k$ is the corresponding  hyper-ideal hyperbolic triangle.}
\label{Figure 1}
\end{figure}
\end{remark}

\begin{remark}\label{Rmk: GM}
We briefly elaborate on the motivation behind Definition \ref{Def: partial edge length} concerning partial edge lengths on surfaces with boundary; further details can be found in \cite{X-Z DCS1}. 
The hyperbolic geodesic passing through $c_{ij}$ and Lorentz orthogonal to $E_{ij}$ is referred to as the \emph{edge perpendicular}, denoted by $P_{ij}$. 
This edge perpendicular $P_{ij}$ can be identified with a unique 2-dimensional time-like subspace $\widetilde{P}_{ij}$ of the Lorentzian space $\mathbb{R}^3$, satisfying $\widetilde{P}_{ij} \cap \mathbb{H}^2 = P_{ij}$. 
Since the 2-dimensional subspaces $\widetilde{P}_{ij}$ and $\widetilde{P}_{jk}$ (associated to $P_{ij}$ and $P_{jk}$, respectively) intersect in a 1-dimensional subspace of $\mathbb{R}^3$, 
their intersection $\widetilde{P}_{ij} \cap \widetilde{P}_{jk}$ corresponds to a point $q$ in the Klein model of $\mathbb{H}^2$. 
This point $q$ is viewed as the intersection of the perpendiculars $P_{ij}$ and $P_{jk}$. 
Note that $q$ may be not in $\mathbb{H}^2$, meaning the perpendiculars $P_{ij}, P_{jk}$ may not intersect within $\mathbb{H}^2$.
Lemma 2.4 in \cite{X-Z DCS1} demonstrates that the perpendiculars $P_{ij}$, $P_{jk}$, and $P_{ki}$ of the ideal face $\{ijk\}$ intersect at a common point $c_{ijk}$ (termed the \emph{face center}) if and only if the condition (\ref{Eq: compatible condition}) is satisfied. 
Consequently, the condition (\ref{Eq: compatible condition}) in Definition \ref{Def: partial edge length} guarantees the existence of a geometric structure on the Poincar\'{e} dual of an ideally triangulated surface with boundary.
Moreover, the geometric motivation behind Definition \ref{Def: DCS} is discussed in detail in \cite[Remark 1.3]{X-Z DCS1}. 
\end{remark}

The main result in \cite{X-Z DCS1} is the following theorem, which provides explicit forms and a classification of the discrete conformal structures (as defined in Definition \ref{Def: DCS}) on surfaces with boundary.

\begin{theorem}[\cite{X-Z DCS1}, Theorem 1.4]\label{Thm: DCS}
Let $(\Sigma,\mathcal{T})$ be an ideally triangulated surface with boundary, and let $d=d(f)$ be a discrete conformal structure on $(\Sigma,\mathcal{T})$.
There exist constant vectors $\alpha\in \mathbb{R}^N$ and $\eta\in \mathbb{R}^{E}$ satisfying $\eta_{ij}=\eta_{ji}$, such that for any right-angled hyperbolic hexagon $\{ijk\}\in F$,
the value of $\frac{\sinh d_{ij}}{\sinh d_{ji}}$ can only take one of the following four forms:
\begin{gather*}
\sqrt{\frac{1+\alpha_ie^{2f_i}}{1+\alpha_je^{2f_j}}},\
-\sqrt{\frac{1+\alpha_ie^{2f_i}}{1+\alpha_je^{2f_j}}},\
e^{\frac{1}{2}C_{ij}(f_i-f_j)},\
-e^{\frac{1}{2}C_{ij}(f_i-f_j)},
\end{gather*}
where $C\in \mathbb{R}^{E_+}$ is a constant vector satisfying
$C_{ij}+C_{jk}+C_{ki}=0$ for any $\{ijk\}\in F$, and $C_{rs}+C_{sr}=0$ for any subset $\{r,s\}\subseteq\{i,j,k\}$.
Furthermore,
\begin{description}
\item[(A1)]
if $\frac{\sinh d_{ij}}{\sinh d_{ji}}=\sqrt{\frac{1+\alpha_ie^{2f_i}}{1+\alpha_je^{2f_j}}}>0$ with $1+\alpha_ie^{2f_i}>0$ and $1+\alpha_je^{2f_j}>0$,
then
\begin{equation}\label{Eq: DCS3}
\cosh l_{ij}
=-\sqrt{(1+\alpha_ie^{2f_i})(1+\alpha_je^{2f_j})}
+\eta_{ij}e^{f_i+f_j};
\end{equation}

\item[(A2)]
if $\frac{\sinh d_{ij}}{\sinh d_{ji}}
=\sqrt{\frac{1+\alpha_ie^{2f_i}}{1+\alpha_je^{2f_j}}}>0$ with $1+\alpha_ie^{2f_i}<0$ and $1+\alpha_je^{2f_j}<0$,
then
\begin{equation}\label{Eq: new 1}
\cosh l_{ij}
=\sqrt{(1+\alpha_ie^{2f_i})(1+\alpha_je^{2f_j})}
+\eta_{ij}e^{f_i+f_j};
\end{equation}

\item[(A3)]
if $\frac{\sinh d_{ij}}{\sinh d_{ji}}
=e^{\frac{1}{2}C_{ij}(f_i-f_j)}>0$,
then
\begin{equation}\label{Eq: DCS1}
\cosh l_{ij}
=-\cosh(f_j-f_i-C_{ij})+\eta_{ij}e^{f_i+f_j};
\end{equation}

\item[(B1)] if $\frac{\sinh d_{ij}}{\sinh d_{ji}}
=-\sqrt{\frac{1+\alpha_ie^{2f_i}}{1+\alpha_je^{2f_j}}}<0$
with $1+\alpha_ie^{2f_i}>0$ and $1+\alpha_je^{2f_j}>0$,
then
\begin{equation}\label{Eq: DCS4}
\cosh l_{ij}
=\sqrt{(1+\alpha_ie^{2f_i})(1+\alpha_je^{2f_j})}
+\eta_{ij}e^{f_i+f_j};
\end{equation}

\item[(B2)]
if $\frac{\sinh d_{ij}}{\sinh d_{ji}}
=-\sqrt{\frac{1+\alpha_ie^{2f_i}}{1+\alpha_je^{2f_j}}}<0$
with $1+\alpha_ie^{2f_i}<0$ and $1+\alpha_je^{2f_j}<0$,
then
\begin{equation}\label{Eq: new 2}
\cosh l_{ij}
=-\sqrt{(1+\alpha_ie^{2f_i})(1+\alpha_je^{2f_j})}
+\eta_{ij}e^{f_i+f_j};
\end{equation}

\item[(B3)]
if $\frac{\sinh d_{ij}}{\sinh d_{ji}}
=-e^{\frac{1}{2}C_{ij}(f_i-f_j)}<0$,
then
\begin{equation}\label{Eq: DCS2}
\cosh l_{ij}
=\cosh(f_j-f_i-C_{ij})+\eta_{ij}e^{f_i+f_j}.
\end{equation}
\end{description}
Moreover, there are exactly six types of combinations of such discrete conformal structures that can exist on the surface with boundary, namely (A1), (A2), (A3), the combination of (A1) and (B1), the combination of (A2) and (B2), and the combination of (A3) and (B3). 
The latter three mixed combinations are referred to as the \emph{mixed discrete conformal structure I}, \emph{mixed discrete conformal structure II}, and \emph{mixed discrete conformal structure III}, respectively.
\end{theorem}

\begin{remark}\label{Rmk: 1}
As noted in \cite[Remark 1.5]{X-Z DCS1}, for the discrete conformal structures (A1) and (B1), the discrete conformal structure $d = d(f)$ allows a reparameterization such that $\alpha: B \to \{-1, 0, 1\}$, with the induced discrete hyperbolic metric $l$ remaining invariant. 
Thus we always assume $\alpha: B \to \{-1, 0, 1\}$ in (A1) and (B1) .
For the discrete conformal structures (A2) and (B2), the conditions $1 + \alpha_i e^{2f_i} < 0$ and $1 + \alpha_j e^{2f_j} < 0$ imply $\alpha_i < 0$ and $\alpha_j < 0$. 
Via reparameterization, we can assume $\alpha \equiv -1$ in (A2) and (B2).
For the discrete conformal structures (A3) and (B3), if $(\Sigma, \mathcal{T})$ has genus zero, a reparameterization allows us to set $C \equiv 0$.
For surfaces of non-zero genus, we take $C \equiv 0$ in (A3) and (B3) as special cases.
\end{remark}

\begin{remark}\label{Rmk: 3}
The right-angled hyperbolic hexagon $\{ijk\} \in F$ with edge lengths $l_{ij}, l_{jk}, l_{ki}$ induced by (\ref{Eq: DCS3}) and (\ref{Eq: DCS4}), (\ref{Eq: new 1}) and (\ref{Eq: new 2}), and (\ref{Eq: DCS1}) and (\ref{Eq: DCS2}) is referred to as \emph{Type-I}, \emph{Type-II}, and \emph{Type-III}, respectively.
Let $\{r, s, t\} = \{i, j, k\}$. 
For a Type-I right-angled hyperbolic hexagon, if $d_{rs} < 0$, then $d_{rt} < 0$ by Lemma \ref{Lem: time-space-like}.
Furthermore, from the conditions $l_{rs} = d_{rs} + d_{sr} > 0$ and $l_{rt} = d_{rt} + d_{tr} > 0$ in Definition \ref{Def: partial edge length}, it follows that $d_{sr} > 0$ and $d_{tr} > 0$. 
For the condition (\ref{Eq: compatible condition}) to hold, $d_{st} > 0$ and $d_{ts} > 0$ are required. 
Thus $l_{rs}$ and $l_{rt}$ are defined by (\ref{Eq: DCS4}), while $l_{st}$ is defined by (\ref{Eq: DCS3}).
Analogous reasoning applies to Type-II and Type-III right-angled hyperbolic hexagons. 
Consequently, for the mixed discrete conformal structure I, all right-angled hyperbolic hexagons are either Type-I or a combination of Type-I right-angled hyperbolic hexagons and those with edge lengths induced by (\ref{Eq: DCS3}). 
Similar arguments also apply for the mixed discrete conformal structures II and III.
\end{remark}

\subsection{Main results}

A central problem in discrete conformal geometry concerns the relationship between discrete conformal structures and their combinatorial curvatures. 
The main result in this paper is the following theorem, which establishes the global rigidity and existence of discrete conformal structures on ideally triangulated surfaces with boundary.

\begin{theorem}\label{Thm: rigidity and image}
Let $(\Sigma,\mathcal{T})$ be an ideally triangulated surface with boundary, and let $d=d(f)$ be a discrete conformal structure on $(\Sigma,\mathcal{T})$.
\begin{description}
\item[(i)]
For the discrete conformal structure (A1), let $\alpha: B\rightarrow \{-1,0,1\}$ and $\eta\in \mathbb{R}_{>0}^E$ be the weights on $(\Sigma,\mathcal{T})$ satisfying $\eta_{ij}>\alpha_i\alpha_j$ for any two adjacent boundary components $i,j\in B$.
Then the discrete conformal factor $f$ is uniquely determined by its generalized combinatorial curvature $K\in \mathbb{R}^N_{>0}$.
Furthermore, if $\alpha: B\rightarrow \{0,1\}$, then
the image of $K$ is $\mathbb{R}^N_{>0}$.

\item[(ii)]
For the discrete conformal structure (A2), let $\eta\in (-1,+\infty)^E$ be the weight on $(\Sigma,\mathcal{T})$.
Then the discrete conformal factor $f$ is uniquely determined by its generalized combinatorial curvature $K\in \mathbb{R}^N_{>0}$.
Furthermore, if $\eta\in (-1,0]^E$, then
the image of $K$ is $\mathbb{R}^N_{>0}$.

\item[(iii)]
For the discrete conformal structure (A3), let $\eta\in \mathbb{R}_{>0}^E$ be the weight on $(\Sigma,\mathcal{T})$.
Then the discrete conformal factor $f$ is uniquely determined by its generalized combinatorial curvature $K\in \mathbb{R}^N_{>0}$.
Furthermore, the image of $K$ is $\mathbb{R}^N_{>0}$.

\item[(iv)]
For the mixed discrete conformal structure I,
let $\alpha: B\rightarrow \{-1,0,1\}$ and $\eta\in (1,+\infty)^E$ be the weights on $(\Sigma,\mathcal{T})$.
Then the discrete conformal factor $f$ is uniquely determined by its generalized combinatorial curvature $K\in \mathbb{R}^N_{>0}$.
Furthermore, if $\alpha: B\rightarrow \{0,1\}$, and $\alpha_j$ and $\alpha_k$ cannot both be 1 simultaneously for any right-angled hyperbolic hexagon with edge lengths $l_{ij},l_{ik}$ given by (\ref{Eq: DCS4}) and $l_{jk}$ given by (\ref{Eq: DCS3}),
then the image of $K$ is $\mathbb{R}^N_{>0}$.

\item[(v)]
For the mixed discrete conformal structure II, let $\eta\in [1,+\infty)^E$ be the weight on $(\Sigma,\mathcal{T})$.
Then the discrete conformal factor $f$ is uniquely determined by its generalized combinatorial curvature $K\in \mathbb{R}^N_{>0}$.

\item[(vi)]
For the mixed discrete conformal structure III, let $\eta\in \mathbb{R}_{>0}^E$ be the weight on $(\Sigma,\mathcal{T})$.
Then the discrete conformal factor $f$ is uniquely determined by its generalized combinatorial curvature $K\in \mathbb{R}^N_{>0}$.
Furthermore, the image of $K$ is $\mathbb{R}^N_{>0}$.
\end{description}
\end{theorem}

\begin{remark}
If $\alpha \equiv 0$, then Theorem \ref{Thm: rigidity and image} (i) generalizes the results obtained by Guo \cite{Guo} and Li-Xu-Zhou \cite{Li-Xu-Zhou}. 
For $\alpha \equiv 1$, Theorem \ref{Thm: rigidity and image} (i) was established by Guo-Luo \cite{GL2}. 
Theorem \ref{Thm: rigidity and image} (i) unifies these existing results and further incorporates mixed-type cases where $\alpha_i = 1$ for some non-empty subset $V_1 \subseteq V$ and $\alpha_j = 0$ for $j \in V \setminus V_1$.
Theorem \ref{Thm: rigidity and image} (ii) generalizes the results of Guo-Luo \cite{GL2}, while Theorem \ref{Thm: rigidity and image} (iii), which was proven by Guo-Luo \cite{GL2}, is included here for completeness.
In Theorem \ref{Thm: rigidity and image} (iv), the weight $\eta \in (1, +\infty)^E$ admits an extension to a broader range, with details provided in Subsection \ref{subsection: AS6}. 
For Theorem \ref{Thm: rigidity and image} (v), the weight $\eta \in [1, +\infty)^E$ can also be extended to a broader range (see Theorem \ref{Thm: Rigidity II}). 
Note that only a rigidity result is presented here, without any existence result; the reason is elaborated in Remark \ref{Rmk: 6}. 
In Theorem \ref{Thm: rigidity and image} (vi), the weight $\eta \in \mathbb{R}_{>0}^E$ similarly allows for extension to a broader range, as discussed in Theorem \ref{Thm: ASC 2}.
\end{remark}

\subsection{Basic ideas of the proof of Theorem \ref{Thm: rigidity and image}}\label{Subsec: 1}

The proof of the existence of discrete conformal structures relies on the continuity method. 
Specifically, we will show that the image of the curvature map $K$ is both an open and a closed subset of $\mathbb{R}^N_{>0}$; consequently, it coincides with the entire space $\mathbb{R}^N_{>0}$.

The proof of the rigidity of discrete conformal structures involves a variational principle introduced by Colin de Verdi\`{e}re \cite{De} for tangential circle packings on triangulated closed surfaces. 
The rigidity proof proceeds in three steps:
First, we characterize the admissible space of discrete conformal factors $f$ for a right-angled hyperbolic hexagon. 
By a variable transformation $u_i=u_i(f_i)$ (as detailed in the proof of Lemma \ref{Lem: angle variations 3}), the admissible space $\Omega_{ijk}$ of $f$ is transformed into the admissible space $\mathcal{U}_{ijk}$ of $u$, which is shown to be convex and simply connected.
Second, we prove that the Jacobian of the generalized angles $\theta$ with respect to $u$ for a right-angled hyperbolic hexagon is symmetric and negative definite.
Consequently, the Jacobian of the generalized combinatorial curvature $K$ with respect to $u$ is symmetric and negative definite on the admissible space $\mathcal{U}=\cap\, \mathcal{U}_{ijk}$.
These two steps collectively enable the definition of a strictly concave function on $\mathcal{U}$.
Third, we invoke the following well-known result from analysis to establish rigidity.

\begin{lemma}\label{Lem: analysis}
If $f:\Omega \rightarrow \mathbb{R}$ is a $C^1$-smooth  strictly concave function on an open convex set $\Omega \subset \mathbb{R}^n$,
then its gradient $\nabla f:\Omega \rightarrow \mathbb{R}^n$ is injective.
Furthermore, $\nabla f$ is a smooth embedding.
\end{lemma}

A critical component of the proof involves constructing the necessary variable transformation $u_i=u_i(f_i)$. 
To this end, we require the following result  concerning the variation of generalized angles.

\begin{theorem}\label{Thm: angle variations}
Let $(\Sigma,\mathcal{T})$ be an ideally triangulated surface with boundary, and let $d=d(f)$ be a discrete conformal structure on $(\Sigma,\mathcal{T})$.
\begin{description}
  \item[(i)] For any right-angled hyperbolic hexagon $\{ijk\}\in F$, the following holds:
\begin{equation*}
\frac{\partial \theta^{jk}_i}{\partial f_j}
=\frac{-1}{\sinh d_{ji}}\frac{\tanh^\beta h_{ij}}{\sinh l_{ij}},
\end{equation*}
where $\beta=1$ if the face center $c_{ijk}$ is time-like, and $\beta=-1$ if $c_{ijk}$ is space-like.
If $c_{ijk}$ is light-like, the formula is interpreted as $\tanh^\beta h_{ij}=\tanh^\beta h_{jk}=\tanh^\beta h_{ki}=1$.
Here, $h_{ij}$ denotes the signed distance from the face center $c_{ijk}$ to the hyperbolic geodesic line $E_{ij}$,
which is positive if $c_{ijk}$ lies on the same side of $E_{ij}$ as the right-angled hyperbolic hexagon $\{ijk\}$, negative otherwise, and zero if $c_{ijk}$ lies on $E_{ij}$.

\item[(ii)]
For a right-angled hyperbolic hexagon $\{ijk\}\in F$, the variation of generalized angles satisfies the following relation: 
\begin{equation*}
\frac{\partial \theta^{jk}_i}{\partial f_i}
=\cosh l_{ij}\frac{\partial \theta^{ik}_j}{\partial f_i}
+\cosh l_{ik}\frac{\partial \theta^{ij}_k}{\partial f_i}.
\end{equation*}
  \item[(iii)]
There exists a change of variables given by $u_i=u_i(f_i)$ such that
\begin{equation*}
\frac{\partial \theta^{jk}_i}{\partial u_j}
=\frac{\partial \theta^{ik}_j}{\partial u_i}.
\end{equation*}
The function $u$ is also referred to as a discrete conformal factor.
\end{description}
\end{theorem}

\subsection{Organization of the paper}

In Section \ref{section 2}, we establish Theorem \ref{Thm: angle variations}. 
In Section \ref{section A1}, we investigate the discrete conformal structure (A1) and prove Theorem \ref{Thm: rigidity and image} (i). 
Section \ref{section A2} is devoted to the discrete conformal structure (A2), where Theorem \ref{Thm: rigidity and image} (ii) is proven. 
In Section \ref{section A3}, we study the discrete conformal structure (A3) and establish Theorem \ref{Thm: rigidity and image} (iii). 
Section \ref{section III} focuses on the mixed discrete conformal structure III, with a proof of a generalization of Theorem \ref{Thm: rigidity and image} (vi). 
Similarly, Section \ref{section II} investigates the mixed discrete conformal structure II and proves a generalization of Theorem \ref{Thm: rigidity and image} (v). 
Finally, Section \ref{section I} studies the mixed discrete conformal structure I and establishes a generalization of Theorem \ref{Thm: rigidity and image} (iv).
\\
\\
\textbf{Acknowledgements}\\[8pt]
The research of X. Xu is supported by National Natural Science Foundation of China under grant no. 12471057.

\section{Discrete conformal variations of generalized angles}\label{section 2}

In this section, we first recall some fundamental results from 2-dimensional hyperbolic geometry (see \cite[Chapters 3 and 6]{Ratcliffe}). 
We then provide a geometric interpretation along with several useful lemmas. 
Finally, we proceed to prove Theorem \ref{Thm: angle variations}.

\subsection{Some fundamental results}

For any vectors $x,y$ in the Lorentzian space $\mathbb{R}^3$,
the Lorentzian inner product $*$ of $x$ and $y$ is defined as $x*y=xJy^\mathrm{T}$,
where $J=\mathrm{diag} \{1,1, -1\}$.
Two vectors $x,y\in \mathbb{R}^3$ are said to be Lorentz orthogonal if $x*y=0$.
The Lorentzian cross product $\otimes$ of $x$ and $y$ is defined as $x\otimes y=J(x\times y)$,
where $\times$ denotes the Euclidean cross product.

The Lorentzian norm of a vector $x\in \mathbb{R}^3$ is defined as the complex number $\Vert x \Vert=\sqrt{x*x}$.
Here $\Vert x \Vert$ is either positive, zero, or positive imaginary.
When $\Vert x \Vert$ is positive imaginary, its absolute value is denoted by $\vert\Vert x \Vert\vert$.
Furthermore, a vector $x\in \mathbb{R}^3$ is termed space-like if $x*x>0$, light-like if $x*x=0$, and time-like if $x*x<0$.
The hyperboloid model of $\mathbb{H}^2$ is defined as
\begin{equation*}
\mathbb{H}^2=\{x=(x_1,x_2,x_3)\in \mathbb{R}^3\mid x*x=-1\ \text{and}\ x_3>0 \},
\end{equation*}
which is embedded in the Lorentzian space $(\mathbb{R}^3,*)$.
A vector subspace of $\mathbb{R}^3$ is termed time-like if it contains at least one time-like vector.
A geodesic (or hyperbolic line) in $\mathbb{H}^2$ is the non-empty intersection of $\mathbb{H}^2$ with a 2-dimensional time-like vector subspace of $\mathbb{R}^3$.

The Klein model of $\mathbb{H}^2$ is defined as the central projection of $\mathbb{H}^2$ onto the plane $x_3=1$ in $\mathbb{R}^3$. 
In this model, the interior, boundary, and exterior points of the unit disk correspond to time-like, light-like, and space-like vectors, respectively, and are referred to as hyperbolic points, ideal points, and hyper-ideal points, correspondingly.  
In the Klein model, a geodesic in $\mathbb{H}^2$ is represented as an open chord of the unit disk.
A triangle $xyz$ with all three vertices $x,y,z$ being hyperbolic is called a hyperbolic triangle.
If at least one of $x,y,z$ is ideal or hyper-ideal,
then the triangle $xyz$ is termed a generalized hyperbolic triangle.
Specifically, if $x,y,z$ are hyper-ideal, this triangle is referred to as a hyper-ideal hyperbolic triangle for brevity.
It should be noted that, in this paper, each edge of a generalized hyperbolic triangle $xyz$ is required to intersect the unit disk.

The following proposition, which follows from properties of the Lorentzian cross product, is used extensively in this paper.

\begin{proposition}\label{Prop: Right angle}
Let $xyz$ be a (hyperbolic or generalized hyperbolic) triangle with $x\in\mathbb{H}^2$ and a right angle at $x$. Then
\begin{equation*}
-(z*y)=(z*x)(x*y).
\end{equation*}
\end{proposition}

For any two vectors $x,y\in \mathbb{H}^2$, there exists a unique geodesic passing through $x$ and $y$.
The arc length of the geodesic segment between $x$ and $y$ is defined as the hyperbolic distance between $x$ and $y$, denoted by $d_{\mathbb{H}}(x,y)$.
We use $\textrm{Span}(x)$ and $\textrm{Span}(x,y)$ to
denote the vector subspace spanned by $x$ and by $x,y$, respectively.
For a space-like vector $x\in \mathbb{R}^3$,
the Lorentzian complement of $\textrm{Span}(x)$, defined as $x^{\bot}:=\{w\in \mathbb{R}^3 \mid x*w=0\}$, is a 2-dimensional time-like vector subspace of $\mathbb{R}^3$, and thus $x^{\bot}\cap \mathbb{H}^2$ is a geodesic in $\mathbb{H}^2$.
The hyperbolic distance between $x$ and $y\in \mathbb{H}^2$ is defined as $d_{\mathbb{H}}(y,x^\bot)=\inf\{d_{\mathbb{H}}(y,z)\mid z\in x^\bot\cap \mathbb{H}^2\}$.
Let $x,y\in \mathbb{R}^3$ be space-like vectors.
If $\textrm{Span}(x,y)\cap \mathbb{H}^2$ is non-empty, then the hyperbolic distance between $x$ and $y$ is defined by $d_{\mathbb{H}}(x^\bot,y^\bot)=\inf\{d_{\mathbb{H}}(z,w)\mid z\in x^\bot\cap \mathbb{H}^2, w\in y^\bot\cap \mathbb{H}^2\}$.
If $\textrm{Span}(x,y)\cap \mathbb{H}^2$ is empty, then $x^\bot$ and $y^\bot$ intersect at an angle $\angle(x^\bot,y^\bot)$ in $\mathbb{H}^2$.

\begin{proposition}[\cite{Ratcliffe}, Chapter 3]\label{Prop: Ratcliffe 1}
Suppose $x,y\in\mathbb{H}^2$ are time-like vectors, and $z,w\in \mathbb{R}^3$ are space-like vectors. Then
\begin{description}
\item[(i)] $x*y=-\vert\Vert x\Vert\vert\cdot\vert\Vert y\Vert\vert\cosh d_{\mathbb{H}}(x,y)$;
\item[(ii)] $\vert x*z\vert=\vert\Vert x\Vert\vert\cdot \Vert z\Vert \sinh d_{\mathbb{H}}(x,z^\bot)$, and $x*z<0$ if and only if $x$ and $z$ lie on opposite sides of the hyperplane $z^\bot$;
\item[(iii)] if $\textrm{Span}(z,w)\cap \mathbb{H}^2\neq \emptyset$, then $|z*w|=\Vert z\Vert\cdot \Vert w\Vert\cosh d_{\mathbb{H}}(z^\bot,w^\bot)$, and $z*w<0$ if and only if $z$ and $w$ are oppositely oriented tangent vectors of $Q$, where $Q$ is the hyperbolic line Lorentz orthogonal to $z^\bot$ and $w^\bot$;
\item[(iv)] if $\textrm{Span}(z,w)\cap \mathbb{H}^2=\emptyset$, then $z*w=\Vert z\Vert \cdot\Vert w\Vert\cos \angle(z^\bot,w^\bot)$.
\end{description}
\end{proposition}

The following two propositions are useful for characterizing perpendicular geodesics in the Klein model.

\begin{proposition}(\cite{Ratcliffe}, Theorem 3.2.11)\label{Prop: Ratcliffe 2}
Let $P$ be a geodesic in $\mathbb{H}^2$, and let $y\in \mathbb{H}^2$ be a time-like vector.
There exists a unique geodesic in $\mathbb{H}^2$ passing through $y$ and Lorentz orthogonal to $P$.
\end{proposition}

\begin{proposition}(\cite{Thomas}, Proposition 4)\label{Prop: Thomas 1}
Let $x\in \mathbb{R}^3$ be a space-like vector, and let $\gamma=x^\bot\cap \mathbb{H}^2$ be the geodesic corresponding to $x^\bot$.
Suppose $\omega$ is another geodesic intersecting $\gamma$, and $P$ is the plane such that $\omega=P\cap \mathbb{H}^2$.
Then $\gamma$ and $\omega$ intersect at a right angle if and only if $x\in P$.
\end{proposition}

\subsection{A geometric interpretation}

\begin{figure}[!ht]
\centering
\includegraphics[scale=0.9]{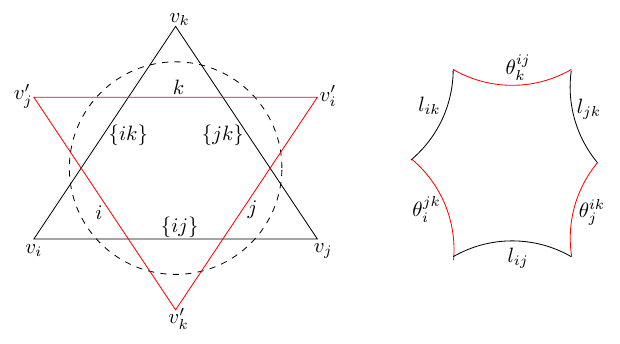}
\caption{A hyper-ideal hyperbolic triangle $v_iv_jv_k$ and its polar triangle $v^\prime_iv^\prime_jv^\prime_k$ in the Klein model (left), and a right-angled hyperbolic hexagon $\{ijk\}\in F$ in the Poincar\'{e} disk model (right).}
\label{Figure_2}
\end{figure}

Let $v_iv_jv_k$ be a hyper-ideal hyperbolic triangle in the Klein model, and let $\{r,s,t\}=\{i,j,k\}$.
By the notations in Remark \ref{Rmk: 5},
$E_{rs}=\textrm{Span}(v_r,v_s)\cap \mathbb{H}^2$ is a geodesic in $\mathbb{H}^2$.
Define $E^\prime_{st}=v_r^{\bot}\cap \mathbb{H}^2$. 
Then $E^\prime_{st}$ is also a geodesic in $\mathbb{H}^2$, and is Lorentz orthogonal to both $E_{rs}$ and $E_{rt}$.
The geodesic arc in $E^\prime_{st}$ bounded by $E_{rs}$ and $E_{rt}$ is a boundary arc, denoted by $r$.
An ideal edge $\{rs\}$ is the geodesic arc between two adjacent boundary arcs $r$ and $s$, lying on $E_{rs}$.
Furthermore, the six geodesic arcs form a right-angled hyperbolic hexagon $\{ijk\}$.
Conversely, for any geodesic $E_{rs}$ in $\mathbb{H}^2$, there exists a unique space-like vector $v^\prime_t$ such that $(v^\prime_t)^{\bot}\cap \mathbb{H}^2=E_{rs}$.
The three points $v^\prime_i, v^\prime_j,v^\prime_k$ are hyper-ideal and form a hyper-ideal hyperbolic triangle $v^\prime_iv^\prime_jv^\prime_k$, 
referred to as the polar triangle of $v_iv_jv_k$.
See Figure \ref{Figure_2} for illustration.

By the construction of the polar triangle $v^\prime_iv^\prime_jv^\prime_k$,
the vertices $v^\prime_i, v^\prime_j, v^\prime_k$ are determined by the vertices $v_i, v_j, v_k$, and vice versa.
Without loss of generality, we assume $\Vert v_i\Vert =\Vert v_j\Vert =\Vert v_k\Vert =1$.
Then
\begin{equation}\label{Eq: polar triangle}
v^\prime_i=\frac{v_k\otimes v_j}{\Vert v_k\otimes v_j\Vert },\
v^\prime_j=\frac{v_i\otimes v_k}{\Vert v_i\otimes v_k\Vert },\
v^\prime_k=\frac{v_j\otimes v_i}{\Vert v_j\otimes v_i\Vert }.
\end{equation}
Therefore, given a hyper-ideal hyperbolic triangle $v_iv_jv_k$, there exists a unique polar triangle $v^\prime_iv^\prime_jv^\prime_k$ whose vertices satisfy (\ref{Eq: polar triangle}).
Let $l_{rs}$ denote the length of the ideal edge $\{rs\}\in E$, and let $\theta^{rs}_{t}$ denote the length of the boundary arc $t\in B$.

\begin{lemma}\label{Lem: coplaner}
The vectors $c_{ijk}$, $c_{rs}$, and $v^\prime_t$ lie in the same 2-dimensional time-like vector subspace of the Lorentzian space $\mathbb{R}^3$,
where the notations $c_{ijk}$ and $c_{rs}$ are defined in Remark \ref{Rmk: GM} and Remark \ref{Rmk: 5}, respectively. 
\end{lemma}
\proof
Without loss of generality, we assume $r=i,\ s=j,\ t=k$.
Since $c_{ij}\in E_{ij}$, by Proposition \ref{Prop: Ratcliffe 2},
there exists a unique geodesic in $\mathbb{H}^2$ passing through $c_{ij}$ and Lorentz orthogonal to $E_{ij}$, which is $P_{ij}=\textrm{Span}(c_{ij},c_{ijk})\cap\mathbb{H}^2$.
Note that $E_{ij}=(v^\prime_k)^\bot\cap\mathbb{H}^2$.
By Proposition \ref{Prop: Thomas 1},
since $E_{ij}$ and $P_{ij}$ intersect at a right angle, it follows that $v^\prime_k\in \textrm{Span}(c_{ij},c_{ijk})$.
\qed

Combining Remark \ref{Rmk: GM} and Lemma \ref{Lem: coplaner}, we obtain
$c_{ijk}=\textrm{Span}(v^\prime_i,c_{jk})
\cap\textrm{Span}(v^\prime_j,c_{ik})
\cap\textrm{Span}(v^\prime_k,c_{ij})$.
Given the equivalence between $v_iv_jv_k$ and $v^\prime_iv^\prime_jv^\prime_k$,
it is natural to consider a dual model.
Note that for any $c_{ijk}\in \mathbb{R}^3$, the intersection $\textrm{Span}(v_r,c_{ijk})\cap v_r^\bot$ may not lie within $\mathbb{H}^2$.
In this paper, we focus on the case that
\begin{equation}\label{Eq: key}
\textrm{Span}(v_r,c_{ijk})\cap v_r^\bot\cap\mathbb{H}^2
=\textrm{Span}(v_r,c_{ijk})\cap E^\prime_{st}
\neq\emptyset,\ \forall r\in \{i,j,k\}.
\end{equation}
The condition (\ref{Eq: key}) constrains the position of $c_{ijk}$ when it is space-like, namely, $c_{ijk}$ can only lie within certain specific open domains. 
Let $c^\prime_{st}=\textrm{Span}(v_r,c_{ijk})\cap E^\prime_{st}$ denote the dual edge center,
which is the unique point on the geodesic $E^\prime_{st}$ such that the signed distance from $c^\prime_{st}$ to $\{rt\}$ is $\theta_{st}$, and the signed distance from $c^\prime_{st}$ to $\{rs\}$ is $\theta_{ts}$.
Thus
\begin{equation}\label{Eq: F20}
\theta_r^{st}=\theta_{st}+\theta_{ts}.
\end{equation}
Note that the signed distance $\theta_{st}$ is positive if $c^\prime_{st}$ and $\{rs\}$ lie on the same side of $\{rt\}$ along the hyperbolic geodesic $E^\prime_{st}$, and negative otherwise.
The above construction yields the following lemma.

\begin{lemma}\label{Lem: coplaner 2}
The vectors $c_{ijk}$, $c^\prime_{rs}$, and $v_t$ lie in the same 2-dimensional time-like vector subspace of $\mathbb{R}^3$.
\end{lemma}

By Proposition \ref{Prop: Ratcliffe 2},
let $P^\prime_{ij}$ denote the edge perpendicular passing through $c^\prime_{ij}$ and Lorentz orthogonal to $E^\prime_{ij}$.
Then $P^\prime_{ij}=\textrm{Span}(v_k,c_{ijk})\cap\mathbb{H}^2$.
Consequently, the three edge perpendiculars $P^\prime_{ij},P^\prime_{jk},P^\prime_{ki}$ of $\{ijk\}$ intersect in a common point $c_{ijk}$.

\begin{figure}[!ht]
\centering
\includegraphics[scale=1]{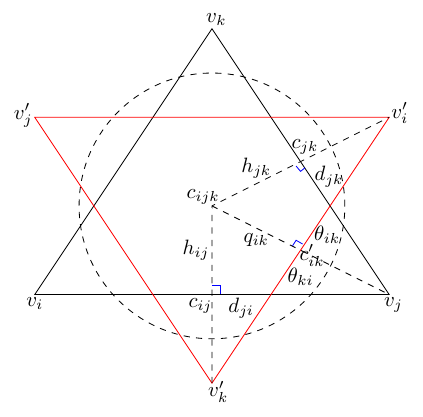}
\caption{Data for two hyper-ideal hyperbolic triangles in the Klein model}
\label{Figure_3}
\end{figure}

Recall the definition of the signed distance $h_{rs}$ from the face center $c_{ijk}$ to the hyperbolic geodesic $E_{rs}$.
It is positive if $c_{ijk}$ lies on the same side of $E_{rs}$ as the right-angled hyperbolic hexagon $\{ijk\}$, negative otherwise, and zero if $c_{ijk}$ lies on $E_{rs}$.
Similarly, let $q_{rs}$ denote the signed distance from $c_{ijk}$ to $E^\prime_{rs}$,
which is defined to be positive if $c_{ijk}$ lies on the same side of the geodesic $E^\prime_{rs}$ as the polar triangle $v^\prime_iv^\prime_jv^\prime_k$, negative otherwise, and zero if $c_{ijk}$ lies on $E^\prime_{rs}$.
Please refer to Figure \ref{Figure_3}.

\begin{figure}[!ht]
\centering
\includegraphics[scale=1]{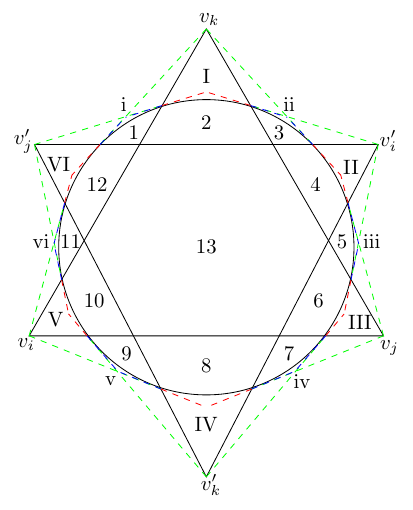}
\caption{Partition of the right-angled hyperbolic hexagon $\{ijk\}$ in the Klein model. If $c_{ijk}$ is space-like, it may lie in the open domains bounded by two red tangents or two blue tangents.}
\label{Figure_4}
\end{figure}

In the Klein model, the hyper-ideal hyperbolic triangle $v_iv_jv_k$ and its polar triangle $v^\prime_iv^\prime_jv^\prime_k$ partition the unit disk into 13 open domains,
denoted by $D_1,D_2,...,D_{13}$ respectively.
Please refer to Figure \ref{Figure_4}.
By the definitions of $h$ and $q$, one can directly obtain the relationships between the position of $c_{ijk}$ and the signs of $h$ and $q$ in these 13 open domains, as shown in Table \ref{table}.
For example, if $c_{ijk}\in D_1$, then $c_{ijk}$ lies on the opposite side of $E_{ik}$ as $v_iv_jv_k$, so $h_{ik}<0$.
Moreover, $c_{ijk}$ lies on the same side of $E_{ij}$ and $ E_{jk}$ as $v_iv_jv_k$, which implies $h_{ij}>0$ and $h_{jk}>0$.
Similarly, $c_{ijk}$ lies on the opposite side of $E^\prime_{ij}$, and the same side of $E^\prime_{ik}$ and $E^\prime_{jk}$ as $v^\prime_iv^\prime_jv^\prime_k$, which implies $q_{ij}<0$, $q_{jk}>0$, and $q_{ik}>0$.
If $c_{ijk}$ is space-like, it can lie only in one of 12 open domains, each bounded by two tangent lines and the unit circle.
These domains are denoted by $D_\mathrm{I}, D_\mathrm{II}, D_\mathrm{III}, D_\mathrm{IV}, D_\mathrm{V}, D_\mathrm{VI}$ and $D_\mathrm{i}, D_\mathrm{ii}, D_\mathrm{iii}, D_\mathrm{iv}, D_\mathrm{v}, D_\mathrm{vi}$ respectively.
Specifically, the domain $D_\mathrm{i}$ is enclosed by the tangents from vertices $v_k$ and $v_j^\prime$ to the unit circle, together with the unit circle itself.
The domain $D_\mathrm{I}$ is enclosed by the tangents from vertices $v_i^\prime$ and $v_j^\prime$ to the unit circle, together with the unit circle itself.
The remaining domains are defined analogously. 
Please refer to Figure \ref{Figure_4}.
If $c_{ijk}$ does not lie in any of these 12 open domains, then either $\textrm{Span}(c_{ijk}, v_r)\cap \mathbb{H}^2=\emptyset$
or $\textrm{Span}(c_{ijk}, v^\prime_r)\cap \mathbb{H}^2=\emptyset$,
which contradicts the assumption that the edge centers $c_{st}\in E_{rs}$ and $c^\prime_{st}\in E^\prime_{rs}$ are in $\mathbb{H}^2$ for $\{r,s,t\}= \{i,j,k\}$.
Similar to the case that $c_{ijk}$ is time-like,
one can also directly obtain the relationships between the position of $c_{ijk}$ and the signs of $h$ and $q$ in these 12 open domains, as shown in Table \ref{table}.
By continuously moving $c_{ijk}$ to the light cone, analogous results are obtained.
Note that if $c_{ijk}$ is light-like,
then it can lie only on the unit disk with 12 points removed.

\begin{table}[ht!]
\resizebox{\textwidth}{!}{
\begin{tabular}{cccccccccccccc}
\hline
Domains & $D_1$($D_\mathrm{i}$) & $D_2$($D_\mathrm{I}$) & $D_3$($D_\mathrm{ii}$) & $D_4$($D_\mathrm{II}$) & $D_5$($D_\mathrm{iii}$) & $D_6$($D_\mathrm{III}$) & $D_7$($D_\mathrm{iv}$) & $D_8$($D_\mathrm{IV}$) & $D_9$($D_\mathrm{v}$) & $D_{10}$($D_\mathrm{V}$) & $D_{11}$($D_\mathrm{vi}$) & $D_{12}$($D_\mathrm{VI}$) & $D_{13}$ \\ \hline
$h_{jk}$ & + & + & - & - & - & + & + & + & + & + & + & + & +   \\
$h_{ik}$ & - & + & + & + & + & + & + & + & + & + & - & - & +   \\
$h_{ij}$ & + & + & + & + & + & + & - & - & - & + & + & + & +   \\
$q_{jk}$ & + & + & + & + & + & + & + & + & - & - & - & + & +   \\
$q_{ik}$ & + & + & + & + & - & - & - & + & + & + & + & + & +   \\
$q_{ij}$ & - & - & - & + & + & + & + & + & + & + & + & + & +   \\
\hline
\end{tabular}}
\caption{Relationships between the position of $c_{ijk}$ and the signs of $h$ and $q$}
\label{table}
\end{table}

In the following, we present an alternative derivation of Table \ref{table}.
Suppose $c_{ijk}$ is time-like.
Without loss of generality, we assume $c_{ijk}\in D_1$.
By Lemma \ref{Lem: coplaner}, the vectors $c_{ijk}$, $c_{ik}$, and $v^\prime_j$ lie in the same 2-dimensional vector subspace of $\mathbb{R}^3$.
Since $c_{ijk}\notin v^\prime_iv^\prime_jv^\prime_k$,
then $c_{ik}\notin \{ik\}$.
Furthermore, since $c_{ik}$ lies on the opposite side of $k$ as $i$.
by the definition of $d$ in Remark \ref{Rmk: 5},
we have $d_{ki}<0$ and $d_{ik}>0$.
Similarly, $d_{kj}<0$, $d_{jk}>0$, $d_{ij}>0$, and $d_{ji}>0$.
By Lemma \ref{Lem: coplaner 2},
the vectors $c_{ijk}$, $c^\prime_{ij}$, and $v_k$ lie in the same 2-dimensional vector subspace of $\mathbb{R}^3$.
Since $c_{ijk}\notin v_iv_jv_k$, it follows that $c^\prime_{ij}\notin k$.
Furthermore, since $c^\prime_{ij}$ is on the opposite side of $\{ik\}$ as $\{jk\}$, it follows that $\theta_{ji}<0$ and $\theta_{ij}>0$.
Similarly, $\theta_{jk}<0$, $\theta_{kj}>0$, $\theta_{ik}>0$, and $\theta_{ki}>0$.
By Corollary \ref{Cor: time-like}, we have $h_{ik}<0$, $h_{jk}>0$, $h_{ij}>0$, $q_{ij}<0$, $q_{ik}>0$, and $q_{jk}>0$.
The case that $c_{ijk}$ is space-like is established via the two lemmas below.

\begin{lemma}\label{Lem: same signs 1}
The signs of $q$ and $h$ in the domains $D_\mathrm{I}, D_\mathrm{II}, D_\mathrm{III}, D_\mathrm{IV}, D_\mathrm{V}, D_\mathrm{VI}$ coincide with those in the domains $D_2,D_4,D_6,D_8,D_{10},D_{12}$, respectively.
\end{lemma}
\proof
We only prove the signs of $q$ and $h$ in the domain $D_{\mathrm{V}}$ coincide with those in the domain $D_{10}$, and the remaining cases follow analogously.

If $c_{ijk}\in D_{10}$, then $q_{jk}<0$ and the rest are positive by Table \ref{table}.
Suppose $c_{ijk}\in D_{\mathrm{V}}$.
By Lemma \ref{Lem: coplaner}, the vectors $c_{ijk}$, $c_{ij}$, and $v^\prime_k$ lie in the same 2-dimensional vector subspace of $\mathbb{R}^3$.
Since $c_{ijk}\notin v^\prime_iv^\prime_jv^\prime_k$,
it follows that $c_{ij}\notin \{ij\}$.
Furthermore, since $c_{ij}$ lies on the opposite side of $j$ as $i$, it follows that $d_{ij}<0$ and $d_{ji}>0$.
Similarly, $d_{ik}<0$, $d_{ki}>0$, $d_{jk}>0$, and $d_{kj}>0$.
By Proposition \ref{Prop: space-like 1},
we have $h_{ij}>0$, $h_{jk}>0$, and $h_{ik}>0$.
Note that $c_{ijk}\in v_iv_jv_k$.
By Lemma \ref{Lem: coplaner 2},
the vectors $c_{ijk}$, $c^\prime_{rs}$, and $v_t$ lie in the same 2-dimensional vector subspace of $\mathbb{R}^3$ for $\{r,s,t\}=\{i,j,k\}$.
Consequently, $c^\prime_{rs}\in t$.
Hence, $\theta_{ij}>0$, $\theta_{ji}>0$, $\theta_{jk}<0$, $\theta_{kj}>0$, $\theta_{ik}>0$, and $\theta_{ki}>0$.
By Proposition \ref{Prop: space-like 1},
we have $q_{ij}>0$, $q_{jk}<0$, and $q_{ik}>0$.
\qed

\begin{lemma}\label{Lem: same signs 2}
The signs of $q$ and $h$ in the domains $D_\mathrm{i}, D_\mathrm{ii}, D_\mathrm{iii}, D_\mathrm{iv}, D_\mathrm{v}, D_\mathrm{vi}$ coincide with those in the domains $D_1,D_3,D_5,D_7,D_{9},D_{11}$, respectively.
\end{lemma}
\proof
We only prove the signs of $q$ and $h$ in the domain $D_{\mathrm{vi}}$ coincide with those in the domain $D_{11}$, and the remaining cases follow analogously.

If $c_{ijk}\in D_{11}$, then $q_{jk}<0$, $h_{ik}<0$, and the rest are positive by Table \ref{table}.
Suppose $c_{ijk}\in D_{\mathrm{vi}}$.
By Lemma \ref{Lem: coplaner}, the vectors $c_{ijk}$, $c_{ij}$, and $v^\prime_k$ lie in the same 2-dimensional vector subspace of $\mathbb{R}^3$.
Since $c_{ijk}\notin v^\prime_iv^\prime_jv^\prime_k$,
it follows that $c_{ij}\notin \{ij\}$.
Furthermore, since $c_{ij}$ lies on the opposite side of $j$ as $i$, it follows that $d_{ij}<0$ and $d_{ji}>0$.
Similarly, $d_{ik}<0$, $d_{ki}>0$, $d_{jk}>0$, and $d_{kj}>0$.
By Proposition \ref{Prop: space-like 2},
we have $h_{ij}<0$, $h_{jk}>0$, and $h_{ik}>0$.
By Lemma \ref{Lem: coplaner 2},
the vectors $c_{ijk}$, $c^\prime_{jk}$, and $v_i$ lie in the same 2-dimensional vector subspace of $\mathbb{R}^3$.
Since $c_{ijk}\notin v_iv_jv_k$, it follows that $c^\prime_{jk}\notin i$.
Furthermore, since $c^\prime_{jk}$ lies on the opposite side of $\{ik\}$ as $\{ij\}$, it follows that $\theta_{jk}<0$ and $\theta_{kj}>0$.
Similarly, $\theta_{ji}<0$, $\theta_{ij}>0$, $\theta_{ik}>0$, and $\theta_{ki}>0$.
By Proposition \ref{Prop: space-like 2}, we have
$q_{ij}>0$, $q_{jk}<0$, and $q_{ik}>0$.
\qed

By direct calculations, we obtain the following results.
As the proofs of Lemma \ref{Lem: time-like}, Proposition \ref{Prop: space-like 1} and Proposition \ref{Prop: space-like 2} are too long, we defer them to Appendix \ref{appendix a}.
For simplicity, we set $h_t=h_{rs}$ and $q_t=q_{rs}$ for $\{r,s,t\}=\{i,j,k\}$.

\begin{lemma}\label{Lem: time-like}
If $c_{ijk}$ is time-like, then
\begin{align*}
&\sinh q_r=\cosh h_s\cdot\sinh d_{rt}=\cosh h_t\cdot\sinh d_{rs},\\
&\sinh h_r=\cosh q_s\cdot\sinh \theta_{rt}=\cosh q_t\cdot\sinh \theta_{rs}.
\end{align*}
\end{lemma}

As a direct corollary of Lemma \ref{Lem: time-like}, we have the following result.
\begin{corollary}\label{Cor: time-like}
If $c_{ijk}$ is time-like,
then the signs of $q_r$, $d_{rt}$, and $d_{rs}$ coincide, and the signs of $h_r$, $\theta_{rs}$ and $\theta_{rt}$ coincide.
\end{corollary}

\begin{proposition}\label{Prop: space-like 1}
If $c_{ijk}\in \mathbb{R}^3$ lies in the open domain $D_\mathrm{V}$, then
\begin{equation*}
\begin{aligned}
&\cosh q_i=-\sinh h_j\cdot\sinh d_{ik}=-\sinh h_k\cdot\sinh d_{ij},\\
&\cosh q_j=\sinh h_i\cdot\sinh d_{jk}=\sinh h_k\cdot\sinh d_{ji},\\
&\cosh q_k=\sinh h_i\cdot\sinh d_{kj}=\sinh h_j\cdot\sinh d_{ki},\\
&\cosh h_i=\sinh q_j\cdot\sinh \theta_{ik}=\sinh q_k\cdot\sinh \theta_{ij},\\
&\cosh h_j=-\sinh q_i\cdot\sinh \theta_{jk}=\sinh q_k\cdot\sinh \theta_{ji},\\
&\cosh h_k=-\sinh q_i\cdot\sinh \theta_{kj}=\sinh q_j\cdot\sinh \theta_{ki}.
\end{aligned}
\end{equation*}
If $c_{ijk}\in \mathbb{R}^3$ lies in the open domain $D_\mathrm{II}$, then we swap $h$ with $q$, and $\theta$ with $d$.
\end{proposition}

If $c_{ijk}$ lies in the domains $D_\mathrm{I}, D_\mathrm{III}, D_\mathrm{IV}$ or $D_\mathrm{VI}$, analogous propositions hold.
These propositions enable us to prove Lemma \ref{Lem: same signs 1} and further derive Table \ref{table}.
Moreover, Table \ref{table} shows that $c_{ijk}$ lies in the open domains  $D_\mathrm{I}, D_\mathrm{II}, D_\mathrm{III}, D_\mathrm{IV}, D_\mathrm{V}$ or $D_\mathrm{VI}$ if and only if one of $q_i,q_j,q_k,h_i,h_j,h_k$ is negative and the rest are positive.
Consequently, these propositions can be written in a unified form.

\begin{lemma}\label{Lem: space-like 1}
Let $c_{ijk}\in \mathbb{R}^3$ lie in one of the open domains $D_\mathrm{I}, D_\mathrm{II}, D_\mathrm{III}, D_\mathrm{IV}, D_\mathrm{V}, D_\mathrm{VI}$.
If $q_r<0$, then
\begin{equation*}
\begin{aligned}
&\cosh q_r=-\sinh h_s\cdot\sinh d_{rt}
=-\sinh h_t\cdot\sinh d_{rs},\\
&\cosh q_s=\sinh h_r\cdot\sinh d_{st}
=\sinh h_t\cdot\sinh d_{sr},\\
&\cosh q_t=\sinh h_r\cdot\sinh d_{ts}
=\sinh h_s\cdot\sinh d_{tr},\\
&\cosh h_r=\sinh q_s\cdot\sinh \theta_{rt}
=\sinh q_t\cdot\sinh \theta_{rs},\\
&\cosh h_s=-\sinh q_r\cdot\sinh \theta_{st}
=\sinh q_t\cdot\sinh \theta_{sr},\\
&\cosh h_t=-\sinh q_r\cdot\sinh \theta_{ts}
=\sinh q_s\cdot\sinh \theta_{tr}.
\end{aligned}
\end{equation*}
If $h_r<0$, then we swap $h$ with $q$, and $\theta$ with $d$.
\end{lemma}

\begin{proposition}\label{Prop: space-like 2}
Let $c_{ijk}\in \mathbb{R}^3$ lie in the open domain $D_\mathrm{vi}$, then
\begin{equation*}
\begin{aligned}
&\cosh q_i=\sinh h_j\cdot\sinh d_{ik}
=-\sinh h_k\cdot\sinh d_{ij},\\
&\cosh q_j=\sinh h_i\cdot\sinh d_{jk}
=\sinh h_k\cdot\sinh d_{ji},\\
&\cosh q_k=\sinh h_i\cdot\sinh d_{kj}
=-\sinh h_j\cdot\sinh d_{ki},\\
&\cosh h_i=\sinh q_j\cdot\sinh \theta_{ik}
=\sinh q_k\cdot\sinh \theta_{ij},\\
&\cosh h_j=\sinh q_i\cdot\sinh \theta_{jk}
=-\sinh q_k\cdot\sinh \theta_{ji},\\
&\cosh h_k=-\sinh q_i\cdot\sinh \theta_{kj}
=\sinh q_j\cdot\sinh \theta_{ki}.
\end{aligned}
\end{equation*}
\end{proposition}

Similarly, if $c_{ijk}$ lies in the domains $D_\mathrm{i}, D_\mathrm{ii}, D_\mathrm{iii}, D_\mathrm{iv}$ or $D_\mathrm{v}$, 
corresponding propositions hold.
These propositions allow us to prove Lemma \ref{Lem: same signs 2} and further derive Table \ref{table}. 
Moreover, Table \ref{table} shows that $c_{ijk}$ lies in the open domains $D_\mathrm{i}, D_\mathrm{ii}, D_\mathrm{iii}, D_\mathrm{iv}, D_\mathrm{v}$ or $D_\mathrm{vi}$ if and only if $q_r<0$, $h_s<0$, and the rest are positive.
Consequently, these propositions also admit a unified form.

\begin{lemma}\label{Lem: space-like 2}
Let $c_{ijk}\in \mathbb{R}^3$ lie in one of the open domains $D_\mathrm{i}, D_\mathrm{ii}, D_\mathrm{iii}, D_\mathrm{iv}, D_\mathrm{v}, D_\mathrm{vi}$.
If $q_r<0$ and $h_s<0$, then
\begin{equation*}
\begin{aligned}
&\cosh q_r=\sinh h_s\cdot\sinh d_{rt}=-\sinh h_t\cdot\sinh d_{rs},\\
&\cosh q_s=\sinh h_r\cdot\sinh d_{st}=\sinh h_t\cdot\sinh d_{sr},\\
&\cosh q_t=\sinh h_r\cdot\sinh d_{ts}=-\sinh h_s\cdot\sinh d_{tr},\\
&\cosh h_r=\sinh q_s\cdot\sinh \theta_{rt}=\sinh q_t\cdot\sinh \theta_{rs},\\
&\cosh h_s=\sinh q_r\cdot\sinh \theta_{st}=-\sinh q_t\cdot\sinh \theta_{sr},\\
&\cosh h_t=-\sinh q_r\cdot\sinh \theta_{ts}=\sinh q_s\cdot\sinh \theta_{tr}.
\end{aligned}
\end{equation*}
\end{lemma}

As a direct corollary of Lemma \ref{Lem: space-like 1} and Lemma \ref{Lem: space-like 2}, we have the following result.

\begin{corollary}\label{Cor: space-like}
If $c_{ijk}$ is space-like,
then the signs of $q_r$, $d_{rt}$, and $d_{rs}$ coincide, and the signs of $h_r$, $\theta_{rs}$, and $\theta_{rt}$ coincide.
\end{corollary}

Combining Corollary \ref{Cor: time-like} and Corollary \ref{Cor: space-like} gives the following lemma.
\begin{lemma}\label{Lem: time-space-like}
For any right-angled hyperbolic  hexagon $\{ijk\}\in F$,
the signs of $q_r$, $d_{rt}$, and $d_{rs}$ coincide, and the signs of $h_r$, $\theta_{rs}$, and $\theta_{rt}$ coincide.
\end{lemma}
\proof
If $c_{ijk}$ is time-like, the conclusion follows from Corollary \ref{Cor: time-like}.
If $c_{ijk}$ is space-like, the conclusion is derived from Corollary \ref{Cor: space-like}.
If $c_{ijk}$ is light-like, the conclusion holds by continuity upon moving $c_{ijk}$ to the light cone.  
\qed

\subsection{Variational formulas of generalized angle}

We decompose Theorem \ref{Thm: angle variations} into Lemma \ref{Lem: angle variations 1}, Lemma \ref{Lem: angle variations 2}, and Lemma \ref{Lem: angle variations 3}, and prove them respectively.

\begin{lemma}\label{Lem: angle variations 1}
For any right-angled hyperbolic hexagon $\{ijk\}\in F$, we have
\begin{equation}\label{Eq: angle 1}
\frac{\partial \theta^{jk}_i}{\partial f_j}
=\frac{-1}{\sinh d_{ji}}\frac{\tanh^\beta h_{ij}}{\sinh l_{ij}},
\end{equation}
where $\beta=1$ if $c_{ijk}$ is time-like, and $\beta=-1$ if $c_{ijk}$ is space-like.
If $c_{ijk}$ is light-like, the formula is interpreted as $\tanh^\beta h_{ij}=\tanh^\beta h_{jk}=\tanh^\beta h_{ki}=1$.
\end{lemma}
\proof
Without loss of generality, we assume that $\Vert v_i\Vert =\Vert v_j\Vert =\Vert v_k\Vert =1$ and $\Vert c_{ijk}\Vert ^2=\pm1$.
For simplicity, we set $l_r=l_{st},\ \theta_r=\theta^{st}_r,\ h_r=h_{st},\ q_r=q_{st}$, and $A=\sinh l_r\sinh l_s\sinh \theta_t$, where $\{r,s,t\}=\{i,j,k\}$.
By generalized hyperbolic cosine law, i.e.,
$\cosh l_r
=\cosh \theta_r\sinh l_s\sinh l_t-\cosh l_s\cosh l_t$,
we derive
\begin{equation}\label{Eq: F24}
\frac{\partial\theta_i}{\partial l_i}=\frac{\sinh l_i}{A},\
\frac{\partial\theta_i}{\partial l_j}=-\frac{\sinh l_i\cosh \theta_k}{A},\
\frac{\partial\theta_i}{\partial l_k}=-\frac{\sinh l_i\cosh \theta_j}{A}.
\end{equation}
Using \eqref{Eq: variation 1} and \eqref{Eq: variation 2}, we obtain
\begin{equation}\label{Eq: F6}
\begin{aligned}
\frac{\partial \theta_i}{\partial f_j}
&=\frac{\partial \theta_i}{\partial l_i}\frac{\partial l_i}{\partial f_j}
+\frac{\partial \theta_i}{\partial l_j}\frac{\partial l_j}{\partial f_j}
+\frac{\partial \theta_i}{\partial l_k}\frac{\partial l_k}{\partial f_j}\\
&=\frac{\sinh l_i}{A}\cdot\coth d_{jk}
+\frac{-\sinh l_i\cosh \theta_j}{A}\cdot\coth d_{ji}\\
&=\frac{1}{\sinh l_k\sinh \theta_j}\cdot\left(\coth d_{jk}-\cosh \theta_j\cdot\coth d_{ji}\right)\\
&=\frac{1}{\sinh l_k\sinh \theta_j}\cdot\left[\coth d_{jk}-\cosh (\theta_{ki}+\theta_{ik})\cdot\coth d_{ji}\right]\\
&=\frac{1}{\sinh l_k\sinh \theta_j}\cdot(\coth d_{jk}-\cosh \theta_{ki}\cosh\theta_{ik}\coth d_{ji}\\
&\ \ \ \ \ \ \ \ \ \ \ \ \ \ \ \ \ \ \ \ \ \ \ \ \ \ \ \ \ \ \ \ \ -\sinh \theta_{ki}\sinh\theta_{ik}\coth d_{ji}),
\end{aligned}
\end{equation}
where (\ref{Eq: F20}) is used in the fourth line.
Define
\begin{equation}\label{Eq: F21}
L_1=\coth d_{jk}-\cosh \theta_{ki}\cosh\theta_{ik}\coth d_{ji},\ L_2=-\sinh \theta_{ki}\sinh\theta_{ik}\coth d_{ji}.
\end{equation}
It follows from (\ref{Eq: F6}) that
\begin{equation}\label{Eq: F22}
\frac{\partial \theta_i}{\partial f_j}=\frac{1}{\sinh l_k\sinh \theta_j}\cdot(L_1+L_2).
\end{equation}

To compute $L_1$ and $L_2$, some preliminary results are required. 
We begin by referring to Figure \ref{Figure_3}.
By Proposition \ref{Prop: Ratcliffe 1} (ii), we have 
$v_j\ast c_{ij}=+\sinh d_\mathbb{H}(v_j^\bot,c_{ij})=+\sinh(-d_{ji})=-\sinh d_{ji}$,
or $v_j\ast c_{ij}=-\sinh d_\mathbb{H}(v_j^\bot,c_{ij})=-\sinh(+d_{ji})=-\sinh d_{ji}$.
In both cases, it follows that $v_j\ast c_{ij}=-\sinh d_{ji}$.
Applying Proposition \ref{Prop: Right angle} to the generalized right-angled hyperbolic triangle $c_{ijk}c_{ij}v_j$ yields
\begin{equation}\label{Eq: F7}
v_j\ast c_{ijk}=-(c_{ij}\ast c_{ijk})(v_j\ast c_{ij})
=\sinh d_{ji}\cdot(c_{ij}\ast c_{ijk}).
\end{equation}
Similarly, in the generalized right-angled hyperbolic  triangle $c_{ijk}c_{jk}v_j$, we derive
\begin{equation}\label{Eq: F10}
v_j\ast c_{ijk}=-(v_j\ast c_{jk})(c_{jk}\ast c_{ijk})=
\sinh d_{jk}\cdot(c_{jk}\ast c_{ijk}),
\end{equation}
where $v_j\ast c_{jk}=-\sinh d_{jk}$ follows from Proposition \ref{Prop: Ratcliffe 1} (ii).
Combining (\ref{Eq: F7}) and (\ref{Eq: F10}) gives
\begin{equation}\label{Eq: F14}
\sinh d_{ji}\cdot(c_{ij}\ast c_{ijk})
=\sinh d_{jk}\cdot(c_{jk}\ast c_{ijk}).
\end{equation}
By Proposition \ref{Prop: Ratcliffe 1} (iii), we have
$(c_{ijk}\otimes v_j)\ast(v_j\otimes c_{jk})
=\Vert v_j\otimes c_{ijk}\Vert \cdot\Vert v_j\otimes c_{jk}\Vert (-\cosh \theta_{ik}).$
Moreover, direct calculation yields
\begin{equation*}
(c_{ijk}\otimes v_j)\ast(v_j\otimes c_{jk})
=(c_{jk}\ast c_{ijk})\Vert v_j\Vert ^2-(v_j\ast c_{ijk})(v_j\ast c_{jk})
=(c_{jk}\ast c_{ijk})\cdot \cosh^2 d_{jk},
\end{equation*}
where (\ref{Eq: F10}) is used in the second equality.
Since $\Vert v_j\otimes c_{jk}\Vert ^2=(v_j\ast c_{jk})^2+1=\cosh^2 d_{jk}$,
it follows that
\begin{equation}\label{Eq: F11}
\cosh \theta_{ik}=-\frac{(c_{ijk}\otimes v_j)\ast(v_j\otimes c_{jk})}{\Vert v_j\otimes c_{ijk}\Vert \cdot\Vert v_j\otimes c_{jk}\Vert }
=-\frac{(c_{jk}\ast c_{ijk})\cdot \cosh d_{jk}}{\Vert v_j\otimes c_{ijk}\Vert }.
\end{equation}
Similarly,
$(c_{ijk}\otimes v_j)\ast(v_j\otimes c_{ij})
=\Vert v_j\otimes c_{ijk}\Vert \cdot\Vert v_j\otimes c_{ij}\Vert (-\cosh \theta_{ki}).$
Moreover,
\begin{equation*}
(c_{ijk}\otimes v_j)\ast(v_j\otimes c_{ij})
=(c_{ij}\ast c_{ijk})\Vert v_j\Vert ^2-(v_j\ast c_{ijk})(v_j\ast c_{ij})
=(c_{ij}\ast c_{ijk})\cdot \cosh^2 d_{ji},
\end{equation*}
where (\ref{Eq: F7}) is used in the second equality.
Since $\Vert v_j\otimes c_{ij}\Vert ^2=(v_j\ast c_{ij})^2+1=\cosh^2 d_{ji}$, it follows that
\begin{equation}\label{Eq: F8}
\cosh \theta_{ki}=-\frac{(c_{ijk}\otimes v_j)\ast(v_j\otimes c_{ij})}{\Vert v_j\otimes c_{ijk}\Vert \cdot\Vert v_j\otimes c_{ij}\Vert }
=-\frac{(c_{ij}\ast c_{ijk})\cdot \cosh d_{ji}}{\Vert v_j\otimes c_{ijk}\Vert }.
\end{equation}
By the properties of Lorentzian cross product, we have
\begin{equation}\label{Eq: F15}
\Vert v_j\otimes c_{ijk}\Vert ^2=(v_j\otimes c_{ijk})*(v_j\otimes c_{ijk})
%=(v_j\ast c_{ijk})^2-\Vert c_{ijk}\Vert ^2\Vert v_j\Vert ^2
=(v_j\ast c_{ijk})^2-\Vert c_{ijk}\Vert ^2,
\end{equation}
and
\begin{equation}\label{Eq: F16}
\Vert c_{ij}\otimes c_{ijk}\Vert ^2=(c_{ij}\otimes c_{ijk})*(c_{ij}\otimes c_{ijk})
%=(c_{ij}\ast c_{ijk})^2-\Vert c_{ijk}\Vert ^2\Vert c_{ij}\Vert ^2
=(c_{ij}\ast c_{ijk})^2+\Vert c_{ijk}\Vert ^2.
\end{equation}
Adding (\ref{Eq: F15}) and (\ref{Eq: F16}) gives
\begin{equation}\label{Eq: F17}
\begin{aligned}
\Vert v_j\otimes c_{ijk}\Vert ^2+\Vert c_{ij}\otimes c_{ijk}\Vert ^2
&=(v_j\ast c_{ijk})^2+(c_{ij}\ast c_{ijk})^2\\
&=(c_{ij}\ast c_{ijk})^2\cdot\sinh^2 d_{ji}+(c_{ij}\ast c_{ijk})^2\\
&=(c_{ij}\ast c_{ijk})^2\cosh^2 d_{ji},
\end{aligned}
\end{equation}
where (\ref{Eq: F7}) is used in the second line.
Then
\begin{equation}\label{Eq: F9}
\begin{aligned}
\sinh^2 \theta_{ki}
&=\cosh^2 \theta_{ki}-1\\
&=\frac{1}{\Vert v_j\otimes c_{ijk}\Vert ^2}[(c_{ij}\ast c_{ijk})^2\cosh^2 d_{ji}-\Vert v_j\otimes c_{ijk}\Vert ^2]\\
&=\frac{\Vert c_{ij}\otimes c_{ijk}\Vert ^2}{\Vert v_j\otimes c_{ijk}\Vert ^2},
\end{aligned}
\end{equation}
where (\ref{Eq: F8}) is used in the second line and (\ref{Eq: F17}) is used in the last line.

We are now ready to compute $L_1$ and $L_2$ in (\ref{Eq: F21}). 
For $L_2$, we proceed as follows:
\begin{align*}
L_2&=-\sinh \theta_{ki}\sinh\theta_{ik}\cdot\coth d_{ji}\\
&=-\sinh \theta_{ki}\cdot\sinh\theta_{ik}\cdot \frac{\cosh d_{ji}}{\sinh d_{ji}}\cdot\frac{\Vert v_j\otimes c_{ijk}\Vert }{c_{ij}\ast c_{ijk}}\cdot\frac{c_{ij}\ast c_{ijk}}{\Vert v_j\otimes c_{ijk}\Vert }\\
&=\sinh \theta_{ki} \sinh \theta_{ik}\cdot \frac{\Vert v_j\otimes c_{ijk}\Vert }{(c_{ij}\ast c_{ijk})\sinh d_{ji}}\cdot\cosh \theta_{ki},
\end{align*}
where (\ref{Eq: F8}) is used in the last line.
For $L_1$, we have
\begin{align*}
L_1&=\coth d_{jk}-\cosh \theta_{ki}\cosh\theta_{ik}\cdot\coth d_{ji}\\
&=\frac{\cosh d_{jk}\cdot (c_{jk}\ast c_{ijk})}{\sinh d_{ji}\cdot (c_{ij}\ast c_{ijk})}
+\frac{(c_{ij}\ast c_{ijk})\cosh d_{ji}}{\Vert v_j\otimes c_{ijk}\Vert }\cdot\cosh\theta_{ik}\cdot \frac{\cosh d_{ji}}{\sinh d_{ji}}\\
&=\frac{\cosh\theta_{ik}}{\sinh d_{ji}\cdot (c_{ij}\ast c_{ijk})\cdot \Vert v_j\otimes c_{ijk}\Vert }\big(-\Vert v_j\otimes c_{ijk}\Vert ^2+(c_{ij}\ast c_{ijk})^2\cosh^2 d_{ji}\big)\\
&=\frac{\cosh\theta_{ik}\cdot\Vert c_{ij}\otimes c_{ijk}\Vert ^2}{\sinh d_{ji}\cdot (c_{ij}\ast c_{ijk})\cdot \Vert v_j\otimes c_{ijk}\Vert }\\
&=\frac{\cosh\theta_{ik}\sinh^2 \theta_{ki}\Vert v_j\otimes c_{ijk}\Vert }{\sinh d_{ji}(c_{ij}\ast c_{ijk})},
\end{align*}
where (\ref{Eq: F14}) and (\ref{Eq: F8}) are used in the second line, (\ref{Eq: F11}) is used in the third line, (\ref{Eq: F17}) is used in the fourth line, and (\ref{Eq: F9}) is used in the last line.
Substituting $L_1$ and $L_2$ into (\ref{Eq: F22}) yields
\begin{equation*}
\begin{aligned}
\frac{\partial \theta_i}{\partial f_j}
&=\frac{1}{\sinh l_k\sinh \theta_j}\cdot(L_1+L_2)\\
&=\frac{1}{\sinh l_k\sinh \theta_j}\cdot\frac{\sinh \theta_{ki}\cdot \Vert v_j\otimes c_{ijk}\Vert }{(c_{ij}\ast c_{ijk})\sinh d_{ji}}\cdot (\cosh \theta_{ik}\sinh \theta_{ki}+\sinh \theta_{ik}\cosh \theta_{ki})\\
&=\frac{1}{\sinh l_k\sinh \theta_j}\cdot\frac{\sinh \theta_{ki}\cdot \Vert v_j\otimes c_{ijk}\Vert }{(c_{ij}\ast c_{ijk})\sinh d_{ji}}\cdot \sinh \theta_j\\
&=\frac{1}{\sinh l_k\sinh d_{ji}}\cdot\frac{\sinh \theta_{ki}\cdot \Vert v_j\otimes c_{ijk}\Vert }{c_{ij}\ast c_{ijk}}.
\end{aligned}
\end{equation*}

\noindent\textbf{(i):}
Assume that $c_{ijk}$ is time-like.
Then $\Vert c_{ijk}\Vert ^2=-1$ and $v_j*c_{ijk}=-\sinh q_j$.
By (\ref{Eq: F15}), we have
$\Vert v_j\otimes c_{ijk}\Vert =\cosh q_j$.
Furthermore, by Lemma \ref{Lem: time-like}, we obtain
$\sinh \theta_{ki}\cdot \Vert v_j\otimes c_{ijk}\Vert =\sinh \theta_{ki}\cdot\cosh q_j=\sinh h_k$.
Since $c_{ij}\ast c_{ijk}=-\cosh h_k$,
it follows that
\begin{equation*}
\frac{\partial \theta_i}{\partial f_j}
=\frac{1}{\sinh d_{ji}\sinh l_k}\cdot\frac{\sinh h_k}{-\cosh h_k}
=\frac{-1}{\sinh d_{ji}}\frac{\tanh h_k}{\sinh l_k}.
\end{equation*}

\noindent\textbf{(ii):}
Assume that $c_{ijk}$ is space-like.
Then $\Vert c_{ijk}\Vert ^2=1$ and $v_j*c_{ijk}=\pm\cosh q_j$.
It follows from (\ref{Eq: F15}) that
$\Vert v_j\otimes c_{ijk}\Vert =\sinh q_j$ if $q_j>0$, and $\Vert v_j\otimes c_{ijk}\Vert =-\sinh q_j$ if $q_j<0$.
We proceed to analyze different cases based on the domain to which $c_{ijk}$ belongs.

If $c_{ijk}\in D_{\mathrm{I}}\cup D_{\mathrm{II}}\cup D_{\mathrm{V}}\cup D_{\mathrm{VI}}$, then Table \ref{table} implies $q_j>0$ and $h_k>0$.
By Lemma \ref{Lem: space-like 1}, we have $\sinh \theta_{ki}\cdot \Vert v_j\otimes c_{ijk}\Vert =\sinh \theta_{ki}\cdot\sinh q_j=\cosh h_k$.
Additionally, it is known that $c_{ij}\ast c_{ijk}=-\sinh h_k$.

If $c_{ijk}\in D_{\mathrm{III}}$, then Table \ref{table} implies $q_j<0$ and $h_k>0$.
By Lemma \ref{Lem: space-like 1}, we have $\sinh \theta_{ki}\cdot \Vert v_j\otimes c_{ijk}\Vert =\sinh \theta_{ki}\cdot(-\sinh q_j)=\cosh h_k$ by Lemma \ref{Lem: space-like 1}.
Additionally, it is known that $c_{ij}\ast c_{ijk}=-\sinh h_k$.

If $c_{ijk}\in D_{\mathrm{IV}}$, then Table \ref{table} implies $q_j>0$ and $h_k<0$.
By Lemma \ref{Lem: space-like 1}, we have $\sinh \theta_{ki}\cdot \Vert v_j\otimes c_{ijk}\Vert =\sinh \theta_{ki}\cdot\sinh q_j=-\cosh h_k$. 
Note that $c_{ij}\ast c_{ijk}=-\sinh(d_\mathbb{H}(c_{ij},(c_{ijk})^\bot))=-\sinh (-h_k)=\sinh h_k$.

If $c_{ijk}\in D_{\mathrm{i}}\cup D_{\mathrm{ii}}\cup D_{\mathrm{vi}}$, then Table \ref{table} implies $q_j>0$ and $h_k>0$.
By Lemma \ref{Lem: space-like 2}, we have $\sinh \theta_{ki}\cdot \Vert v_j\otimes c_{ijk}\Vert =\sinh \theta_{ki}\cdot\sinh q_j=\cosh h_k$.
Additionally, $c_{ij}\ast c_{ijk}=-\sinh h_k$.

If $c_{ijk}\in D_{\mathrm{iii}}$, then Table \ref{table} implies $q_j<0$ and $h_k>0$.
By Lemma \ref{Lem: space-like 2}, we have $\sinh \theta_{ki}\cdot \Vert v_j\otimes c_{ijk}\Vert =\sinh \theta_{ki}\cdot(-\sinh q_j)=\cosh h_k$.
Additionally, $c_{ij}\ast c_{ijk}=-\sinh h_k$.

If $c_{ijk}\in D_{\mathrm{v}}$, then Table \ref{table} implies $q_j>0$ and $h_k<0$.
By Lemma \ref{Lem: space-like 2}, we have $\sinh \theta_{ki}\cdot \Vert v_j\otimes c_{ijk}\Vert =\sinh \theta_{ki}\cdot\sinh q_j=-\cosh h_k$.
Additionally, $c_{ij}\ast c_{ijk}=\sinh h_k$.

If $c_{ijk}\in D_{\mathrm{iv}}$, then Table \ref{table} implies $q_j<0$ and $h_k<0$.
By Lemma \ref{Lem: space-like 2}, we have $\sinh \theta_{ki}\cdot \Vert v_j\otimes c_{ijk}\Vert =\sinh \theta_{ki}\cdot(-\sinh q_j)=-\cosh h_k$.
Additionally, $c_{ij}\ast c_{ijk}=\sinh h_k$.

In all the above cases for space-like $c_{ijk}$, we can conclude that
\begin{equation*}
\frac{\partial \theta_i}{\partial f_j}
=\frac{-1}{\sinh d_{ji}}\frac{\coth h_k}{\sinh l_k}.
\end{equation*}

\noindent\textbf{(iii):}
The case where $c_{ijk}$ is light-like can be derived by continuity, specifically by considering the limit as the face center approaches the light cone.
In detail, for a light-like $c_{ijk}$, we have $\Vert c_{ijk}\Vert ^2=0$.
Combining (\ref{Eq: F9}) and (\ref{Eq: F16}) gives
\begin{equation*}
\frac{\sinh^2 \theta_{ki}\cdot \Vert v_j\otimes c_{ijk}\Vert ^2}{(c_{ij}\ast c_{ijk})^2}
=\frac{\Vert c_{ij}\otimes c_{ijk}\Vert ^2}{(c_{ij}\ast c_{ijk})^2}=1.
\end{equation*}
Furthermore, since $\Vert c_{ij}\otimes c_{ijk}\Vert =-(c_{ij}\ast c_{ijk})$, it follows that
\begin{equation*}
\frac{\partial \theta_i}{\partial f_j}
=\frac{-1}{\sinh d_{ji}\sinh l_k}.
\end{equation*}
This completes the proof.
\qed

\begin{remark}
The formulas derived in the proof of Lemma \ref{Lem: angle variations 1} essentially correspond to the sine and cosine laws for generalized right-angled hyperbolic triangles. 
To avoid the need for cumbersome case-by-case analyses, 
we adopt the Lorentzian inner product and Lorentzian cross product as a unifying framework for the calculations.
\end{remark}

\begin{lemma}\label{Lem: angle variations 2}
For a right-angled hyperbolic  hexagon $\{ijk\}\in F$, we have
\begin{equation}\label{Eq: angle 2}
\frac{\partial \theta^{jk}_i}{\partial f_i}
=\cosh l_{ij}\frac{\partial \theta^{ik}_j}{\partial f_i}
+\cosh l_{ik}\frac{\partial \theta^{ij}_k}{\partial f_i}.
\end{equation}
\end{lemma}
\proof
Under the same notations as used in the proof of Lemma \ref{Lem: angle variations 1}, analogous calculations yield
\begin{equation*}
\begin{aligned}
\frac{\partial \theta_i}{\partial f_i}
&=\frac{\partial \theta_i}{\partial l_i}\frac{\partial l_i}{\partial f_i}
+\frac{\partial \theta_i}{\partial l_j}\frac{\partial l_j}{\partial f_i}
+\frac{\partial \theta_i}{\partial l_k}\frac{\partial l_k}{\partial f_i}\\
&=\frac{-\sinh l_i\cosh \theta_k}{A}\cdot\coth d_{ik}
+\frac{-\sinh l_i\cosh \theta_j}{A}\cdot\coth d_{ij}\\
&=\frac{-\sinh l_i}{A}\cdot\left(\cosh \theta_k\cdot\coth d_{ik}+\cosh \theta_j\cdot\coth d_{ij}\right).
\end{aligned}
\end{equation*}
Similar to (\ref{Eq: F6}), we derive
\begin{gather*}
\frac{\partial \theta_j}{\partial f_i}
=\frac{\sinh l_j}{A}(\coth d_{ik}-\cosh \theta_i\coth d_{ij}), \\
\frac{\partial \theta_k}{\partial f_i}
=\frac{\sinh l_k}{A}(\coth d_{ij}-\cosh \theta_i\coth d_{ik}).
\end{gather*}
Then
\begin{equation*}
\begin{aligned}
&\, \cosh l_k\frac{\partial \theta_j}{\partial f_i}+\cosh l_j\frac{\partial \theta_k}{\partial f_i}\\
&=\frac{1}{A}(\cosh l_k\sinh l_j\coth d_{ik}
-\cosh l_k\sinh l_j\cosh \theta_i\coth d_{ij}\\
&\ \ \ \ \ \ \ \ \ \ +\cosh l_j\sinh l_k\coth d_{ij}
-\cosh l_j\sinh l_k\cosh \theta_i\coth d_{ik})\\
&=\frac{1}{A}[\coth d_{ik}\cdot(\cosh l_k\sinh l_j-\cosh l_j\sinh l_k\cdot\cosh \theta_i)\\
&\ \ \ \ \ \ \ \ \ \ +\coth d_{ij}\cdot(\cosh l_j\sinh l_k-\cosh l_k\sinh l_j\cosh \theta_i)]\\
&=\frac{1}{A}[\coth d_{ik}\cdot(\cosh l_k\sinh l_j-\cosh l_j\sinh l_k\cdot\frac{\cosh l_i+\cosh l_j\cosh l_k}{\sinh l_j\sinh l_k})\\
&\ \ \ \ \ \ \ \ \  \ +\coth d_{ij}\cdot(\cosh l_j\sinh l_k-\cosh l_k\sinh l_j\cdot\frac{\cosh l_i+\cosh l_j\cosh l_k}{\sinh l_j\sinh l_k})]\\
&=\frac{1}{A}(-\coth d_{ik}\sinh l_i\cosh \theta_k
-\coth d_{ij}\sinh l_i\cosh \theta_j)\\
&=\frac{\partial \theta_i}{\partial f_i}.
\end{aligned}
\end{equation*}
\qed

\begin{remark}
Formula (\ref{Eq: angle 2}) was first established by Glickenstein-Thomas \cite{GT} in the context of discrete conformal structures on closed surfaces, and it has since been widely applied.
For instance, relevant discussions can be found in \cite{Xu21, Xu 21a, Xu 21b, XZ22} and related works. 
In \cite{Xu22}, Xu showed that formula (\ref{Eq: angle 2}) remains valid for a specific class of discrete conformal structures on surfaces with boundary, and conjectured that this formula holds for generic discrete conformal structures on surfaces with boundary. 
Lemma \ref{Lem: angle variations 2} provides an affirmative resolution to Xu's conjecture.
\end{remark}

Given that the variational formula for the generalized angle (\ref{Eq: angle 1}) lacks symmetry in $i$ and $j$, we aim to derive a symmetric variational formula for the generalized angle via reparameterization.

\begin{lemma}\label{Lem: angle variations 3}
There exists a change of variables $u_i=u_i(f_i)$ such that
\begin{equation*}
\frac{\partial \theta^{jk}_i}{\partial u_j}
=\frac{\partial \theta^{ik}_j}{\partial u_i}.
\end{equation*}
The function $u$ is also referred to as a discrete conformal factor.
\end{lemma}
\proof
Suppose there exists a change of variables $u_i=u_i(f_i)$ such that $\partial \theta^{jk}_i/\partial u_j=\partial \theta^{ik}_j/\partial u_i$.
It follows from (\ref{Eq: angle 1}) that
\begin{equation*}
\frac{\partial \theta^{jk}_i}{\partial u_j}
=\frac{\partial \theta^{jk}_i}{\partial f_j}\frac{\partial f_j}{\partial u_j}
=\frac{-\partial f_j/\partial u_j}{\sinh d_{ji}}\cdot\frac{\tanh^\beta h_{ij}}{\sinh l_{ij}},
\end{equation*}
and
\begin{equation*}
\frac{\partial \theta^{ik}_j}{\partial u_i}
=\frac{-\partial f_i/\partial u_i}{\sinh d_{ij}}\cdot\frac{\tanh^\beta h_{ij}}{\sinh l_{ij}}.
\end{equation*}
This implies
\begin{equation*}
\frac{\partial f_i/\partial u_i}{\partial f_j/\partial u_j}=\frac{\sinh d_{ij}}{\sinh d_{ji}}.
\end{equation*}
From the proof of \cite[Theorem \ref{Thm: DCS}]{X-Z DCS1}, we recall the identities
\begin{equation*}
\frac{\sinh^2 d_{ij}}{\sinh^2 d_{ji}}=e^{C_{ij}}e^{2f_i-2f_j}  \quad \text{or} \quad
\frac{\sinh^2 d_{ij}}{\sinh^2 d_{ji}}
=\frac{1+\alpha_{i}e^{2f_i}}{1+\alpha_{j}e^{2f_j}}.
\end{equation*}
\begin{description}
\item[(i)]
By Remark \ref{Rmk: 1}, we set $e^{C_{ij}}=1$.
If $\frac{\sinh d_{ij}}{\sinh d_{ji}}=e^{f_i-f_j}$,
then
\begin{equation*}
\frac{\partial f_i}{\partial u_i}=e^{f_i}
\quad \text{and}\ \quad
\frac{\partial f_j}{\partial u_j}=e^{f_j}.
\end{equation*}
Suppose $\frac{\sinh d_{ij}}{\sinh d_{ji}}=-e^{f_i-f_j}$.
If $d_{ij}<0$, then $d_{ik}<0$ by Lemma \ref{Lem: time-space-like}, which implies $d_{ji}>0$ and $d_{ki}>0$.
For the condition (\ref{Eq: compatible condition}) to hold, we require $d_{jk}>0$ and $d_{kj}>0$.
It follows that $\frac{\sinh d_{ik}}{\sinh d_{ki}}=-e^{f_i-f_k}$ and $\frac{\sinh d_{jk}}{\sinh d_{ki}}=e^{f_j-f_k}$.
Hence
\begin{equation}\label{Eq: F26}
\frac{\partial f_i}{\partial u_i}=-e^{f_i},\ \quad
\frac{\partial f_j}{\partial u_j}=e^{f_j}
\quad \text{and}\ \quad
\frac{\partial f_k}{\partial u_k}=e^{f_k}.
\end{equation}
Similarly, if $d_{ji}<0$, then $\frac{\partial f_i}{\partial u_i}=e^{f_i}$, $\frac{\partial f_j}{\partial u_j}=-e^{f_j}$ and $\frac{\partial f_k}{\partial u_k}=e^{f_k}$.

\item[(ii)]
If $\frac{\sinh d_{ij}}{\sinh d_{ji}}
=\sqrt{\frac{1+\alpha_{i}e^{2f_i}}{1+\alpha_{j}e^{2f_j}}}$ with $1+\alpha_{i}e^{2f_i}>0$ and $1+\alpha_{j}e^{2f_j}> 0$,
then
\begin{equation*}
\frac{\partial f_i}{\partial u_i}=\sqrt{1+\alpha_ie^{2f_i}} \quad \text{and}\ \quad
\frac{\partial f_j}{\partial u_j}=\sqrt{1+\alpha_je^{2f_j}}.
\end{equation*}

If $\frac{\sinh d_{ij}}{\sinh d_{ji}}
=\sqrt{\frac{1+\alpha_{i}e^{2f_i}}{1+\alpha_{j}e^{2f_j}}}$ with $1+\alpha_{i}e^{2f_i}<0$ and $1+\alpha_{j}e^{2f_j}<0$,
then $\alpha\equiv -1$ and
\begin{equation*}
\frac{\partial f_i}{\partial u_i}=\sqrt{e^{2f_i}-1} \quad \text{and}\ \quad
\frac{\partial f_j}{\partial u_j}=\sqrt{e^{2f_j}-1}.
\end{equation*}

\item[(iii)]
If $\frac{\sinh d_{ij}}{\sinh d_{ji}}
=-\sqrt{\frac{1+\alpha_{i}e^{2f_i}}{1+\alpha_{j}e^{2f_j}}}$
with $1+\alpha_{i}e^{2f_i}>0$ and $1+\alpha_{j}e^{2f_j}>0$,
then by arguments analogous to those in the case (i),
when $d_{ij}<0$,
\begin{equation*}
\frac{\partial f_i}{\partial u_i}=-\sqrt{1+\alpha_ie^{2f_i}}, \quad
\frac{\partial f_j}{\partial u_j}=\sqrt{1+\alpha_je^{2f_j}} \quad \text{and}\ \quad
\frac{\partial f_k}{\partial u_k}=\sqrt{1+\alpha_ke^{2f_k}};
\end{equation*}
when $d_{ji}<0$,
\begin{equation*}
\frac{\partial f_i}{\partial u_i}=\sqrt{1+\alpha_ie^{2f_i}}, \quad
\frac{\partial f_j}{\partial u_j}=-\sqrt{1+\alpha_je^{2f_j}} \quad \text{and}\  \quad
\frac{\partial f_k}{\partial u_k}=\sqrt{1+\alpha_ke^{2f_k}}.
\end{equation*}
If $\frac{\sinh d_{ij}}{\sinh d_{ji}}
=-\sqrt{\frac{1+\alpha_{i}e^{2f_i}}{1+\alpha_{j}e^{2f_j}}}$
with $1+\alpha_{i}e^{2f_i}<0$ and $1+\alpha_{j}e^{2f_j}<0$,
then
when $d_{ij}<0$,
\begin{equation}\label{Eq: F50}
\frac{\partial f_i}{\partial u_i}=-\sqrt{e^{2f_i}-1}, \quad
\frac{\partial f_j}{\partial u_j}=\sqrt{e^{2f_j}-1} \quad \text{and}\ \quad
\frac{\partial f_k}{\partial u_k}=\sqrt{e^{2f_k}-1};
\end{equation}
when $d_{ji}<0$,
\begin{equation*}
\frac{\partial f_i}{\partial u_i}=\sqrt{e^{2f_i}-1}, \quad
\frac{\partial f_j}{\partial u_j}=-\sqrt{e^{2f_j}-1} \quad \text{and}\  \quad
\frac{\partial f_k}{\partial u_k}=\sqrt{e^{2f_k}-1}.
\end{equation*}
\end{description}
\qed

\begin{remark}\label{Rmk: 2}
In fact, the function $u_i=u_i(f_i)$ admits an explicit expression.
\begin{description}
\item[(i)] 
If $\frac{\partial f_i}{\partial u_i}=e^{f_i}$, then $u_i=-e^{-f_i}$.
If $\frac{\partial f_i}{\partial u_i}=-e^{f_i}$,
then $u_i=e^{-f_i}$.

\item[(ii)]
Suppose $\frac{\partial f_i}{\partial u_i}=\sqrt{1+\alpha_ie^{2f_i}}$,
if $\alpha_i=0$, then $u_i=f_i\in \mathbb{R}$;
if $\alpha_i\neq0$, then
\begin{equation*}
u_i=\frac{1}{2}\ln\left|\frac{\sqrt{1+\alpha_ie^{2f_i}}-1}
{\sqrt{1+\alpha_ie^{2f_i}}+1}\right|.
\end{equation*}
In particular, if $\alpha_i=-1$, then
\begin{equation}\label{Eq: F33}
-\tanh u_i=\sqrt{1-e^{2f_i}},
\end{equation}
which implies $e^{f_i}=\frac{1}{\cosh u_i}$.
Note that $\sqrt{1-e^{2f_i}}\neq0$ and thus $f_i\neq0$.
Consequently, $f_i\in \mathbb{R}_{<0}$ and $u_i\in \mathbb{R}_{<0}$.
If $\alpha_i=1$, then
\begin{equation}\label{Eq: F34}
-\coth u_i=\sqrt{1+e^{2f_i}},
\end{equation}
which implies $e^{f_i}=-\frac{1}{\sinh u_i}$ with $f_i\in \mathbb{R}$ and $u_i\in \mathbb{R}_{<0}$.

\item[(iii)]
Suppose $\frac{\partial f_i}{\partial u_i}=-\sqrt{1+\alpha_ie^{2f_i}}$,
if $\alpha_i=0$, then $u_i=-f_i\in \mathbb{R}$;
if $\alpha_i\neq0$, then
\begin{equation*}
u_i=-\frac{1}{2}\ln\left|\frac{\sqrt{1+\alpha_ie^{2f_i}}-1}
{\sqrt{1+\alpha_ie^{2f_i}}+1}\right|.
\end{equation*}
In particular, if $\alpha_i=-1$, then
\begin{equation*}
\tanh u_i=\sqrt{1-e^{2f_i}},
\end{equation*}
which implies $e^{f_i}=\frac{1}{\cosh u_i}$.
Similarly, the case that $f_i=0$ is excluded,
so $f_i\in \mathbb{R}_{<0}$ and $u_i\in \mathbb{R}_{>0}$.
If $\alpha_i=1$, then
\begin{equation*}
\coth u_i=\sqrt{1+e^{2f_i}},
\end{equation*}
which implies $e^{f_i}=\frac{1}{\sinh u_i}$ with $f_i\in \mathbb{R}$ and $u_i\in \mathbb{R}_{>0}$.

\item[(iv)]
If $\frac{\partial f_i}{\partial u_i}=\sqrt{e^{2f_i}-1}$,
then
\begin{gather*}
\tan u_i=\sqrt{e^{2f_i}-1},\ \text{for}\ u_i\in (0,\frac{\pi}{2}),\\
-\frac{1}{\tan u_i}=\sqrt{e^{2f_i}-1},\ \text{for}\ u_i\in (-\frac{\pi}{2},0).
\end{gather*}
This implies that $e^{f_i}=\frac{1}{\cos u_i}$ for $u_i\in (0,\frac{\pi}{2})$, and $e^{f_i}=-\frac{1}{\sin u_i}$ for $u_i\in (-\frac{\pi}{2},0)$.
If $\frac{\partial f_i}{\partial u_i}=-\sqrt{e^{2f_i}-1}$,
then
\begin{gather*}
-\tan u_i=\sqrt{e^{2f_i}-1},\ \text{for}\ u_i\in (-\frac{\pi}{2},0),\\
\frac{1}{\tan u_i}=\sqrt{e^{2f_i}-1},\ \text{for}\ u_i\in (0,\frac{\pi}{2}).
\end{gather*}
This implies that $e^{f_i}=\frac{1}{\cos u_i}$ for $u_i\in (-\frac{\pi}{2},0)$, and $e^{f_i}=\frac{1}{\sin u_i}$ for $u_i\in (0,\frac{\pi}{2})$.
\end{description}
\end{remark}

\section{Rigidity and existence of the discrete conformal structure (A1)}\label{section A1}

\subsection{Admissible space of the discrete conformal structure (A1)}\label{section a4}

Suppose $(\Sigma,\mathcal{T},\alpha,\eta)$ is a weighted triangulated surface with boundary, 
where the weights are given by $\alpha: B\rightarrow \{-1,0,1\}$ and $\eta\in \mathbb{R}^E_{>0}$.
For a right-angled hyperbolic hexagon $\{ijk\}\in F$ with edge lengths $l_{ij}, l_{jk}, l_{ki}$ given by (\ref{Eq: DCS3}), 
the admissible space of the discrete conformal factors $f$ is defined as
\begin{equation}\label{Eq: admissible space a1}
\begin{aligned}
\Omega_{ijk}(\eta)=\{(f_i,f_j,f_k)\in \mathbb{R}^3\mid &\cosh l_{rs}=-\sqrt{(1+\alpha_re^{2f_r})(1+\alpha_se^{2f_s})}\\
& \ \ \ \ \ \ \ \ \ \ \ \ \ \ \ \ \
+\eta_{rs}e^{f_r+f_s}>1,\  \forall\{r,s\}\subseteq\{i,j,k\}\}.
\end{aligned}
\end{equation}
By a change of variables $u_i=u_i(f_i)$ in Lemma \ref{Lem: angle variations 3},
the admissible space $\Omega_{ijk}(\eta)$ is transformed into the following admissible space $\mathcal{U}_{ijk}(\eta)$ of the discrete conformal factors $u$ for a right-angled hyperbolic hexagon $\{ijk\}$:
\begin{equation*}
\mathcal{U}_{ijk}(\eta)=u(\Omega_{ijk}(\eta))
=\{(u_i,u_j,u_k)\in \mathbb{R}^{N_1}\times \mathbb{R}_{<0}^{3-N_1}\mid (f_i,f_j,f_k)\in \Omega_{ijk}(\eta)\},
\end{equation*}
where $N_1$ is the number of the boundary components $i\in B$ with $\alpha_i=0$.

Since $\frac{\sinh d_{rs}}{\sinh d_{sr}}
=\sqrt{\frac{1+\alpha_re^{2f_r}}{1+\alpha_se^{2f_s}}}>0$ with $1+\alpha_{r}e^{2f_r}>0$ and $1+\alpha_{s}e^{2f_s}>0$ in Theorem \ref{Thm: DCS},
it follows from the proof of Lemma \ref{Lem: angle variations 3} that
\begin{equation}\label{Eq: F32}
\frac{\partial f_r}{\partial u_r}
=\sqrt{1+\alpha_re^{2f_r}} \quad \text{and}\ \quad
\frac{\partial f_s}{\partial u_s}=\sqrt{1+\alpha_se^{2f_s}}.
\end{equation}
By Remark \ref{Rmk: 2}, the following results are obtained:
\begin{description}
\item[(i)]
If $\alpha\equiv0$, then $\cosh l_{rs}=-1+\eta_{rs}e^{u_r+u_s}>1$,
which implies $u_r+u_s>\log \frac{2}{\eta_{rs}}$.
\item[(ii)]
If $\alpha\equiv-1$,
then $\cosh l_{rs}=-\sqrt{(1-e^{2f_r})(1-e^{2f_s})}
+\eta_{rs}e^{f_r+f_s}
=-\tanh u_r\tanh u_s+\eta_{rs}\frac{1}{\cosh u_r}\frac{1}{\cosh u_s}>1$,
which implies $\cosh\, (u_r+u_s)<\eta_{rs}$.
Since $u\in \mathbb{R}^3_{<0}$,
it follows that $u_r+u_s>-\mathrm{arccosh}\, \eta_{rs}$.
Note that $\eta_{rs}>1$.
\item[(iii)]
If $\alpha\equiv1$, then $\cosh l_{rs}=-\sqrt{(1+e^{2f_r})(1+e^{2f_s})}
+\eta_{rs}e^{f_r+f_s}=-\coth u_r\coth u_s+\eta_{rs}\frac{1}{\sinh u_r}\frac{1}{\sinh u_s}>1$,
which implies $\cosh\, (u_r+u_s)<\eta_{rs}$.
Since $u\in \mathbb{R}^3_{<0}$,
it follows that $u_r+u_s>-\mathrm{arccosh}\, \eta_{rs}$.
Note that $\eta_{rs}>1$.
\item[(iv)]
If $\alpha_r=0$ and $\alpha_s=-1$, then
$\cosh l_{rs}=-\sqrt{1-e^{2f_s}}
+\eta_{rs}e^{f_r+f_s}=\tanh u_s+\eta_{rs}e^{u_r}\frac{1}{\cosh u_s}>1$,
which implies $u_r+u_s>\log \frac{1}{\eta_{rs}}$.
\item[(v)]
If $\alpha_r=0$ and $\alpha_s=1$, then
$\cosh l_{rs}=-\sqrt{1+e^{2f_s}}
+\eta_{rs}e^{f_r+f_s}=\coth u_s+\eta_{rs}e^{u_r}\frac{-1}{\sinh u_s}>1$.
Since $u_s\in \mathbb{R}_{<0}$,
it follows that $u_r+u_s>\log \frac{1}{\eta_{rs}}$.
\item[(vi)]
If $\alpha_r=1$ and $\alpha_s=-1$, then
$\cosh l_{rs}=-\sqrt{(1+e^{2f_r})(1-e^{2f_s})}
+\eta_{rs}e^{f_r+f_s}=-\coth u_r\tanh u_s+\eta_{rs}\frac{-1}{\sinh u_r}\frac{1}{\cosh u_s}>1$.
Since $u_r\in \mathbb{R}_{<0}$ and $u_s\in \mathbb{R}_{<0}$,
it follows that $\sinh\,(u_r+u_s)>-\eta_{rs}$, and hence $u_r+u_s>\mathrm{arcsinh}\, (-\eta_{rs})$.
\end{description}
In summary, the admissible space $\mathcal{U}_{ijk}(\eta)$ can be rewritten as
\begin{equation}\label{Eq: admissible space a2}
\mathcal{U}_{ijk}(\eta)
=\{(u_i,u_j,u_k)\in \mathbb{R}^{N_1}\times \mathbb{R}_{<0}^{3-N_1}\mid u_r+u_s>C(\eta_{rs}),\ \forall\{r,s\}\subseteq\{i,j,k\}\},
\end{equation}
where $C(\eta_{rs})$ is a constant depending on $\eta_{rs}$.

\begin{theorem}\label{Thm: ASC 1}
Suppose $(\Sigma,\mathcal{T},\alpha,\eta)$ is a weighted triangulated surface with boundary, where the weighs are given by $\alpha: B\rightarrow \{-1,0,1\}$ and $\eta\in \mathbb{R}_{>0}^E$ satisfying $\eta_{ij}>\alpha_i\alpha_j$ for any two adjacent boundary components $i,j\in B$.
Then the admissible space $\mathcal{U}_{ijk}(\eta)$ is a convex polytope.
As a result, the admissible space
$\mathcal{U}(\eta)=\bigcap_{\{ijk\}\in F}\mathcal{U}_{ijk}(\eta)$
is also a convex polytope.
\end{theorem}

\subsection{Negative definiteness of the Jacobian}
\label{subsec: matrix}

For a hyper-ideal hyperbolic triangle $\{ijk\}\in F$,
we introduce some simplifying notations: $l_r=l_{st}$, $\theta_r=\theta^{st}_r$, and $A=\sinh l_r\sinh l_s\sinh \theta_t$, where $\{r,s,t\}=\{i,j,k\}$.
Define
\begin{equation}\label{Eq: Q_1}
Q_1=\left(
   \begin{array}{ccc}
     -1 & \cosh\theta_k & \cosh\theta_j \\
     \cosh\theta_k & -1 & \cosh\theta_i \\
     \cosh\theta_j & \cosh\theta_i & -1 \\
   \end{array}
 \right).
\end{equation}
Then $\det Q_1=-1+2\cosh \theta_i\cosh \theta_j\cosh \theta_k+\cosh^2 \theta_i+\cosh^2 \theta_j+\cosh^2 \theta_k>0$.
Following calculations analogous to those in \eqref{Eq: F24}, we obtain
\begin{equation*}
\begin{aligned}
\frac{\partial(\theta_i,\theta_j,\theta_k)}
{\partial(l_i,l_j,l_k)}
=&-\frac{1}{A}
 \left(
   \begin{array}{ccc}
     \sinh l_i & 0 & 0 \\
     0 & \sinh l_j & 0 \\
     0 & 0 & \sinh l_k \\
   \end{array}
 \right)\times
 Q_1.
\end{aligned}
\end{equation*}
This implies
\begin{equation}\label{Eq: F35}
\det\bigg(\frac{\partial(\theta_i,\theta_j,\theta_k)}
{\partial(l_i,l_j,l_k)}\bigg)
=-\frac{1}{A^3}\sinh l_i\sinh l_j\sinh l_k\cdot \det Q_1<0.
\end{equation}
By combining \eqref{Eq: variation 1} and \eqref{Eq: variation 2}, we derive
\begin{equation}\label{Eq: Q_2}
\begin{aligned}
\frac{\partial(l_i,l_j,l_k)}
{\partial(f_i,f_j,f_k)}=
\left(
   \begin{array}{ccc}
     0 & \coth d_{jk} & \coth d_{kj} \\
     \coth d_{ik} & 0 & \coth d_{ki} \\
     \coth d_{ij} & \coth d_{ji} & 0 \\
   \end{array}
 \right)
:=Q_2.
\end{aligned}
\end{equation}
Then
\begin{equation}\label{Eq: F23}
\det Q_2=\coth d_{jk}\coth d_{ki}\coth d_{ij}
+\coth d_{kj}\coth d_{ik}\coth d_{ji}.
\end{equation}
For any $\alpha_i\in \{-1,0,1\}$,
set $S_i=e^{f_i}$ and $C_i=\sqrt{1+\alpha_ie^{2f_i}}$,
so $S_i>0$, $C_i>0$, and $C_i^2-\alpha_iS_i^2=1$.
It follows from (\ref{Eq: F32}) that
\begin{equation}\label{Eq: Q_3}
\frac{\partial(f_i,f_j,f_k)}
{\partial(u_i,u_j,u_k)}
=\left(
   \begin{array}{ccc}
     C_i & 0 & 0 \\
     0 & C_j & 0 \\
     0 & 0 & C_k \\
   \end{array}
 \right)
 :=Q_3.
\end{equation}

\begin{theorem}\label{Thm: matrix negative 1}
Under the same assumptions as those in Theorem \ref{Thm: ASC 1},
for a right-angled hyperbolic hexagon $\{ijk\}\in F$ on $\mathcal{U}_{ijk}(\eta)$,
the Jacobian  $\Lambda_{ijk}=\frac{\partial(\theta^{jk}_i,\theta^{ik}_j,\theta^{ij}_k)}
{\partial(u_i,u_j,u_k)}$ is symmetric and negative definite.
\end{theorem}
\proof
Since $d>0$ for the discrete conformal structure (A1) in Theorem \ref{Thm: DCS},
then $\det Q_2>0$ by (\ref{Eq: F23}).
Consequently,
\begin{equation*}
\det\bigg(\frac{\partial(l_i,l_j,l_k)}
{\partial(u_i,u_j,u_k)}\bigg)
=\det Q_2\cdot \det Q_3>0.
\end{equation*}
Combining with (\ref{Eq: F35}) gives
\begin{equation}\label{Eq: F87}
\det\Lambda_{ijk}
=\det\bigg(\frac{\partial(\theta^{jk}_i,\theta^{ik}_j,\theta^{ij}_k)}
{\partial(l_i,l_j,l_k)}\bigg)\cdot\det\bigg(\frac{\partial(l_i,l_j,l_k)}
{\partial(u_i,u_j,u_k)}\bigg)<0.
\end{equation}
This implies the matrix $\Lambda_{ijk}$ is non-degenerate.
The symmetry of $\Lambda_{ijk}$ follows from Lemma \ref{Lem: angle variations 3}.
Set $a_{rrs}:=C_r\cdot\coth d_{rs}>0$.
Define
\begin{equation*}
Q:=Q_1Q_2Q_3=\left(
   \begin{array}{ccc}
     -1 & \cosh\theta_k & \cosh\theta_j \\
     \cosh\theta_k & -1 & \cosh\theta_i \\
     \cosh\theta_j & \cosh\theta_i & -1 \\
   \end{array}
 \right)\\
 \times
  \left(
   \begin{array}{ccc}
     0 & a_{jjk} & a_{kkj} \\
     a_{iik} & 0 & a_{kki} \\
     a_{iij} & a_{jji} & 0 \\
   \end{array}
   \right).
\end{equation*}
Then $\Lambda_{ijk}$ can be expressed as
\begin{equation*}
\Lambda_{ijk}=-\frac{1}{A}
 \left(
   \begin{array}{ccc}
     \sinh l_i & 0 & 0 \\
     0 & \sinh l_j & 0 \\
     0 & 0 & \sinh l_k \\
   \end{array}
 \right)\times Q.
\end{equation*}
To prove the negative definiteness of $\Lambda_{ijk}$,
it suffices to show that $Q$ is positive definite.
Formula (\ref{Eq: F87}) shows that $\det Q>0$.
It remains to verify that the $1\times 1$ and $2\times 2$ leading principal minors of $Q$ are positive. 
Let $m_{rs}$ denote the entry of $Q$ at $r$-th row and $s$-th column.
Then
\begin{equation*}
\begin{aligned}
&m_{11}=\cosh \theta_k\cdot a_{iik}+\cosh \theta_j\cdot a_{iij}\\
&m_{22}=\cosh \theta_k\cdot a_{jjk}+\cosh \theta_i\cdot a_{jji}\\
&m_{12}=-a_{jjk}+\cosh \theta_j\cdot a_{jji}\\
&m_{21}=-a_{iik}+\cosh \theta_i\cdot a_{iij}.
\end{aligned}
\end{equation*}
Thus the $1\times 1$ leading principal minor of $Q$ is $m_{11}>0$ due to $a_{iik}>0$ and $a_{iij}>0$.
And the $2\times 2$ leading principal minor of $Q$ is
\begin{equation*}
\begin{aligned}
m_{11}m_{22}-m_{12}m_{21}
&=\sinh^2 \theta_ka_{iik}a_{jjk}
+(\cosh \theta_i\cosh\theta_k+\cosh \theta_j)a_{iik}a_{jji}\\
&\ \ \ \ +(\cosh \theta_j\cosh\theta_k+\cosh \theta_i)a_{iij}a_{jjk}
>0.
\end{aligned}
\end{equation*}
This completes the proof.
\qed

As a consequence, we obtain the following corollary.
\begin{corollary}\label{Cor: matrix 1}
Under the same assumptions as those in Theorem \ref{Thm: ASC 1},
the Jacobian $\Lambda=\frac{\partial (K_i,..., K_N)}{\partial(u_i,...,u_N)}$ is symmetric and negative definite on $\mathcal{U}(\eta)$.
\end{corollary}
\proof
The conclusion follows directly from Theorem \ref{Thm: matrix negative 1} and the identity $\Lambda=\sum_{\{ijk\}\in F}\Lambda_{ijk}$.
Specifically, each $\Lambda_{ijk}$ is extended to an $N\times N$ matrix via zero-padding, such that it acts non-trivially on a vector $(1,...,N)$ exclusively through the coordinates corresponding to the boundary components $i,j,k$ of the right-angled hyperbolic hexagon $\{ijk\}\in F$.
\qed

\subsection{Rigidity of the discrete conformal structure (A1)}

The following theorem establishes the rigidity of the discrete conformal structure (A1), 
which corresponds to the rigidity part of Theorem \ref{Thm: rigidity and image} (i).

\begin{theorem}\label{Thm: rigidity 1}
Under the same assumptions as those in Theorem \ref{Thm: ASC 1},
the discrete conformal factor $f$ is uniquely determined by its generalized combinatorial curvature $K\in\mathbb{R}^N_{>0}$.
%In particular, the map from $f$ to $K$ is a smooth embedding.
\end{theorem}
\proof
Theorem \ref{Thm: ASC 1} and Theorem \ref{Thm: matrix negative 1} imply that for a right-angled hyperbolic  hexagon $\{ijk\}\in F$,
the following energy function
\begin{equation*}
\mathcal{F}_{ijk}(u_i,u_j,u_k)=\int^{(u_i,u_j,u_k)}
\theta^{jk}_idu_i+\theta^{ik}_jdu_j+\theta^{ij}_kdu_k
\end{equation*}
is a well-defined smooth function on $\mathcal{U}_{ijk}(\eta)$.
Furthermore, $\mathcal{F}_{ijk}(u_i,u_j,u_k)$ is strictly concave on $\mathcal{U}_{ijk}(\eta)$ with $\nabla_{u_i}\mathcal{F}_{ijk}=\theta^{jk}_i$.
Define the function $\widetilde{\mathcal{F}}:\mathcal{U}(\eta)\rightarrow \mathbb{R}$ by
\begin{equation*}
\widetilde{\mathcal{F}}(u)=\sum_{\{ijk\}\in F}\mathcal{F}_{ijk}(u_i,u_j,u_k).
\end{equation*}
By Corollary \ref{Cor: matrix 1},
$\widetilde{\mathcal{F}}(u)$ is strictly concave on $\mathcal{U}(\eta)$ with $\nabla_{u_i}\widetilde{\mathcal{F}}(u)=K_i$.
By Lemma \ref{Lem: analysis}, the gradient map $\nabla \widetilde{\mathcal{F}}: \mathcal{U}(\eta)\rightarrow \mathbb{R}^N$ is a smooth embedding.
Note that $\mathcal{U}(\eta)$ and $\Omega(\eta)=\bigcap_{\{ijk\}\in F}\Omega_{ijk}(\eta)$ are diffeomorphic.
Consequently, the map from $\{f\in \Omega(\eta)\}$ to $\{K\in \mathbb{R}^N_{>0}\}$ is a smooth injective map.
\qed

\subsection{Existence of the discrete conformal structure (A1)}

The following theorem characterizes the image of $K$ for the discrete conformal structure (A1), 
which corresponds to the existence part of Theorem \ref{Thm: rigidity and image} (i).

\begin{theorem}\label{Thm: image 1}
Suppose $(\Sigma,\mathcal{T},\alpha,\eta)$ is a weighted triangulated surface with boundary, where the weights are given by $\alpha: B\rightarrow \{0,1\}$ and $\eta\in \mathbb{R}_{>0}^E$ satisfying $\eta_{ij}>\alpha_i\alpha_j$ for any two adjacent boundary components $i,j\in B$.
The image of $K$ is $\mathbb{R}^N_{>0}$.
\end{theorem}

To prove Theorem \ref{Thm: image 1}, we need the following two lemmas.

\begin{lemma}\label{Lem: limit a1}
Under the same assumptions as those in Theorem \ref{Thm: ASC 1},
for a right-angled hyperbolic hexagon $\{ijk\}\in F$ on $\Omega_{ijk}(\eta)$,
if one of the following conditions is satisfied
\begin{description}
\item[(i)] $\lim f_i=+\infty,\ \lim f_j=+\infty,\ \lim f_k=+\infty$;
\item[(ii)] $\lim f_i=+\infty,\ \lim f_j=+\infty,\ \lim f_k=c$;
\item[(iii)] $\lim f_i=+\infty,\ \lim f_j=a,\ \lim f_k=b$.
\end{description}
where $a, b,c$ are constants.
Then $\lim\theta^{jk}_i(f_i,f_j,f_k)=0$.
\end{lemma}
\proof
For the case (i), it follows from (\ref{Eq: DCS3}) that
\begin{equation}\label{Eq: F57}
\begin{aligned}
\cosh l_{ij}=&\bigg[\eta_{ij}-\sqrt{(e^{-2f_i}+\alpha_i)
(e^{-2f_j}+\alpha_j)}\bigg]e^{f_i+f_j},\\
\cosh l_{jk}=&\bigg[\eta_{jk}-\sqrt{(e^{-2f_j}+\alpha_j)
(e^{-2f_k}+\alpha_k)}\bigg]e^{f_j+f_k},\\
\cosh l_{ik}=&\bigg[\eta_{ik}-\sqrt{(e^{-2f_i}+\alpha_i)
(e^{-2f_k}+\alpha_k)}\bigg]e^{f_i+f_k}.
\end{aligned}
\end{equation}
Then $\lim\cosh l_{ij}:=\lim c_ke^{f_i+f_j}$,\ $\lim\cosh l_{jk}:=\lim c_ie^{f_j+f_k}$ and $\lim\cosh l_{ik}:=\lim c_je^{f_i+f_k}$,
where $c_i, c_j, c_k$ are positive constants.
Set $L_1=\frac{\cosh l_{jk}}{\sinh l_{ij}\sinh l_{ik}}$ and $L_2=\frac{\cosh l_{ij}\cosh l_{ik}}{\sinh l_{ij}\sinh l_{ik}}$,
then $\cosh \theta^{jk}_i=L_1+L_2$.
Since
\begin{equation*}
L_1\rightarrow
\frac{c_ie^{f_j+f_k}}{c_ke^{f_i+f_j}\cdot c_je^{f_i+f_k}}
\rightarrow 0 \quad \ \text{and}\ \quad
L_2\rightarrow 1,
\end{equation*}
it follows that $\cosh \theta^{jk}_i\rightarrow 1$, which implies $\theta_i^{jk}\rightarrow 0$.
The same arguments apply for the cases (ii) and (iii), the details of which are omitted here for brevity.
\qed

\begin{lemma}\label{Lem: limit a2}
Under the same assumptions as those in Theorem \ref{Thm: ASC 1},
for a right-angled hyperbolic  hexagon $\{ijk\}\in F$ on $\Omega_{ijk}(\eta)$, if $\lim f_i=-\infty$,
then $\lim f_j=\lim f_k=+\infty$.
Furthermore, $\lim\theta^{jk}_i(f_i,f_j,f_k)=+\infty$.
\end{lemma}
\proof
The definition of the admissible space $\mathcal{U}_{ijk}(\eta)$ in (\ref{Eq: admissible space a2}) implies $u_r+u_s>C(\eta_{rs})$ for any $\{r,s\}\subseteq\{i,j,k\}$.
Hence, there exists no subset $\{r,s\}\subseteq\{i,j,k\}$ such that $\lim u_r=-\infty$ and either $\lim u_s=-\infty$ or $\lim u_s=c$, where $c$ is a constant.
By combining (\ref{Eq: F32}) and Remark \ref{Rmk: 2}, the following hold: 
if $\lim f_i=-\infty$, then $\lim u_i=-\infty$;
if $\lim f_j=+\infty$, then $\lim u_j=+\infty$ when $\alpha_j=0$ and $\lim u_j=0^-$ when $\alpha_j=1$.
Therefore, if $\lim f_i=-\infty$, then $\lim u_i=-\infty$, which implies $\lim u_j=\lim u_k=+\infty$. 
This further yields $\alpha_j=\alpha_k=0$ and $\lim f_j=\lim f_k=+\infty$.

If $\alpha_i=0$, then
$\cosh l_{ij}=-1+\eta_{ij}e^{u_i+u_j}$, 
$\cosh l_{jk}=-1+\eta_{jk}e^{u_j+u_k}$, and
$\cosh l_{ik}=-1+\eta_{ik}e^{u_i+u_k}$.
Thus
\begin{equation*}
\begin{aligned}
L_1&=\frac{-1+\eta_{jk}e^{u_j+u_k}}
{\sqrt{(-1+\eta_{ij}e^{u_i+u_j})^2-1}\cdot
\sqrt{(-1+\eta_{ik}e^{u_i+u_k})^2-1}}\\
&=\frac{-1+\eta_{jk}e^{u_j+u_k}}
{\sqrt{\eta_{ij}\eta_{ik}}e^{u_j+u_k}}\cdot
\frac{1}{\sqrt{\eta_{ij}e^{2u_i}-2e^{u_i-u_j}}}\cdot
\frac{1}{\sqrt{\eta_{ik}e^{2u_i}-2e^{u_i-u_k}}}
\rightarrow +\infty.
\end{aligned}
\end{equation*}
Since $L_2>0$, it follows that $\cosh \theta^{jk}_i=L_1+L_2\rightarrow +\infty$,
which implies $\theta^{jk}_i\rightarrow +\infty$.

If $\alpha_i=1$, then $\sqrt{1+e^{2f_i}}=-\coth u_i$
and $e^{f_i}=-\frac{1}{\sinh u_i}$.
Hence,
\begin{equation*}
\cosh l_{ij}=-\sqrt{1+e^{2f_i}}+\eta_{ij}e^{f_i+f_j}
=\frac{-\cosh u_i+\eta_{ij}e^{u_j}}{-\sinh u_i}\rightarrow-1+2\eta_{ij}e^{u_i+u_j}.
\end{equation*}
Similarly, $\cosh l_{ik}\rightarrow-1+2\eta_{ik}e^{u_i+u_k}$ and
$\cosh l_{jk}=-1+\eta_{jk}e^{u_j+u_k}$.
By reasoning analogous to the case $\alpha_i=0$,
we have $\cosh \theta^{jk}_i=L_1+L_2\rightarrow +\infty$, which implies $\theta^{jk}_i\rightarrow +\infty$.

If $\alpha_i=-1$, then
$\sqrt{1-e^{2f_i}}=-\tanh u_i$ and $e^{f_i}=\frac{1}{\cosh u_i}$.
Hence,
\begin{equation*}
\cosh l_{ij}=-\sqrt{1-e^{2f_i}}+\eta_{ij}e^{f_i+f_j}
=\frac{\sinh u_i+\eta_{ij}e^{u_j}}{\cosh u_i}\rightarrow-1+2\eta_{ij}e^{u_i+u_j}.
\end{equation*}
Similarly, $\cosh l_{ik}\rightarrow-1+2\eta_{ik}e^{u_i+u_k}$ and
$\cosh l_{jk}=-1+\eta_{jk}e^{u_j+u_k}$.
Similar to the proof of the case $\alpha_i=0$,
we have $\cosh \theta^{jk}_i=L_1+L_2\rightarrow +\infty$,
which implies $\theta^{jk}_i\rightarrow +\infty$.
\qed

\noindent\textbf{Proof of Theorem \ref{Thm: image 1}:}
Let $X$ be the image of $K$.
To prove $X=\mathbb{R}^N_{>0}$, we aim to show that $X$ is both an open and a closed subset of $\mathbb{R}^N_{>0}$.
By the definition of the curvature $K$, it follows directly that $X\subseteq \mathbb{R}^N_{>0}$.
Since the map from $f\in \Omega(\eta)$ to $\mathbb{R}^N_{>0}$ is a smooth injective map,
its image $X$ must be an open subset of  $\mathbb{R}^N_{>0}$.
It remains to verify that $X$ is closed in $\mathbb{R}^N_{>0}$.
Suppose $\{f^{(m)}\}$ is a sequence in $\Omega(\eta)$ such that $\lim_{m\rightarrow +\infty}K^{(m)}\in \mathbb{R}^N_{>0}$.
We shall show that there exists a subsequence,
still denoted by $\{f^{(m)}\}$, such that $\lim_{m\rightarrow +\infty}f^{(m)}\in\Omega(\eta)$.

Suppose, for contradiction, that there exists a subsequence $\{f^{(m)}\}$ such that $\lim_{m\rightarrow +\infty}f^{(m)}$ lies on the boundary of $\Omega(\eta)$.
By (\ref{Eq: admissible space a1}), the boundary of the admissible space $\Omega(\eta)$ in $[-\infty,+\infty]^N$ consists of three distinct components, specified as follows:
\begin{description}
\item[(i)]
There exists an ideal edge $\{rs\}\in E$ such that $l_{rs}=0$, a condition equivalent to $u_r+u_s=C(\eta_{rs})$ for the admissible space $\mathcal{U}_{ijk}(\eta)$ defined by (\ref{Eq: admissible space a2});
\item[(ii)]
For $\alpha\in \{0,1\}$, there exists $i\in B$ such that $f_i^{(m)}\rightarrow +\infty$;
\item[(iii)]
For $\alpha\in \{0,1\}$, there exists $i\in B$ such that $f_i^{(m)}\rightarrow -\infty$.
\end{description}

For the case (i), consider a right-angled hyperbolic hexagon $\{ijk\}\in F$ with edge lengths $l_{ij},l_{jk},l_{ik}$ and opposite boundary lengths $\theta^{ij}_k,\theta^{jk}_i,\theta^{ik}_j$.
Without loss of generality, we assume $l_{ij}\rightarrow 0$. By the cosine law, we have
\begin{equation*}
\cosh \theta^{ik}_j=\frac{\cosh l_{ik}+\cosh l_{ij}\cosh l_{jk}}{\sinh l_{ij}\sinh l_{jk}}>\frac{\cosh l_{ij}\cosh l_{jk}}{\sinh l_{ij}\sinh l_{jk}}\geq\frac{\cosh l_{ij}}{\sinh l_{ij}}\rightarrow +\infty,
\end{equation*}
which contradicts the assumption that $\lim_{m\rightarrow +\infty}K^{(m)}\in \mathbb{R}^N_{>0}$.
For the case (ii), by Lemma \ref{Lem: limit a1}, each generalized angle incident to the boundary component $i\in B$ converges to 0.
Hence, $\lim_{m\rightarrow +\infty}K_i^{(m)}=0$.
For the case (iii), Lemma \ref{Lem: limit a2} implies that $\lim_{m\rightarrow +\infty}K_i^{(m)}=+\infty$.
Both cases contradict the assumption that $\lim_{m\rightarrow +\infty}K^{(m)}\in \mathbb{R}^N_{>0}$.
\qed

\begin{remark}
Note that Lemma \ref{Lem: limit a1} and Lemma \ref{Lem: limit a2} are valid for $\alpha: B\rightarrow \{-1,0,1\}$.
In contrast, Theorem \ref{Thm: image 1} only holds for $\alpha: B\rightarrow \{0,1\}$.
Specifically, when $\alpha_i=-1$, the boundary of the admissible space $\Omega(\eta)$ in $[-\infty,+\infty]^N$ will include the case that
there exists $i\in B$ such that $f_i^{(m)}\rightarrow 0^-$.
For this case, we are unable to prove the closedness of $X$.
\end{remark}

\section{Rigidity and existence of the discrete conformal structure (A2)}\label{section A2}

By Remark \ref{Rmk: 1}, we set $\alpha\equiv -1$.
Then (\ref{Eq: new 1}) can be written as
\begin{equation}\label{Eq: F40}
\cosh l_{ij}=\sqrt{(e^{2f_i}-1)(e^{2f_j}-1)}
+\eta_{ij}e^{f_i+f_j}.
\end{equation}
For a right-angled hyperbolic hexagon $\{ijk\}\in F$ with edge lengths $l_{ij}, l_{jk}, l_{ki}$ given by (\ref{Eq: F40}),
the admissible space of the discrete conformal factors $f$ is defined as
\begin{equation}\label{Eq: AS_A2}
\Omega_{ijk}(\eta)=\{(f_i,f_j,f_k)\in \mathbb{R}_{> 0}^3 \mid l_{ij}>0,\  l_{jk}>0,\ l_{ik}>0\}.
\end{equation}
Since $\frac{\partial f_i}{\partial u_i}=\sqrt{e^{2f_i}-1}$,
it follows from Remark \ref{Rmk: 2} (iv) that $\tan u_i=\sqrt{e^{2f_i}-1}$ for $u_i\in (0,\frac{\pi}{2})$ and
$-\frac{1}{\tan u_i}=\sqrt{e^{2f_i}-1}$ for $u_i\in (-\frac{\pi}{2},0)$.
For simplicity, we consider $u\in (-\frac{\pi}{2},0)^N$.
Consequently,
\begin{equation}\label{Eq: F51}
\cosh l_{ij}=\frac{1}{\tan u_i\tan u_j}
+\eta_{ij}\frac{1}{\sin u_i\sin u_j}>1,
\end{equation}
which implies $\cos (u_i+u_j)>-\eta_{ij}$.
This further yields $u_i+u_j\in (-\pi,0)$ when $\eta_{ij}\in [1,+\infty)$, and
$u_i+u_j>-\arccos \eta_{ij}$ when $\eta_{ij}\in (-1,1)$.
Hence, the admissible space $\Omega_{ijk}(\eta)$ is transformed into the following admissible space $\mathcal{U}_{ijk}(\eta)$ of $u$ for a right-angled hyperbolic hexagon $\{ijk\}\in F$:
\begin{equation}\label{Eq: F13}
\mathcal{U}_{ijk}(\eta)=u(\Omega_{ijk}(\eta))
=\{(u_i,u_j,u_k)\in (-\frac{\pi}{2},0)^3 \mid \cos (u_r+u_s)>-\eta_{rs}, \, \forall\{r,s\}\subseteq\{i,j,k\} \},
\end{equation}
where $\eta_{rs}\in (-1,+\infty)$.
Hence, $\mathcal{U}_{ijk}(\eta)$ is a convex polytope.
As a result, the admissible space $\mathcal{U}(\eta)=\bigcap_{\{ijk\}\in F}\mathcal{U}_{ijk}(\eta)$ is also a convex polytope on $(\Sigma,\mathcal{T},\eta)$.

\begin{theorem}\label{Thm: matrix 1}
Suppose $(\Sigma,\mathcal{T},\eta)$ is a weighted triangulated surface with boundary, where $\eta\in (-1,+\infty)^E$.
For a right-angled hyperbolic hexagon $\{ijk\}\in F$ on $\mathcal{U}_{ijk}(\eta)$,
the Jacobian $\Lambda_{ijk}=\frac{\partial(\theta^{jk}_i,\theta^{ik}_j,\theta^{ij}_k)}
{\partial(u_i,u_j,u_k)}$ is symmetric and negative definite.
Consequently, the Jacobian $\Lambda=\frac{\partial (K_i,..., K_N)}{\partial(u_i,...,u_N)}$ is symmetric and negative definite on $\mathcal{U}(\eta)$.
\end{theorem}
\proof
For simplicity, we adopt the notations introduced in Subsection \ref{subsec: matrix}.
Define
\begin{equation*}
\frac{\partial(f_i,f_j,f_k)}
{\partial(u_i,u_j,u_k)}
=\left(
   \begin{array}{ccc}
     \sqrt{e^{2f_i}-1} & 0 & 0 \\
     0 & \sqrt{e^{2f_j}-1} & 0 \\
     0 & 0 & \sqrt{e^{2f_k}-1} \\
   \end{array}
 \right)
 :=Q_4.
\end{equation*}
Then $\Lambda_{ijk}$ can be expressed as
\begin{equation*}
\Lambda_{ijk}=
\frac{\partial(\theta_i,\theta_j,\theta_k)}
{\partial(u_i,u_j,u_k)}
=-\frac{1}{A}
 \left(
   \begin{array}{ccc}
     \sinh l_i & 0 & 0 \\
     0 & \sinh l_j & 0 \\
     0 & 0 & \sinh l_k \\
   \end{array}
 \right)Q_1Q_2Q_4,
\end{equation*}
where $Q_1$ and $Q_2$ are defined by (\ref{Eq: Q_1}) and (\ref{Eq: Q_2}), respectively. 
Furthermore, $\det Q_1>0$.
Since $d>0$ for the discrete conformal structure (A2) in Theorem \ref{Thm: DCS},
it follows from \eqref{Eq: F23} that $\det Q_2>0$.
Consequently, $\det\Lambda_{ijk}<0$,
which implies $\Lambda_{ijk}$ is non-degenerate.
The symmetry of $\Lambda_{ijk}$ follows from Lemma \ref{Lem: angle variations 3}.
Set $a_{rrs}=\sqrt{e^{2f_r}-1}\cdot\coth d_{rs}>0$.
Define
\begin{equation*}
Q:=Q_1Q_2Q_4=\left(
   \begin{array}{ccc}
     -1 & \cosh\theta_k & \cosh\theta_j \\
     \cosh\theta_k & -1 & \cosh\theta_i \\
     \cosh\theta_j & \cosh\theta_i & -1 \\
   \end{array}
 \right)\\
 \times
  \left(
   \begin{array}{ccc}
     0 & a_{jjk} & a_{kkj} \\
     a_{iik} & 0 & a_{kki} \\
     a_{iij} & a_{jji} & 0 \\
   \end{array}
   \right).
\end{equation*}
The remaining part of the proof is analogous to that of Theorem \ref{Thm: matrix negative 1}, and thus is omitted here.

\qed

The following theorem establishes the rigidity of the discrete conformal structure (A2), 
which corresponds to the rigidity part of Theorem \ref{Thm: rigidity and image} (ii). 
As its proof is analogous to that of Theorem \ref{Thm: rigidity 1}, it is omitted here.

\begin{theorem}\label{Thm: rigidity A2}
Suppose $(\Sigma,\mathcal{T},\eta)$ is a weighted triangulated surface with boundary, where $\eta\in (-1,+\infty)^E$.
The discrete conformal factor $f$ is uniquely determined by its generalized combinatorial curvature $K\in \mathbb{R}^N_{>0}$.
%In particular, the map from $f$ to $K$ is a smooth embedding.
\end{theorem}

The following theorem characterizes the image of $K$ for the discrete conformal structure (A2), 
which corresponds to the existence part of Theorem \ref{Thm: rigidity and image} (ii).

\begin{theorem}\label{Thm: image A2}
Suppose $(\Sigma,\mathcal{T},\eta)$ is a weighted triangulated surface with boundary, where $\eta\in (-1,0]^E$.
The image of $K$ is $\mathbb{R}^N_{>0}$.
\end{theorem}

To prove Theorem \ref{Thm: image A2}, we need the following lemma. 
The proof is nearly identical to that of Lemma \ref{Lem: limit a1}, and is thus omitted here.

\begin{lemma}\label{Lem: limit A2}
Under the same assumptions as those in Theorem \ref{Thm: rigidity A2},
for a right-angled hyperbolic hexagon $\{ijk\}\in F$ on $\Omega_{ijk}(\eta)$,
if one of the following conditions is satisfied
\begin{description}
\item[(i)] $\lim f_i=+\infty,\ \lim f_j=+\infty,\ \lim f_k=+\infty$,
\item[(ii)] $\lim f_i=+\infty,\ \lim f_j=+\infty,\ \lim f_k=c$,
\item[(iii)] $\lim f_i=+\infty,\ \lim f_j=a,\ \lim f_k=b$,
\end{description}
where $a,b,c$ are positive constants.
Then $\lim\theta^{jk}_i(f_i,f_j,f_k)=0$.
\end{lemma}

\noindent\textbf{Proof of Theorem \ref{Thm: image A2}:}
Let $X$ be the image of $K$.
The openness of $X$ in $\mathbb{R}^N_{>0}$ is a consequence of the smooth injectivity of the map, 
a property that plays a pivotal role in establishing the rigidity theorem.
To prove the closedness of $X$ in $\mathbb{R}^N_{>0}$, we show that for any sequence $\{f^{(m)}\}\in\Omega(\eta)=\bigcap_{\{ijk\}\in F}\Omega_{ijk}(\eta)$ satisfying $\lim_{m\rightarrow +\infty}K^{(m)}\in \mathbb{R}^N_{>0}$, 
there exists a subsequence,
still denoted by $\{f^{(m)}\}$, such that $\lim_{m\rightarrow +\infty}f^{(m)}\in\Omega(\eta)$.

To derive a contradiction, suppose there exists a subsequence $\{f^{(m)}\}$ such that $\lim_{m\rightarrow +\infty}f^{(m)}$ lies on the boundary of $\Omega(\eta)$.
By (\ref{Eq: AS_A2}), the boundary of the admissible space $\Omega(\eta)$ in $[-\infty,+\infty]^N$ consists of the following three parts:
\begin{description}
\item[(i)]
There exists an ideal edge $\{rs\}\in E$ such that $l_{rs}=0$;
\item[(ii)]
There exists $i\in B$ such that $f_i^{(m)}\rightarrow +\infty$;
\item[(iii)]
There exists $i\in B$ such that $f_i^{(m)}\rightarrow 0^+$.
\end{description}
The proof for the case (i) is identical to that of Theorem \ref{Thm: image 1} and is thus omitted here.
For the case (ii), by Lemma \ref{Lem: limit A2}, each generalized angle incident to the boundary component $i\in B$ converges to 0.
Consequently, $\lim_{m\rightarrow +\infty}K_i^{(m)}=0$.
This contradicts the assumption that $\lim_{m\rightarrow +\infty}K^{(m)}\in \mathbb{R}^N_{>0}$.
For the case (iii), if $f_i^{(m)}\rightarrow 0^+$, then
$\cosh l_{ij}=\sqrt{(e^{2f_i}-1)(e^{2f_j}-1)}
+\eta_{ij}e^{f_i+f_j}\rightarrow \eta_{ij}e^{f_j}\geq 1$.
However, since $\eta_{ij}\in (-1,0]$, it follows that $\eta_{ij}e^{f_j}\leq0$, which is a contradiction.
\qed

\begin{remark}
If $\eta=-\cos \Phi\in (-1,1]^E$ and $e^f=\cosh r$,
then the discrete conformal structure (A2) reduces to the $(-1,-1,1)$-type generalized circle packings introduced by Guo-Luo in \cite{GL2},
where the rigidity of such generalized circle packings was also established. 
Consequently, Theorem \ref{Thm: rigidity A2} generalizes the rigidity result for Guo-Luo's  $(-1,-1,1)$-type generalized circle packings in \cite{GL2}.
Guo-Luo further asserted that the existence holds for $\eta\in (-1,1]^E$.
However, their proof in \cite{GL2} contains a gap when  $\eta\in (0,1)^E$,
because they did not account for the boundary behavior corresponding to $r\rightarrow 0^+$.
As a result, Theorem \ref{Thm: image A2} holds for $\eta\in (-1,0]^E$ rather than $\eta\in (-1,1]^E$.
\end{remark}

\section{Rigidity and existence of the discrete conformal structure (A3)}\label{section A3}

By Remark \ref{Rmk: 1}, we set $C\equiv 0$ in the discrete conformal structure (A3).
In this case, the discrete conformal structure (A3) coincides precisely with the $(-1,-1,0)$-type generalized circle packings introduced by Guo-Luo in \cite{GL2}. 
Guo-Luo \cite{GL2} have established the rigidity and existence results for such structures, which can be stated as follows.

\begin{theorem}(\cite{GL2})\label{Thm: GL-rigidity}
Suppose $(\Sigma,\mathcal{T},\eta)$ is a weighted triangulated surface with boundary, where $\eta\in \mathbb{R}_{>0}^E$.
The discrete conformal factor $f$ is uniquely determined by its generalized combinatorial curvature $K\in \mathbb{R}^N_{>0}$.
Furthermore, the image of $K$ is $\mathbb{R}_{>0}^N$.
\end{theorem}

In this section, we present a concise proof for the rigidity part of Theorem \ref{Thm: GL-rigidity}, 
which simplifies the argument given by Guo-Luo in \cite{GL2}. 
For the proof of the existence part, we refer the reader to \cite{GL2}.

For a right-angled hyperbolic hexagon $\{ijk\}\in F$ with edge lengths $l_{ij}, l_{jk}, l_{ki}$ given by (\ref{Eq: DCS1}),
the admissible space of the discrete conformal factor $f$ is defined as
\begin{equation*}
\Omega_{ijk}(\eta)=\{(f_i,f_j,f_k)\in \mathbb{R}^3 \mid \cosh l_{rs}=-\cosh\,(f_r-f_s)+\eta_{rs}e^{f_r+f_s}>1,\,  \forall\{r,s\}\subseteq\{i,j,k\}\}.
\end{equation*}
By Theorem \ref{Thm: DCS}, for the discrete conformal structure (A3), we have $d>0$.
By the variable transformation in Remark \ref{Rmk: 2}, specifically, $\frac{d f_r}{d u_r}=e^{f_r}$ and $u_r=-e^{-f_r}$, we derive
\begin{equation*}
\eta_{rs}=\frac{\cosh l_{rs}+\cosh(f_r-f_s)}{e^{f_r+f_s}}
>\frac{1+\frac{1}{2}(e^{f_r-f_s}+e^{f_s-f_r})}{e^{f_r+f_s}}
=\frac{1}{2}(e^{-f_r}+e^{-f_s})^2
=\frac{1}{2}(u_r+u_s)^2.
\end{equation*}
Since $u_r\in \mathbb{R}_{<0}$ and $u_s\in \mathbb{R}_{<0}$, 
it follows that $\eta_{rs}>0$ and $u_r+u_s>-\sqrt{2\eta_{rs}}$.
Consequently, for a right-angled hyperbolic hexagon $\{ijk\}\in F$,
the admissible space $\Omega_{ijk}(\eta)$ of $f$ is transformed into the following admissible space $\mathcal{U}_{ijk}(\eta)$ of $u$:
\begin{equation}\label{Eq: F12}
\mathcal{U}_{ijk}(\eta)=u(\Omega_{ijk}(\eta))
=\{(u_i,u_j,u_k)\in \mathbb{R}_{<0}^3 \mid u_r+u_s>-\sqrt{2\eta_{rs}}, \, \forall\{r,s\}\subseteq\{i,j,k\} \},
\end{equation}
which is a convex polytope.
Consequently, the admissible space $\mathcal{U}(\eta)=\bigcap_{\{ijk\}\in F}\mathcal{U}_{ijk}(\eta)$ is also a convex polytope on $(\Sigma,\mathcal{T},\eta)$.

\begin{theorem}\label{Thm: matrix negative 2}
Suppose $(\Sigma,\mathcal{T},\eta)$ is a weighted triangulated surface with boundary, where $\eta\in \mathbb{R}_{>0}^E$.
For a right-angled hyperbolic hexagon $\{ijk\}\in F$ on $\mathcal{U}_{ijk}(\eta)$,
the Jacobian $\Lambda_{ijk}=\frac{\partial(\theta^{jk}_i,\theta^{ik}_j,\theta^{ij}_k)}
{\partial(u_i,u_j,u_k)}$ is symmetric and negative definite.
Consequently, the Jacobian $\Lambda=\frac{\partial (K_i,..., K_N)}{\partial(u_i,...,u_N)}$ is symmetric and negative definite on $\mathcal{U}(\eta)$.
\end{theorem}
\proof
For simplicity, we adopt the notations introduced in Subsection \ref{subsec: matrix}.
Since $\frac{d f_r}{d u_r}=e^{f_r}$, it follows that
\begin{equation*}
\frac{\partial(f_i,f_j,f_k)}
{\partial(u_i,u_j,u_k)}
=\left(
   \begin{array}{ccc}
     e^{f_i} & 0 & 0 \\
     0 & e^{f_j} & 0 \\
     0 & 0 & e^{f_k} \\
   \end{array}
 \right)
 :=Q_5.
\end{equation*}
Then $\Lambda_{ijk}$ can be expressed as
\begin{equation*}
\Lambda_{ijk}=
\frac{\partial(\theta_i,\theta_j,\theta_k)}
{\partial(u_i,u_j,u_k)}
=-\frac{1}{A}
 \left(
   \begin{array}{ccc}
     \sinh l_i & 0 & 0 \\
     0 & \sinh l_j & 0 \\
     0 & 0 & \sinh l_k \\
   \end{array}
 \right)Q_1Q_2Q_5,
\end{equation*}
where $Q_1$ and $Q_2$ are defined by (\ref{Eq: Q_1}) and (\ref{Eq: Q_2}), respectively. 
Furthermore, $\det Q_1>0$.
Since $d>0$ for the discrete conformal structure (A3) in Theorem \ref{Thm: DCS},
it follows from (\ref{Eq: F23}) that $\det Q_2>0$.
Consequently, $\det\Lambda_{ijk}<0$,
which implies $\Lambda_{ijk}$ is non-degenerate.
The symmetry of $\Lambda_{ijk}$ follows from Lemma \ref{Lem: angle variations 3}.
Set $a_{rrs}=e^{f_r}\cdot\coth d_{rs}>0$.
Define
\begin{equation*}
Q:=Q_1Q_2Q_5=\left(
   \begin{array}{ccc}
     -1 & \cosh\theta_k & \cosh\theta_j \\
     \cosh\theta_k & -1 & \cosh\theta_i \\
     \cosh\theta_j & \cosh\theta_i & -1 \\
   \end{array}
 \right)\\
 \times
  \left(
   \begin{array}{ccc}
     0 & a_{jjk} & a_{kkj} \\
     a_{iik} & 0 & a_{kki} \\
     a_{iij} & a_{jji} & 0 \\
   \end{array}
   \right).
\end{equation*}
The remaining part of the proof is analogous to that of Theorem \ref{Thm: matrix negative 1}, and thus is omitted here.

\qed

The following theorem establishes the rigidity of the discrete conformal structure (A3),
which corresponds to the rigidity part of Theorem \ref{Thm: rigidity and image} (iii).
As its proof is similar to that of Theorem \ref{Thm: rigidity 1}, it is omitted here.

\begin{theorem}\label{Thm: rigidity 2}
Suppose $(\Sigma,\mathcal{T},\eta)$ is a weighted triangulated surface with boundary, where $\eta\in \mathbb{R}_{>0}^E$.
The discrete conformal factor $f$ is uniquely determined by its generalized combinatorial curvature $K\in \mathbb{R}^N_{>0}$.
\end{theorem}

\section{Rigidity and existence of the mixed discrete conformal structure III}\label{section III}

In this section, we investigate the rigidity and existence of the mixed discrete conformal structure III, specifically, an ideally triangulated surface with boundary admitting both discrete conformal structures (A3) and (B3). 
All right-angled hyperbolic hexagons are either Type-III or a combination of Type-III right-angled hyperbolic hexagons and those with edge lengths given by (\ref{Eq: DCS1}). 
Here, a right-angled hyperbolic hexagon is Type-III if its edge lengths are given by both (\ref{Eq: DCS1}) and (\ref{Eq: DCS2}). 
By Remark \ref{Rmk: 1}, we set $C\equiv 0$ for both (\ref{Eq: DCS1}) and (\ref{Eq: DCS2}).

As noted in Remark \ref{Rmk: 3},
if $d_{rs}<0$, then $d_{rt}<0,\ d_{sr}>0,\ d_{tr}>0,\ d_{st}>0,\ d_{ts}>0$,
and $l_{rs}, l_{rt}$ are given by (\ref{Eq: DCS2}) and $l_{st}$ is given by (\ref{Eq: DCS1}),
where $\{r,s,t\}=\{i,j,k\}$.
Without loss of generality, we assume $d_{ij}<0$ in this section.
This implies $d_{ik}<0$, $d_{ji}>0$, $d_{ki}>0$, $d_{jk}>0$, and $d_{kj}>0$.
For the Type-III right-angled hyperbolic hexagon $\{ijk\}\in F$,
the edge lengths $l_{ij},l_{ik}$ are given by (\ref{Eq: DCS2}), and $l_{jk}$ is given by (\ref{Eq: DCS1}), namely,
\begin{equation}\label{Eq: F37}
\begin{aligned}
\cosh l_{ij}
&=\cosh(f_j-f_i)+\eta_{ij}e^{f_i+f_j},\\
\cosh l_{ik}
&=\cosh(f_k-f_i)+\eta_{ik}e^{f_i+f_k},\\
\cosh l_{jk}
&=-\cosh(f_j-f_k)+\eta_{jk}e^{f_j+f_k}.
\end{aligned}
\end{equation}
Here $\eta_{ij}\in \mathbb{R}$, $\eta_{ik} \in \mathbb{R}$, and $\eta_{jk}\in \mathbb{R}_{>0}$.

\subsection{Rigidity of the mixed discrete conformal structure III}

For a right-angled hyperbolic hexagon $\{ijk\}\in F$ with edge lengths given by (\ref{Eq: F37}),
the admissible space of the discrete conformal factor $f$ is defined as
\begin{equation*}
\Omega_{ijk}(\eta)=\{(f_i,f_j,f_k)\in \mathbb{R}^3 \mid l_{ij}>0,\  l_{jk}>0,\ l_{ik}>0\}.
\end{equation*}
Combining (\ref{Eq: F26}) and Remark \ref{Rmk: 2} yields $u_i=e^{-f_i}\in \mathbb{R}_{>0}$, $u_j=-e^{-f_j}\in \mathbb{R}_{<0}$, and $u_k=-e^{-f_k}\in \mathbb{R}_{<0}$.
By (\ref{Eq: F37}), we have
\begin{equation*}
\eta_{jk}=\frac{\cosh l_{jk}+\cosh(f_j-f_k)}{e^{f_j+f_k}}
>\frac{1+\frac{1}{2}(e^{f_j-f_k}+e^{f_k-f_j})}{e^{f_j+f_k}}
=\frac{1}{2}(e^{-f_j}+e^{-f_k})^2
=\frac{1}{2}(u_j+u_k)^2.
\end{equation*}
Since $u_j\in \mathbb{R}_{<0}$ and $u_k\in \mathbb{R}_{<0}$, it follows that
$\eta_{jk}>0$ and $u_j+u_k>-\sqrt{2\eta_{jk}}$.
Similarly,
\begin{equation*}
\eta_{ij}=\frac{\cosh l_{ij}-\cosh(f_j-f_i)}{e^{f_i+f_j}}
>\frac{1-\frac{1}{2}(e^{f_j-f_i}+e^{f_i-f_j})}{e^{f_i+f_j}}
=-\frac{1}{2}(e^{-f_i}-e^{-f_j})^2=-\frac{1}{2}(u_i+u_j)^2,
\end{equation*}
which implies $(u_i+u_j)^2>-2\eta_{ij}$.
Moreover, if $\eta_{ij}>0$, then $u_i+u_j\in \mathbb{R}$;
if $\eta_{ij}\leq0$, then $u_i+u_j>\sqrt{-2\eta_{ij}}$ or $u_i+u_j<-\sqrt{-2\eta_{ij}}$.
To exclude the solution $u_i+u_j<-\sqrt{-2\eta_{ij}}$,
we add the condition that $\eta_{jk}+\eta_{ij}\leq0$.
In fact, $u_i+u_j<-\sqrt{-2\eta_{ij}}
\leq-\sqrt{2\eta_{jk}}<u_j+u_k$,
which contradicts the fact $u_i\in \mathbb{R}_{>0}$ and $u_k\in \mathbb{R}_{<0}$.
Similar arguments also apply for $\eta_{ik}$.
Consequently, we obtain the following theorem, which characterizes the admissible space of $u$.

\begin{theorem}\label{Thm: ASC 2}
Let $\{ijk\}\in F$ be a right-angled hyperbolic  hexagon with edge lengths given by (\ref{Eq: F37}).
If one of the following conditions is satisfied
\begin{description}
\item[(i)] $\eta_{jk}>0$, $\eta_{ij}>0$, and $\eta_{ik}>0$;
\item[(ii)] $\eta_{jk}>0$, $\eta_{ij}>0$, $\eta_{ik}\leq0$, and $\eta_{jk}+\eta_{ik}\leq0$;
\item[(iii)] $\eta_{jk}>0$, $\eta_{ij}\leq0$, $\eta_{ik}>0$, and $\eta_{jk}+\eta_{ij}\leq0$;
\item[(iv)] $\eta_{jk}>0$, $\eta_{ij}\leq0$, $\eta_{ik}\leq0$, $\eta_{jk}+\eta_{ij}\leq0$, and $\eta_{jk}+\eta_{ij}\leq0$.
\end{description}
Then the admissible space
\begin{equation}\label{Eq: admissible space c1}
\begin{aligned}
\mathcal{U}_{ijk}(\eta)
&=u(\Omega_{ijk}(\eta))\\
&=\{(u_i,u_j,u_k)\in \mathbb{R}_{>0}\times\mathbb{R}_{<0}^2 \mid u_j+u_k>-\sqrt{2\eta_{jk}},\ \\
&\ \ \ \ \ \ \ \ (u_i+u_j)^2>-2\eta_{ij}, (u_i+u_k)^2>-2\eta_{ik} \}
\end{aligned}
\end{equation}
of $u$ is a convex polytope.
\end{theorem}
\proof
According to the arguments above Theorem \ref{Thm: ASC 2}, we can explicitly characterize the admissible space $\mathcal{U}_{ijk}(\eta)$.
For the case (i), the admissible space (\ref{Eq: admissible space c1}) reduces to
\begin{equation*}
\mathcal{U}_{ijk}(\eta)
=\{(u_i,u_j,u_k)\in \mathbb{R}_{>0}\times\mathbb{R}_{<0}^2 \mid u_j+u_k>-\sqrt{2\eta_{jk}} \}.
\end{equation*}
For the case (ii), the admissible space (\ref{Eq: admissible space c1}) reduces to
\begin{align*}
\mathcal{U}_{ijk}(\eta)
=\{(u_i,u_j,u_k)\in \mathbb{R}_{>0}\times\mathbb{R}_{<0}^2 \mid u_j+u_k>-\sqrt{2\eta_{jk}}, u_i+u_k>\sqrt{-2\eta_{ik}} \}.
\end{align*}
For the case (iii), the admissible space (\ref{Eq: admissible space c1}) reduces to
\begin{align*}
\mathcal{U}_{ijk}(\eta)
=\{(u_i,u_j,u_k)\in \mathbb{R}_{>0}\times\mathbb{R}_{<0}^2 \mid u_j+u_k>-\sqrt{2\eta_{jk}}, u_i+u_j>\sqrt{-2\eta_{ij}} \}.
\end{align*}
For the case (iv), the admissible space (\ref{Eq: admissible space c1}) reduces to
\begin{align*}
\mathcal{U}_{ijk}(\eta)
=\{(u_i,u_j,u_k)\in \mathbb{R}_{>0}\times\mathbb{R}_{<0}^2 \mid &\ u_j+u_k>-\sqrt{2\eta_{jk}}, u_i+u_j>\sqrt{-2\eta_{ij}}, \\ &\ u_i+u_k>\sqrt{-2\eta_{ik}} \}.
\end{align*}
Each of these four cases shows that $\mathcal{U}_{ijk}(\eta)$ is a convex polytope.
\qed

\begin{corollary}\label{Cor: ASC 2}
Suppose $(\Sigma,\mathcal{T},\eta)$ is a weighted triangulated surface with boundary.
If one of the following conditions holds:
\begin{description}
\item[(1)]
All right-angled hyperbolic hexagons are Type-III (i.e., the edge lengths of any right-angled hyperbolic hexagons are given by (\ref{Eq: F37})),
and the weight $\eta$ satisfies one of the conditions (i)--(iv) in Theorem \ref{Thm: ASC 2};

\item[(2)]
Not all right-angled hyperbolic hexagons are of Type-III (i.e., there exists at least one right-angled hyperbolic hexagon with edge lengths given by \eqref{Eq: DCS1}). 
In this case, the weight $\eta$ satisfies one of the conditions (i)--(iv) in Theorem \ref{Thm: ASC 2} for those right-angled hyperbolic hexagons with edge lengths given by (\ref{Eq: F37}),
and $\eta\in \mathbb{R}_{>0}^E$ for those with edge lengths given by (\ref{Eq: DCS1}).
\end{description}
Then the admissible space 
$\mathcal{U}(\eta)=\bigcap_{\{ijk\}\in F}\mathcal{U}_{ijk}(\eta)$ is a convex polytope on $(\Sigma,\mathcal{T},\eta)$,
\end{corollary}
\proof
The case (1) follows from Theorem \ref{Thm: ASC 2}.
The case (2) follows from (\ref{Eq: F12}) and Theorem \ref{Thm: ASC 2}.
Indeed, by (\ref{Eq: F12}), the admissible space $\mathcal{U}_{ijk}(\eta)$ of a right-angled hyperbolic hexagon $\{ijk\}$ with edge lengths given by (\ref{Eq: DCS1}) is a convex polytope.
\qed

\begin{theorem}\label{Thm: matrix negative 3}
Under the same assumptions as those in Theorem \ref{Thm: ASC 2},
the Jacobian $\Lambda_{ijk}=\frac{\partial(\theta^{jk}_i,\theta^{ik}_j,\theta^{ij}_k)}
{\partial(u_i,u_j,u_k)}$ is symmetric and negative definite on $\mathcal{U}_{ijk}(\eta)$.
\end{theorem}
\proof
For simplicity, we adopt the notations introduced in Subsection \ref{subsec: matrix}.
Formula (\ref{Eq: F26}) implies
\begin{equation*}
\frac{\partial(f_i,f_j,f_k)}
{\partial(u_i,u_j,u_k)}
=\left(
   \begin{array}{ccc}
     -e^{f_i} & 0 & 0 \\
     0 & e^{f_j} & 0 \\
     0 & 0 & e^{f_k} \\
   \end{array}
 \right)
 :=Q_6.
\end{equation*}
Then $\Lambda_{ijk}$ can be expressed as
\begin{equation*}
\Lambda_{ijk}=
\frac{\partial(\theta_i,\theta_j,\theta_k)}
{\partial(u_i,u_j,u_k)}
=-\frac{1}{A}
 \left(
   \begin{array}{ccc}
     \sinh l_i & 0 & 0 \\
     0 & \sinh l_j & 0 \\
     0 & 0 & \sinh l_k \\
   \end{array}
 \right)Q_1Q_2Q_6,
\end{equation*}
where $Q_1$ and $Q_2$ are defined by (\ref{Eq: Q_1}) and (\ref{Eq: Q_2}), respectively. 
Furthermore, $\det Q_1>0$.

For any right-angled hyperbolic hexagon $\{ijk\}\in F$ with edge lengths given by (\ref{Eq: F37}),
since $d_{ij}<0$, $d_{ik}<0$, $d_{ik}<0$, $d_{ji}>0$, $d_{ki}>0$, $d_{jk}>0$, and $d_{kj}>0$, it follows that
\begin{align*}
\det Q_2&=\coth d_{jk}\coth d_{ki}\coth d_{ij}+\coth d_{kj}\coth d_{ik}\coth d_{ji}\\
&=\frac{\cosh d_{jk}\cosh d_{ki}\cosh d_{ij}+\cosh d_{kj}\cosh d_{ik}\cosh d_{ji}}{\sinh d_{ij}\sinh d_{jk}\sinh d_{ki}}<0,
\end{align*}
where (\ref{Eq: compatible condition}) is used in the second line.
Consequently,
\begin{equation*}
\det\Lambda_{ijk}
=-\frac{1}{A^3}\sinh l_i\sinh l_j\sinh l_k\cdot\det Q_1
\cdot\det Q_2
\cdot\det Q_6
<0,
\end{equation*}
which implies the matrix $\Lambda_{ijk}$ is non-degenerate.
The symmetry of $\Lambda_{ijk}$ follows from Lemma \ref{Lem: angle variations 3}.
Set $a_{rrs}=e^{f_r}\cdot\coth d_{rs}$.
Then $a_{iij}<0$, $a_{iik}<0$, $a_{jjk}>0$, $a_{kkj}>0$, $a_{jji}>0$, and $a_{kki}>0$.
Define
\begin{equation*}
Q:=Q_1Q_2Q_6
=\left(
   \begin{array}{ccc}
     -1 & \cosh\theta_k & \cosh\theta_j \\
     \cosh\theta_k & -1 & \cosh\theta_i \\
     \cosh\theta_j & \cosh\theta_i & -1 \\
   \end{array}
 \right)
 \times
  \left(
   \begin{array}{ccc}
     0 & a_{jjk} & a_{kkj} \\
     -a_{iik} & 0 & a_{kki} \\
     -a_{iij} & a_{jji} & 0 \\
   \end{array}
   \right).
\end{equation*}
By arguments analogous to those in the proof of Theorem \ref{Thm: matrix negative 1}, 
we show that the $1\times 1$ and $2\times 2$ leading principal minors of $Q$ are positive. 
Let $m_{rs}$ denote the entry of $Q$ at $r$-th row and $s$-th column.
Then
\begin{equation*}
\begin{aligned}
&m_{11}=-\cosh \theta_k\cdot a_{iik}-\cosh \theta_j\cdot a_{iij}\\
&m_{22}=\cosh \theta_k\cdot a_{jjk}+\cosh \theta_i\cdot a_{jji}\\
&m_{12}=-a_{jjk}+\cosh \theta_j\cdot a_{jji}\\
&m_{21}=a_{iik}-\cosh \theta_i\cdot a_{iij}.
\end{aligned}
\end{equation*}
Thus the $1\times 1$ leading principal minor of $Q$ is $m_{11}>0$ due to $a_{iik}<0$ and $a_{iij}<0$.
And the $2\times 2$ leading principal minor of $Q$ is
\begin{equation*}
\begin{aligned}
m_{11}m_{22}-m_{12}m_{21}
=&-[\sinh^2 \theta_ka_{iik}a_{jjk}
+(\cosh \theta_i\cosh\theta_k+\cosh \theta_j)a_{iik}a_{jji}\\
&\ \ \ \ \ +(\cosh \theta_j\cosh\theta_k+\cosh \theta_i)a_{iij}a_{jjk}]
>0.
\end{aligned}
\end{equation*}
This completes the proof.
\qed

\begin{corollary}
Under the same assumptions as those in Corollary \ref{Cor: ASC 2},
the Jacobian $\Lambda=\frac{\partial (K_i,..., K_N)}{\partial(u_i,...,u_N)}$ is symmetric and negative definite on $\mathcal{U}(\eta)$.
\end{corollary}
\proof
The case (1) follows directly from Theorem \ref{Thm: matrix negative 3}.
The proof of the case (2) is analogous to that of Corollary \ref{Cor: matrix 1}, and thus we omit the details here.
Note that $\Lambda=\sum_{\{ijk\}\in F}\Lambda_{ijk}$,
where each $\Lambda_{ijk}$ is derived from either Theorem \ref{Thm: matrix negative 2} or Theorem \ref{Thm: matrix negative 3}.
\qed

The following theorem establishes the rigidity of the mixed discrete conformal structure III,
which generalizes the rigidity part of Theorem \ref{Thm: rigidity and image} (vi).
As its proof is analogous to that of Theorem \ref{Thm: rigidity 1}, it is omitted here.

\begin{theorem}\label{Thm: rigidity 3}
Under the same assumptions as those in Corollary \ref{Cor: ASC 2},
the discrete conformal factor $f$ is uniquely determined by its generalized combinatorial curvature $K\in \mathbb{R}_{>0}^N$.
%In particular, the map from $f$ to $K$ is a smooth embedding.
\end{theorem}

\subsection{Existence of the mixed discrete conformal structure III}

The following theorem characterizes the image of $K$ for the mixed discrete conformal structure III,
which generalizes the existence part of Theorem \ref{Thm: rigidity and image} (vi).

\begin{theorem}\label{Thm: image 3}
Under the same assumptions as Corollary \ref{Cor: ASC 2}, the image of $K$ is $\mathbb{R}^N_{>0}$.
\end{theorem}

To prove Theorem \ref{Thm: image 3}, we need the following three lemmas.

\begin{lemma}(\cite{GL2}, Lemma 4.6)
\label{Lem: GL-limit}
Let $\{ijk\}\in F$ be a right-angled hyperbolic  hexagon with edge lengths given by (\ref{Eq: DCS1}).
Then $\lim_{f_j\rightarrow+\infty}\theta^{ik}_j(f_i,f_j,f_k)=0$, and the converge is uniform.
\end{lemma}

\begin{lemma}\label{Lem: limit c1}
Under the same assumptions as those in Theorem \ref{Thm: ASC 2},
if one of the following conditions is satisfied
\begin{description}
\item[(i)] $\lim f_i=+\infty,\ \lim f_j=+\infty,\ \lim f_k=+\infty$;
\item[(ii)] $\lim f_i=+\infty,\ \lim f_j=+\infty,\ \lim f_k=c_1$;
\item[(iii)] $\lim f_i=+\infty,\ \lim f_j=c_2,\ \lim f_k=c_3$;
\item[(iv)] $\lim f_i=-\infty,\ \lim f_j=+\infty,\ \lim f_k=+\infty$;
\item[(v)] $\lim f_i=-\infty,\ \lim f_j=+\infty,\ \lim f_k=c_4$;
\item[(vi)] $\lim f_i=-\infty,\ \lim f_j=c_5,\ \lim f_k=c_6$,
\end{description}
where $c_1,...,c_6\in \mathbb{R}$ are constants.
Then $\lim\theta^{jk}_i(f_i,f_j,f_k)=0$.
\end{lemma}
\proof
For the case (i), by (\ref{Eq: F37}), we derive
\begin{equation}\label{Eq: F18}
\begin{aligned}
\cosh l_{ij}
&=e^{f_i+f_j}(\eta_{ij}+\frac{1}{2}e^{-2f_i}
+\frac{1}{2}e^{-2f_j}),\\
\cosh l_{ik}
&=e^{f_i+f_k}(\eta_{ik}+\frac{1}{2}e^{-2f_i}
+\frac{1}{2}e^{-2f_k}),\\
\cosh l_{jk}
&=e^{f_j+f_k}(\eta_{jk}-\frac{1}{2}e^{-2f_j}
-\frac{1}{2}e^{-2f_k}).
\end{aligned}
\end{equation}
Since $\lim f_i=+\infty$, $\lim f_j=+\infty$, and $\lim f_k=+\infty$,
it follows that $\lim\cosh l_{ij}:=\lim c_ke^{f_i+f_j}$, $\lim\cosh l_{jk}:=\lim c_ie^{f_j+f_k}$, and $\lim\cosh l_{ik}:=\lim c_je^{f_i+f_k}$, where $c_i, c_j, c_k$ are positive constants. By the cosine law, we have
\begin{equation*}
\cosh \theta^{jk}_i=\frac{\cosh l_{jk}+\cosh l_{ij}\cosh l_{ik}}{\sinh l_{ij}\sinh l_{ik}}\rightarrow
\frac{c_ie^{f_j+f_k}}{c_ke^{f_i+f_j}\cdot c_je^{f_i+f_k}}+1
\rightarrow 1,
\end{equation*}
which implies $\theta^{jk}_i\rightarrow 0$.
The proofs for the cases (ii) and (iii) follow analogously and are omitted here.

For the case (iv), by (\ref{Eq: F37}), we derive
\begin{equation}\label{Eq: F19}
\begin{aligned}
&\cosh l_{ij}=e^{f_j-f_i}(\eta_{ij}e^{2f_i}
+\frac{1}{2}e^{2f_i-2f_j}+\frac{1}{2}),\\
&\cosh l_{ik}
=e^{f_k-f_i}(\eta_{ik}e^{2f_i}
+\frac{1}{2}e^{2f_i-2f_k}+\frac{1}{2}),\\
&\cosh l_{jk}
=e^{f_j+f_k}(\eta_{jk}-\frac{1}{2}e^{-2f_j}
-\frac{1}{2}e^{-2f_k}).
\end{aligned}
\end{equation}
Since $\lim f_i=-\infty$, $\lim f_j=+\infty$, and $\lim f_k=+\infty$,
it follows that $\lim\cosh l_{ij}:=\lim c_ke^{f_j-f_i}$, $\lim\cosh l_{ik}:=\lim c_je^{f_k-f_i}$, and $\lim\cosh l_{jk}:=\lim c_ie^{f_j+f_k}$,
where $c_i, c_j, c_k$ are positive constants.
By the cosine law, we have
\begin{equation*}
\cosh \theta^{jk}_i=\frac{\cosh l_{jk}+\cosh l_{ij}\cosh l_{ik}}{\sinh l_{ij}\sinh l_{ik}}\rightarrow
\frac{c_ie^{f_j+f_k}}{c_ke^{f_j-f_i}\cdot c_je^{f_k-f_i}}+1
\rightarrow 1,
\end{equation*}
which implies $\theta^{jk}_i\rightarrow 0$.
The proofs for the cases (ii) and (iii) follow analogously and are omitted here.
\qed

\begin{lemma}\label{Lem: limit c2}
Under the same assumptions as those in Theorem \ref{Thm: ASC 2},
if one of the following conditions is satisfied
\begin{description}
\item[(i)] $\lim f_j=+\infty,\ \lim f_i=+\infty,\ \lim f_k=+\infty$;
\item[(ii)] $\lim f_j=+\infty,\ \lim f_i=+\infty,\ \lim f_k=c_1$;
\item[(iii)] $\lim f_j=+\infty,\ \lim f_i=-\infty,\ \lim f_k=+\infty$;
\item[(iv)] $\lim f_j=+\infty,\ \lim f_i=-\infty,\ \lim f_k=c_2$;
\item[(v)] $\lim f_j=+\infty,\ \lim f_i=c_3,\ \lim f_k=+\infty$;
\item[(vi)] $\lim f_j=+\infty,\ \lim f_i=c_4,\ \lim f_k=c_5$,
\end{description}
where $c_1,...,c_5\in \mathbb{R}$ are constants.
Then $\lim\theta^{ik}_j(f_i,f_j,f_k)=0$.
\end{lemma}
\proof
For the case (i), it follows from (\ref{Eq: F18}) that $\lim\cosh l_{ij}:=\lim c_ke^{f_i+f_j}$,\ $\lim\cosh l_{jk}:=\lim c_ie^{f_j+f_k}$, and $\lim\cosh l_{ik}:=\lim c_je^{f_i+f_k}$, where $c_i, c_j, c_k$ are positive constants.
By the cosine law, we have
\begin{equation*}
\cosh \theta^{ik}_j=\frac{\cosh l_{ik}+\cosh l_{ij}\cosh l_{jk}}{\sinh l_{ij}\sinh l_{jk}}\rightarrow
\frac{c_je^{f_i+f_k}}{c_ke^{f_i+f_j}\cdot c_ie^{f_i+f_j}}+1
\rightarrow 1,
\end{equation*}
which implies $\theta^{ik}_j\rightarrow 0$.
The proofs for the cases (ii), (v) and (vi) follow analogously and are omitted here.

For the case (iii), it follows from (\ref{Eq: F19}) that $\lim\cosh l_{ij}:=\lim c_ke^{f_j-f_i}$,\ $\lim\cosh l_{ik}:=\lim c_je^{f_k-f_i}$, and $\lim\cosh l_{jk}:=\lim c_ie^{f_j+f_k}$,
where $c_i, c_j, c_k$ are positive constants.
By the cosine law, we have
\begin{equation*}
\cosh \theta^{ik}_j=\frac{\cosh l_{ik}+\cosh l_{ij}\cosh l_{jk}}{\sinh l_{ij}\sinh l_{jk}}\rightarrow
\frac{c_je^{f_k-f_i}}{c_ke^{f_j-f_i}\cdot c_ie^{f_j+f_k}}+1
\rightarrow 1,
\end{equation*}
which implies $\theta^{ik}_j\rightarrow 0$.
The proof for the case (iv) follow analogously and is omitted here.

\qed

\noindent\textbf{Proof of Theorem \ref{Thm: image 3}:}
Let $X$ be the image of $K$.
The openness of $X$ in $\mathbb{R}^N_{>0}$ follows directly from the smooth injectivity of the map.
To establish the closedness of $X$ in $\mathbb{R}^N_{>0}$, 
we prove that every sequence $\{f^{(m)}\}\in\Omega(\eta)=\bigcap_{\{ijk\}\in F}\Omega_{ijk}(\eta)$ with $\lim_{m\rightarrow +\infty}K^{(m)}\in \mathbb{R}^N_{>0}$ admits a subsequence, still denoted by $\{f^{(m)}\}$, such that $\lim_{m\rightarrow +\infty}f^{(m)}\in\Omega(\eta)$.

Assume, for the sake of contradiction, that there exists a subsequence $\{f^{(m)}\}$ such that $\lim_{m\rightarrow +\infty}f^{(m)}$ lies on the boundary of $\Omega(\eta)$.
The boundary of the admissible space $\Omega(\eta)$ in $[-\infty,+\infty]^N$ consists of the following three parts:
\begin{description}
\item[(i)]
There exists an ideal edge $\{rs\}\in E$ such that $l_{rs}=0$;
\item[(ii)]
There exists $i\in B$ in a right-angled hyperbolic hexagon $\{ijk\}\in F$ with edge lengths given by (\ref{Eq: F37}), such that $f_i^{(m)}\rightarrow \pm\infty$;
\item[(iii)]
There exists $j\in B$ in a right-angled hyperbolic hexagon $\{ijk\}\in F$ with edge lengths given by (\ref{Eq: F37}) or (\ref{Eq: DCS1}), such that $f_j^{(m)}\rightarrow +\infty$.
\end{description}
The proof for the case (i) is identical to that of Theorem \ref{Thm: image 1} and is thus omitted. 
For the case (ii), by Lemma \ref{Lem: limit c1}, each generalized angle incident to the boundary component $i\in B$ converges to 0.
Consequently, $\lim_{m\rightarrow +\infty}K_i^{(m)}=0$,
which contradicts the assumption that $\lim_{m\rightarrow +\infty}K^{(m)}\in \mathbb{R}^N_{>0}$.
For the case (iii),
if all right-angled hyperbolic hexagons are Type-III,
then Lemma \ref{Lem: limit c2} implies $\lim_{m\rightarrow +\infty}K_j^{(m)}=0$;
if not all right-angled hyperbolic hexagons are Type-III,
then by Lemma \ref{Lem: GL-limit} and Lemma \ref{Lem: limit c2}, we have $\lim_{m\rightarrow +\infty}K_j^{(m)}=0$.
Both scenarios contradict the assumption that $\lim_{m\rightarrow +\infty}K^{(m)}\in \mathbb{R}^N_{>0}$.
\qed

\begin{remark}
Note that in a right-angled hyperbolic hexagon $\{ijk\}\in F$ with edge lengths given by (\ref{Eq: F37}),
there is no $j\in B$ such that $f_j^{(m)}\rightarrow -\infty$.
Suppose, to the contrary, that such a $j$ exists, 
then $u_j=-e^{-f_j}\rightarrow -\infty$.
By the definition of the admissible space $\mathcal{U}_{ijk}(\eta)$ in (\ref{Eq: admissible space c1}),
the inequality $u_j+u_k>-2\sqrt{\eta_{jk}}$ implies $u_k\rightarrow +\infty$, which contradicts the fact that $u_k\in \mathbb{R}_{<0}$.
An analogous argument applies to a right-angled hyperbolic hexagon $\{ijk\}\in F$ with edge lengths given by (\ref{Eq: DCS1}).
\end{remark}

\section{Rigidity of the mixed discrete conformal structure II}\label{section II}

In this section, we investigate the rigidity of the mixed discrete conformal structure II, specifically,
an ideally triangulated surface with boundary admitting both discrete conformal structures (A2) and (B2). 
All right-angled hyperbolic hexagons are either Type-II or a combination of Type-II right-angled hyperbolic hexagons and those with edge lengths given by (\ref{Eq: new 1}).
Here, a right-angled hyperbolic hexagon is Type-II if its edge lengths are given by both (\ref{Eq: new 1}) and (\ref{Eq: new 2}).
By Remark \ref{Rmk: 1}, we have $\alpha\equiv -1$ for both (\ref{Eq: new 1}) and (\ref{Eq: new 2}).

As noted in Remark \ref{Rmk: 3},
if $d_{rs}<0$, then $d_{rt}<0,\ d_{sr}>0,\ d_{tr}>0,\ d_{st}>0,\ d_{ts}>0$,
and $l_{rs}, l_{rt}$ are given by (\ref{Eq: new 2}) and $l_{st}$ is given by (\ref{Eq: new 1}),
where $\{r,s,t\}=\{i,j,k\}$.
Without loss of generality,
we assume $d_{ij}<0$ in this section.
This implies $d_{ik}<0$, $d_{ji}>0$, $d_{ki}>0$, $d_{jk}>0$, and $d_{kj}>0$.
For the Type-II right-angled hyperbolic hexagon $\{ijk\}\in F$, the edge lengths $l_{ij},l_{ik}$ are given by (\ref{Eq: new 2}), and $l_{jk}$ is given by (\ref{Eq: new 1}), namely,
\begin{equation}\label{Eq: F39}
\begin{aligned}
\cosh l_{ij}
&=-\sqrt{(e^{2f_i}-1)(e^{2f_j}-1)}
+\eta_{ij}e^{f_i+f_j},\\
\cosh l_{ik}
&=-\sqrt{(e^{2f_i}-1)(e^{2f_k}-1)}+\eta_{ik}e^{f_i+f_k},\\
\cosh l_{jk}
&=\sqrt{(e^{2f_j}-1)(e^{2f_k}-1)}+\eta_{jk}e^{f_j+f_k}.
\end{aligned}
\end{equation}
Here $\eta_{ij}\in \mathbb{R}_{>0}$, $\eta_{ik}\in \mathbb{R}_{>0}$, and $\eta_{jk}\in (-1,+\infty)$.

For a right-angled hyperbolic hexagon $\{ijk\}\in F$ with edge lengths given by (\ref{Eq: F39}), 
the admissible space of the discrete conformal factor $f$ is defined as
\begin{equation*}
\Omega_{ijk}(\eta)=\{(f_i,f_j,f_k)\in \mathbb{R}_{>0}^3 \mid l_{ij}>0,\  l_{jk}>0,\ l_{ik}>0\}.
\end{equation*}
By combining (\ref{Eq: F50}) and Remark \ref{Rmk: 2} (iv), 
we consider $u\in (-\frac{\pi}{2},0)^N$ for simplicity.
It follows that
$\cosh l_{ij}=-(-\tan u_i)\frac{-1}{\tan u_j}+\eta_{ij}\frac{1}{\cos u_i}\frac{-1}{\sin u_j}>1$,
which implies $\sin(u_i+u_j)>-\eta_{ij}$.
Consequently, $u_i+u_j\in (-\pi,0)$ for $\eta_{ij}\in [1,+\infty)$, and $u_i+u_j>-\arcsin \eta_{ij}$ for $\eta_{ij}\in [0, 1)$.
Moreover, from (\ref{Eq: F51}),
$\cosh l_{jk}>1$ implies $u_j+u_k\in (-\pi,0)$ for $\eta_{jk}\in [1,+\infty)$ and
$u_j+u_k>-\arccos \eta_{jk}$ for $\eta_{jk}\in (-1,1)$.

\begin{theorem}\label{Thm: ASC new}
Let $\{ijk\}\in F$ be a right-angled hyperbolic  hexagon with edge lengths given by (\ref{Eq: F39}),
and the weight $\eta_{ij}\in [0,+\infty)$, $\eta_{ik}\in [0,+\infty)$, $\eta_{jk}\in (-1,+\infty)$.
Then the admissible space
\begin{equation*}
\begin{aligned}
\mathcal{U}_{ijk}(\eta)=u(\Omega_{ijk}(\eta))
=\{&(u_i,u_j,u_k)\in (-\frac{\pi}{2},0)^3 \mid \sin(u_i+u_j)>-\eta_{ij}, \\
&\sin(u_i+u_k)>-\eta_{ik},
\cos(u_j+u_k)>-\eta_{jk}\}.
\end{aligned}
\end{equation*}
of $u$ is a convex polytope.
\end{theorem}

\begin{corollary}\label{Cor: ASC new}
Suppose $(\Sigma,\mathcal{T},\eta)$ is a weighted triangulated surface with boundary.
If one of the following conditions holds: 
\begin{description}
\item[(1)]
All right-angled hyperbolic hexagons are Type-II (i.e., the edge lengths of any right-angled hyperbolic hexagons are given by (\ref{Eq: F39})),
and the weight $\eta$ satisfies the condition in Theorem \ref{Thm: ASC new};
\item[(2)]
Not all right-angled hyperbolic hexagons are Type-II (i.e., there exist some right-angled hyperbolic hexagons with edge lengths given by (\ref{Eq: new 1})).
In this case, the weight $\eta$ satisfies the condition in Theorem \ref{Thm: ASC new} for those right-angled hyperbolic hexagons with edge lengths given by (\ref{Eq: F39}),
and $\eta\in (-1,+\infty)^E$ for those with edge lengths given by (\ref{Eq: new 1}).
\end{description}
Then the admissible space
$\mathcal{U}(\eta)=\bigcap_{\{ijk\}\in F}\mathcal{U}_{ijk}(\eta)$ is a convex polytope on $(\Sigma,\mathcal{T},\eta)$.
\end{corollary}
\proof
The case (1) follows from Theorem \ref{Thm: ASC new}, and the case (2) follows from (\ref{Eq: F13}) and Theorem \ref{Thm: ASC new}.
\qed

\begin{theorem}\label{Thm: matrix 2}
Under the same assumptions as those in Theorem \ref{Thm: ASC new},
the Jacobian $\Lambda_{ijk}=\frac{\partial(\theta^{jk}_i,\theta^{ik}_j,\theta^{ij}_k)}
{\partial(u_i,u_j,u_k)}$ is symmetric and negative definite on $\mathcal{U}_{ijk}(\eta)$.
\end{theorem}
\proof
For simplicity, we adopt the notations introduced in Subsection \ref{subsec: matrix}.
Define
\begin{equation*}
\frac{\partial(f_i,f_j,f_k)}
{\partial(u_i,u_j,u_k)}
=\left(
   \begin{array}{ccc}
     -\sqrt{e^{2f_i}-1} & 0 & 0 \\
     0 & \sqrt{e^{2f_j}-1} & 0 \\
     0 & 0 & \sqrt{e^{2f_k}-1} \\
   \end{array}
 \right)
 :=Q_7.
\end{equation*}
Then $\Lambda_{ijk}$ can be expressed as
\begin{equation*}
\Lambda_{ijk}=
\frac{\partial(\theta_i,\theta_j,\theta_k)}
{\partial(u_i,u_j,u_k)}
=-\frac{1}{A}
 \left(
   \begin{array}{ccc}
     \sinh l_i & 0 & 0 \\
     0 & \sinh l_j & 0 \\
     0 & 0 & \sinh l_k \\
   \end{array}
 \right)Q_1Q_2Q_7,
\end{equation*}
where $Q_1$ and $Q_2$ are defined by (\ref{Eq: Q_1}) and (\ref{Eq: Q_2}), respectively.
Furthermore, $\det Q_1>0$.
For any right-angled hyperbolic hexagon $\{ijk\}\in F$ with edge lengths given by (\ref{Eq: F39}),
since $d_{ij}<0$, $d_{ik}<0$, $d_{ji}>0$, $d_{ki}>0$, $d_{jk}>0$, and $d_{kj}>0$, it follows from (\ref{Eq: F23}) that $\det Q_2<0$.
The remaining part of the proof is analogously to that of Theorem \ref{Thm: matrix negative 3}, and is thus omitted here.

\qed

Combining Theorem \ref{Thm: matrix 1} and Theorem \ref{Thm: matrix 2},
we obtain the following result.
\begin{corollary}
Under the same assumptions as those in Corollary \ref{Cor: ASC new},
the Jacobian $\Lambda=\frac{\partial (K_i,..., K_N)}{\partial(u_i,...,u_N)}$ is symmetric and negative definite on $\mathcal{U}(\eta)$.
\end{corollary}

The following theorem establishes the rigidity of the mixed discrete conformal structure II,
which generalizes Theorem \ref{Thm: rigidity and image} (v).
As its proof is analogous to that of Theorem \ref{Thm: rigidity 1}, it is omitted here.

\begin{theorem}\label{Thm: Rigidity II}
Under the same assumptions as those in Corollary \ref{Cor: ASC new},
the discrete conformal factor $f$ is uniquely determined by its generalized combinatorial curvature $K\in \mathbb{R}^N_{>0}$.
\end{theorem}

\begin{remark}\label{Rmk: 6}
The existence of the mixed discrete conformal structure II cannot be established. 
Specifically, the boundary of the admissible space $\Omega(\eta)=\bigcap_{\{ijk\}\in F}\Omega_{ijk}(\eta)$ in $[-\infty,+\infty]^N$ include the case that
there exists $i\in B$ such that $f_i^{(m)}\rightarrow 0^+$.
For this case, we are unable to prove the closedness of the curvature map.
\end{remark}

\section{Rigidity and existence of the mixed discrete conformal structure I}\label{section I}

In this section, we investigate the rigidity and existence of the mixed discrete conformal structure I, specifically, an ideally triangulated surface with boundary admitting both discrete conformal structures (A1) and (B1). 
All right-angled hyperbolic hexagons are either Type-I or a combination of Type-I right-angled hyperbolic hexagons and those with edge lengths given by (\ref{Eq: DCS3}).
Here, a right-angled hyperbolic hexagon is Type-I if its edge lengths are given by both (\ref{Eq: DCS3}) and (\ref{Eq: DCS4}).

As noted in Remark \ref{Rmk: 3},
if $d_{rs}<0$, then $d_{rt}<0,\ d_{sr}>0,\ d_{tr}>0,\ d_{st}>0,\ d_{ts}>0$,
and $l_{rs}, l_{rt}$ are given by (\ref{Eq: DCS4}) and $l_{st}$ is given by (\ref{Eq: DCS3}),
where $\{r,s,t\}=\{i,j,k\}$.
Without loss of generality, we assume $d_{ij}<0$ in this section.
This implies $d_{ik}<0$, $d_{ji}>0$, $d_{ki}>0$, $d_{jk}>0$, and $d_{kj}>0$.
For the Type-I right-angled hyperbolic hexagon $\{ijk\}\in F$, the edge lengths $l_{ij},l_{ik}$ are given by (\ref{Eq: DCS4}), and $l_{jk}$ is given by (\ref{Eq: DCS3}), namely,
\begin{equation}\label{Eq: F38}
\begin{aligned}
\cosh l_{ij}
&=\sqrt{(1+\alpha_ie^{2f_i})(1+\alpha_je^{2f_j})}
+\eta_{ij}e^{f_i+f_j},\\
\cosh l_{ik}
&=\sqrt{(1+\alpha_ie^{2f_i})(1+\alpha_ke^{2f_k})}
+\eta_{ik}e^{f_i+f_k},\\
\cosh l_{jk}
&=-\sqrt{(1+\alpha_je^{2f_j})(1+\alpha_ke^{2f_k})}
+\eta_{jk}e^{f_j+f_k}.
\end{aligned}
\end{equation}

\subsection{Admissible space of the mixed discrete conformal structure I}\label{subsection: AS6}

For a right-angled hyperbolic  hexagon $\{ijk\}\in F$ with edge lengths given by (\ref{Eq: F38}),
the admissible space of the discrete conformal factor $f$ is defined as
\begin{equation*}
\Omega_{ijk}(\eta)=\{(f_i,f_j,f_k)\in \mathbb{R}^3 \mid l_{ij}>0,\  l_{jk}>0,\ l_{ik}>0\}.
\end{equation*}
Since $\frac{\sinh d_{ij}}{\sinh d_{ji}}
=-\sqrt{\frac{1+\alpha_ie^{2f_i}}{1+\alpha_je^{2f_j}}}$,
$\frac{\sinh d_{jk}}{\sinh d_{kj}}
=\sqrt{\frac{1+\alpha_je^{2f_j}}{1+\alpha_ke^{2f_k}}}$, and
$\frac{\sinh d_{ik}}{\sinh d_{ki}}
=-\sqrt{\frac{1+\alpha_ie^{2f_i}}{1+\alpha_ke^{2f_k}}}$,
it follows from the proof of Lemma \ref{Lem: angle variations 3} that 
\begin{equation}\label{Eq: F27}
\frac{\partial f_i}{\partial u_i}=-\sqrt{1+\alpha_ie^{2f_i}},\
\frac{\partial f_j}{\partial u_j}=\sqrt{1+\alpha_je^{2f_j}},\
\frac{\partial f_k}{\partial u_k}=\sqrt{1+\alpha_ke^{2f_k}}.
\end{equation}

To characterizes the admissible space $\mathcal{U}_{ijk}(\eta)$ of the discrete conformal factors $u$,
it is necessary to derive the conditions under which the edge lengths $l_{ij}>0$ and $l_{jk}>0$.
Since the case that $l_{jk}>0$ has already been addressed in Subsection \ref{section a4}, 
we focus on analyzing $l_{ij}>0$ here.
By Remark \ref{Rmk: 2}, the following results are obtained:
\begin{description}
\item[(i)] If $\alpha_i=0$ and $\alpha_j=0$, then $\cosh l_{ij}=1+\eta_{ij}e^{u_j-u_i}>1$. 
    This implies
    $\eta_{ij}>0$ and $(u_j-u_i)\in \mathbb{R}$.
\item[(ii)] If $\alpha_i=1$ and $\alpha_j=0$, then $\cosh l_{ij}=\sqrt{1+e^{2f_i}}+\eta_{ij}e^{f_i+f_j}=\coth u_i+\eta_{ij}e^{u_j}\frac{1}{\sinh u_i}>1$.
    Since $u_i>0$, it follows that $\eta_{ij}e^{u_i+u_j}+1>0$. 
    Consequently, if $\eta_{ij}\geq0$, then $(u_i+u_j)\in \mathbb{R}$; if $\eta_{ij}<0$, then $u_i+u_j<\log(-\frac{1}{\eta_{ij}})$.
    For simplicity, we only consider the case $\eta_{ij}\geq 0$ in the following.
\item[(iii)] If $\alpha_i=-1$ and $\alpha_j=0$, then $\cosh l_{ij}=\sqrt{1-e^{2f_i}}+\eta_{ij}e^{f_i+f_j}=\tanh u_i+\eta_{ij}e^{u_j}\frac{1}{\cosh u_i}>1$. 
    This implies $\eta_{ij}e^{u_i+u_j}>1$, so $\eta_{ij}>0$ and $u_i+u_j>\log\frac{1}{\eta_{ij}}$.
\item[(iv)] If $\alpha_i=0$ and $\alpha_j=1$, then $\cosh l_{ij}=\sqrt{1+e^{2f_j}}+\eta_{ij}e^{f_i+f_j}=-\coth u_j+\eta_{ij}e^{-u_i}\frac{-1}{\sinh u_j}>1$.
    Since $u_j<0$, it follows that $e^{u_i+u_j}+\eta_{ij}>0$.
    Consequently, if $\eta_{ij}\geq0$, then $(u_i+u_j)\in \mathbb{R}$; if $\eta_{ij}<0$, then $u_i+u_j>\log(-\eta_{ij})$.
\item[(v)] If $\alpha_i=1$ and $\alpha_j=1$, then $\cosh l_{ij}=\sqrt{(1+e^{2f_i})(1+e^{2f_j})}+\eta_{ij}e^{f_i+f_j}=-\coth u_i\coth u_j+\eta_{ij}\frac{1}{\sinh u_i}\frac{-1}{\sinh u_j}>1$.
    Since $u_i>0$ and $u_j<0$, it follows that $\cosh(u_i+u_j)+\eta_{ij}>0$. 
    Consequently, if $\eta_{ij}> -1$, then $(u_i+u_j)\in \mathbb{R}$; if $\eta_{ij}<-1$, then $u_i+u_j>\mathrm{arccosh}\, (-\eta_{ij})$ or $u_i+u_j<-\mathrm{arccosh}\, (-\eta_{ij})$. 
    For simplicity, we only consider the case $\eta_{ij}>-1$ in the following.
\item[(vi)] If $\alpha_i=-1$ and $\alpha_j=1$, then $\cosh l_{ij}=\sqrt{(1-e^{2f_i})(1+e^{2f_j})}+\eta_{ij}e^{f_i+f_j}=-\tanh u_i\coth u_j+\eta_{ij}\frac{1}{\cosh u_i}\frac{-1}{\sinh u_j}>1$.
    Since $u_j<0$, it follows that $\sinh(u_i+u_j)+\eta_{ij}>0$.
    Thus $\eta_{ij}\in \mathbb{R}$ and $u_i+u_j>\mathrm{arcsinh}\ (-\eta_{ij})$.
\item[(vii)] If $\alpha_i=0$ and $\alpha_j=-1$, then $\cosh l_{ij}=\sqrt{1-e^{2f_j}}+\eta_{ij}e^{f_i+f_j}=-\tanh u_j+\eta_{ij}e^{-u_i}\frac{1}{\cosh u_j}>1$. 
    This implies $e^{u_i+u_j}<\eta_{ij}$, so $\eta_{ij}>0$ and $u_i+u_j<\log \eta_{ij}$.
\item[(viii)] If $\alpha_i=1$ and $\alpha_j=-1$, then $\cosh l_{ij}=\sqrt{(1+e^{2f_i})(1-e^{2f_j})}+\eta_{ij}e^{f_i+f_j}
    =-\coth u_i\tanh u_j+\eta_{ij}\frac{1}{\sinh u_i}\frac{1}{\cosh u_j}>1$.
    This implies $\sinh(u_i+u_j)<\eta_{ij}$, so $\eta_{ij}\in \mathbb{R}$ and $u_i+u_j<\mathrm{arcsinh}\, \eta_{ij}$.
\item[(ix)] If $\alpha_i=-1$ and $\alpha_j=-1$, then $\cosh l_{ij}=\sqrt{(1-e^{2f_i})(1-e^{2f_j})}+\eta_{ij}e^{f_i+f_j}
    =-\tanh u_i\tanh u_j+\eta_{ij}\frac{1}{\cosh u_i}\frac{1}{\cosh u_j}>1$.
    This implies $\cosh(u_i+u_j)<\eta_{ij}$, so $\eta_{ij}>1$ and $-\mathrm{arccosh}\,  \eta_{ij}<u_i+u_j<\mathrm{arccosh}\, \eta_{ij}$.
\end{description}

For simplicity, we adopt the triple notation $(\alpha_i,\alpha_j,\alpha_k)$ to classify the distinct types of admissible spaces of $u$ for a right-angled hyperbolic hexagon $\{ijk\}\in F$ with edge lengths given by (\ref{Eq: F38}).
Through permutations and combinations, the admissible space $\Omega_{ijk}(\eta)$ of $f$ can be transformed into 27 distinct types of admissible spaces of $u$.
However, due to the interchangeability of the positions of $j$ and $k$ in (\ref{Eq: F38}),
the triples $(\alpha_i,\alpha_j,\alpha_k)$ and $(\alpha_i,\alpha_k,\alpha_j)$ are considered equivalent and thus belong to the same type.
This reduction results in a total of 18 distinct types of admissible spaces of $u$, denoted by $\mathcal{U}^{(n)}_{ijk}(\eta)$, where $n=1,2,...,18$.
These types are enumerated as follows.

\begin{description}
\item[(I)]
Assume $(\alpha_i,\alpha_j,\alpha_k)=(0,0,0)$.
In Subsection \ref{section a4}, it follows that $u_j+u_k>\log \frac{2}{\eta_{jk}}$ and $\eta_{jk}>0$.
Building upon the preceding analysis, we deduce $\eta_{ij}>0$ with $(u_j-u_i)\in \mathbb{R}$, and $\eta_{ik}>0$ with $(u_k-u_i)\in \mathbb{R}$.
Consequently, the admissible space of $u$ is defined as
\begin{equation}\label{Eq: F4}
\mathcal{U}^{(1)}_{ijk}(\eta)
=\{(u_i,u_j,u_k)\in \mathbb{R}^3\mid u_j+u_k>\log \frac{2}{\eta_{jk}}\},
\end{equation}
where $\eta_{jk}>0$, $\eta_{ij}>0$, and $\eta_{ik}>0$.

\item[(II)]
Assume $(\alpha_i,\alpha_j,\alpha_k)=(1,0,0)$.
It is known that $u_j+u_k>\log \frac{2}{\eta_{jk}}$ and $\eta_{jk}>0$.
Moreover, $(u_i+u_j)\in \mathbb{R}$ for $\eta_{ij}\geq0$, and $(u_i+u_k)\in \mathbb{R}$ for $\eta_{ik}\geq0$.
Consequently, the admissible space of $u$ is defined as
\begin{equation}\label{Eq: F5}
\begin{aligned}
\mathcal{U}^{(2)}_{ijk}(\eta)
=\{&(u_i,u_j,u_k)\in \mathbb{R}_{>0}\times\mathbb{R}^2 \mid u_j+u_k>\log \frac{2}{\eta_{jk}}\},
\end{aligned}
\end{equation}
where $\eta_{jk}>0$, $\eta_{ij}\geq 0$, and $\eta_{ik}\geq0$.

\item[(III)]
Assume $(\alpha_i,\alpha_j,\alpha_k)=(-1,0,0)$.
It is known that $u_j+u_k>\log \frac{2}{\eta_{jk}}$ and $\eta_{jk}>0$.
Moreover, $\eta_{ij}>0$ with $u_i+u_j>\log\frac{1}{\eta_{ij}}$, and
$\eta_{ik}>0$ with $u_i+u_k>\log\frac{1}{\eta_{ik}}$.
Consequently, the admissible space of $u$ is defined as
\begin{equation*}
\mathcal{U}^{(3)}_{ijk}(\eta)
=\{(u_i,u_j,u_k)\in \mathbb{R}_{>0}\times\mathbb{R}^2\mid u_j+u_k>\log \frac{2}{\eta_{jk}},\ u_i+u_j>\log\frac{1}{\eta_{ij}},\ u_i+u_k>\log\frac{1}{\eta_{ik}}\},
\end{equation*}
where $\eta_{jk}>0$, $\eta_{ij}>0$, and $\eta_{ik}>0$.

\item[(IV)]
Assume $(\alpha_i,\alpha_j,\alpha_k)=(0,0,1)$.
In Subsection \ref{section a4}, it follows that $u_j+u_k>\log \frac{1}{\eta_{jk}}$ and $\eta_{jk}>0$.
Moreover, $\eta_{ij}>0$ and $(u_j-u_i)\in \mathbb{R}$.
Similarly, $(u_i+u_k)\in \mathbb{R}$ for $\eta_{ik}\geq0$, and $u_i+u_k>\log(-\eta_{ik})$ for $\eta_{ik}<0$.
Consequently, the admissible space of $u$ is defined as
\begin{equation}\label{AS4}
\mathcal{U}^{(4)}_{ijk}(\eta)
=\{(u_i,u_j,u_k)\in \mathbb{R}^2\times\mathbb{R}_{<0}\mid u_j+u_k>\log \frac{1}{\eta_{jk}},\  e^{u_i+u_k}>-\eta_{ik}\},
\end{equation}
where $\eta_{jk}>0$, $\eta_{ij}>0$, and $\eta_{ik}\in \mathbb{R}$.

The admissible space (\ref{AS4}) can be explicitly characterized as follows:
\begin{description}
\item[(i)]For $\eta_{jk}>0$, $\eta_{ij}>0$, and $\eta_{ik}\geq0$, the admissible space is given by
\begin{equation}\label{Eq: F48}
\mathcal{U}^{(4)}_{ijk}(\eta)
=\{(u_i,u_j,u_k)\in \mathbb{R}^2\times\mathbb{R}_{<0}\mid u_j+u_k>\log \frac{1}{\eta_{jk}}\}.
\end{equation}
\item[(ii)]For $\eta_{jk}>0$, $\eta_{ij}>0$, and $\eta_{ik}<0$, the admissible space is given by
\begin{equation}\label{Eq: F42}
\mathcal{U}^{(4)}_{ijk}(\eta)
=\{(u_i,u_j,u_k)\in \mathbb{R}^2\times\mathbb{R}_{<0}\mid u_j+u_k>\log \frac{1}{\eta_{jk}},\  u_i+u_k>\log(-\eta_{ik})\}.
\end{equation}
\end{description}

\item[(V)]
Assume $(\alpha_i,\alpha_j,\alpha_k)=(1,0,1)$.
It is known that $u_j+u_k>\log \frac{1}{\eta_{jk}}$ and $\eta_{jk}>0$.
Moreover, $(u_i+u_j)\in \mathbb{R}$ for $\eta_{ij}\geq0$, and $(u_i+u_k)\in \mathbb{R}$ for $\eta_{ik}> -1$.
Consequently, the admissible space of $u$ is defined as
\begin{equation}\label{AS5}
\mathcal{U}^{(5)}_{ijk}(\eta)
=\{(u_i,u_j,u_k)\in \mathbb{R}\times\mathbb{R}_{>0}\times\mathbb{R}_{<0}
\mid u_j+u_k>\log \frac{1}{\eta_{jk}}\},
\end{equation}
where $\eta_{jk}>0$, $\eta_{ij}\geq0$, and $\eta_{ik}>-1$.

\item[(VI)]
Assume $(\alpha_i,\alpha_j,\alpha_k)=(-1,0,1)$.
It is known that $u_j+u_k>\log \frac{1}{\eta_{jk}}$ and $\eta_{jk}>0$.
Moreover, $\eta_{ij}>0$ and $u_i+u_j>\log\frac{1}{\eta_{ij}}$. 
Similarly, $\eta_{ik}\in \mathbb{R}$ and $u_i+u_k>\mathrm{arcsinh}\, (-\eta_{ik})$.
Consequently, the admissible space of $u$ is defined as
\begin{equation*}
\begin{aligned}
\mathcal{U}^{(6)}_{ijk}(\eta)
=\{(u_i,u_j,u_k)\in \mathbb{R}_{>0}\times\mathbb{R}\times\mathbb{R}_{<0}
\mid &\ u_j+u_k>\log \frac{1}{\eta_{jk}},\
u_i+u_j>\log\frac{1}{\eta_{ij}}, \\ &\ u_i+u_k>\mathrm{arcsinh}\, (-\eta_{ik}) \},
\end{aligned}
\end{equation*}
where $\eta_{jk}>0$, $\eta_{ij}>0$, and $\eta_{ik}\in \mathbb{R}$.

\item[(VII)]
Assume $(\alpha_i,\alpha_j,\alpha_k)=(0,0,-1)$.
In Subsection \ref{section a4}, it follows that $u_j+u_k>\log \frac{1}{\eta_{jk}}$ and $\eta_{jk}>0$.
Moreover, $\eta_{ij}>0$ and $(u_j-u_i)\in \mathbb{R}$. 
Similarly, $\eta_{ik}>0$ and $u_i+u_k<\log \eta_{ik}$.
Consequently, the admissible space of $u$ is defined as
\begin{equation*}
\mathcal{U}^{(7)}_{ijk}(\eta)
=\{(u_i,u_j,u_k)\in \mathbb{R}^2\times\mathbb{R}_{<0}
\mid u_j+u_k>\log \frac{1}{\eta_{jk}},\ u_i+u_k<\log \eta_{ik}\},
\end{equation*}
where $\eta_{jk}>0$, $\eta_{ij}>0$, and $\eta_{ik}>0$.

\item[(VIII)]
Assume $(\alpha_i,\alpha_j,\alpha_k)=(1,0,-1)$.
It is known that $u_j+u_k>\log \frac{1}{\eta_{jk}}$ and $\eta_{jk}>0$.
Moreover, $(u_i+u_j)\in \mathbb{R}$ for $\eta_{ij}\geq0$.
Similarly, $\eta_{ik}\in \mathbb{R}$ and $u_i+u_k<\mathrm{arcsinh}\, \eta_{ik}$.
Consequently, the admissible space of $u$ is defined as
\begin{equation*}
\begin{aligned}
\mathcal{U}^{(8)}_{ijk}(\eta)
=\{&(u_i,u_j,u_k)\in \mathbb{R}_{>0}\times\mathbb{R}\times\mathbb{R}_{<0}
\mid u_j+u_k>\log \frac{1}{\eta_{jk}},\ u_i+u_k<\mathrm{arcsinh}\, \eta_{ik} \},
\end{aligned}
\end{equation*}
where $\eta_{jk}>0$, $\eta_{ij}\geq0$, and $\eta_{ik}\in \mathbb{R}$.

\item[(IX)]
Assume $(\alpha_i,\alpha_j,\alpha_k)=(-1,0,-1)$.
It is known that $u_j+u_k>\log \frac{1}{\eta_{jk}}$ and $\eta_{jk}>0$.
Moreover, $\eta_{ij}>0$ and $u_i+u_j>\log\frac{1}{\eta_{ij}}$. 
Similarly, $\eta_{ik}>1$ and $-\mathrm{arccosh}\, \eta_{ik}<u_i+u_k<\mathrm{arccosh}\, \eta_{ik}$.
Consequently, the admissible space of $u$ is defined as
\begin{equation*}
\begin{aligned}
\mathcal{U}^{(9)}_{ijk}(\eta)
=\{(u_i,u_j,u_k)\in \mathbb{R}_{>0}\times\mathbb{R}\times\mathbb{R}_{<0}
\mid &\ u_j+u_k>\log \frac{1}{\eta_{jk}},\ u_i+u_j>\log\frac{1}{\eta_{ij}},\ \\
&\ -\mathrm{arccosh}\, \eta_{ik}<u_i+u_k<\mathrm{arccosh}\, \eta_{ik} \},
\end{aligned}
\end{equation*}
where $\eta_{jk}>0$, $\eta_{ij}>0$, and $\eta_{ik}>1$.

\item[(X)]
Assume $(\alpha_i,\alpha_j,\alpha_k)=(0,1,1)$.
In Subsection \ref{section a4}, it follows that $u_j+u_k>-\mathrm{arccosh}\, \eta_{jk}$ and $\eta_{jk}>1$.
Moreover, $(u_i+u_j)\in \mathbb{R}$ for $\eta_{ij}\geq0$, and $u_i+u_j>\log(-\eta_{ij})$ for $\eta_{ij}<0$.
Similarly, $(u_i+u_k)\in \mathbb{R}$ for $\eta_{ik}\geq0$, and $u_i+u_k>\log(-\eta_{ik})$ for $\eta_{ik}<0$.
Consequently, the admissible space of $u$ is defined as
\begin{equation}\label{AS10}
\begin{aligned}
\mathcal{U}^{(10)}_{ijk}(\eta)
=\{(u_i,u_j,u_k)\in \mathbb{R}\times\mathbb{R}^2_{<0}
\mid &\ u_j+u_k>-\mathrm{arccosh}\, \eta_{jk},\ e^{u_i+u_j}>-\eta_{ij},\\
&\ e^{u_i+u_k}>-\eta_{ik} \},
\end{aligned}
\end{equation}
where $\eta_{jk}>1$, $\eta_{ij}\in \mathbb{R}$, and $\eta_{ik}\in \mathbb{R}$.

The admissible space (\ref{AS10}) can be explicitly characterized as follows:
\begin{description}
\item[(i)]For $\eta_{jk}>1$, $\eta_{ij}\geq0$, and $\eta_{ik}\geq0$, the admissible space is given by
\begin{equation*}
\mathcal{U}^{(10)}_{ijk}(\eta)
=\{(u_i,u_j,u_k)\in \mathbb{R}\times\mathbb{R}^2_{<0}
\mid u_j+u_k>-\mathrm{arccosh}\, \eta_{jk} \}.
\end{equation*}
\item[(ii)]For $\eta_{jk}>1$, $\eta_{ij}\geq0$, and $\eta_{ik}<0$, the admissible space is given by
\begin{equation*}
\mathcal{U}^{(10)}_{ijk}(\eta)
=\{(u_i,u_j,u_k)\in \mathbb{R}\times\mathbb{R}^2_{<0}
\mid u_j+u_k>-\mathrm{arccosh}\,  \eta_{jk},\ u_i+u_k>\log(-\eta_{ik}) \}.
\end{equation*}
\item[(iii)]For $\eta_{jk}>1$, $\eta_{ij}<0$, and $\eta_{ik}\geq0$, the admissible space is given by
\begin{equation*}
\mathcal{U}^{(10)}_{ijk}(\eta)
=\{(u_i,u_j,u_k)\in \mathbb{R}\times\mathbb{R}^2_{<0}
\mid u_j+u_k>-\mathrm{arccosh}\,  \eta_{jk},\ u_i+u_j>\log(-\eta_{ij}) \},
\end{equation*}
\item[(iv)]For $\eta_{jk}>1$, $\eta_{ij}<0$, and $\eta_{ik}<0$, the admissible space is given by
\begin{equation*}
\begin{aligned}
\mathcal{U}^{(10)}_{ijk}(\eta)
=\{&(u_i,u_j,u_k)\in \mathbb{R}\times\mathbb{R}^2_{<0}
\mid u_j+u_k>-\mathrm{arccosh}\,\eta_{jk},\ \\ &u_i+u_j>\log(-\eta_{ij}),\ u_i+u_k>\log(-\eta_{ik}) \}.
\end{aligned}
\end{equation*}
\end{description}

\item[(XI)]
Assume $(\alpha_i,\alpha_j,\alpha_k)=(1,1,1)$.
It is known that $u_j+u_k>-\mathrm{arccosh}\, \eta_{jk}$ and $\eta_{jk}>1$.
Moreover, $\eta_{ij}>-1$ with $(u_i+u_j)\in \mathbb{R}$, and $\eta_{ik}>-1$ with $(u_i+u_k)\in \mathbb{R}$.
Consequently, the admissible space of $u$ is defined as
\begin{equation}\label{AS11}
\mathcal{U}^{(11)}_{ijk}(\eta)
=\{(u_i,u_j,u_k)\in \mathbb{R}_{>0}\times\mathbb{R}^2_{<0}
\mid u_j+u_k>-\mathrm{arccosh}\, \eta_{jk}\},
\end{equation}
where $\eta_{jk}>1$, $\eta_{ij}>-1$, and $\eta_{ik}>-1$.

\item[(XII)]
Assume $(\alpha_i,\alpha_j,\alpha_k)=(-1,1,1)$.
It is known that $u_j+u_k>-\mathrm{arccosh}\, \eta_{jk}$ and $\eta_{jk}>1$.
Moreover, $\eta_{ij}\in \mathbb{R}$ with $u_i+u_j>\mathrm{arcsinh}\, (-\eta_{ij})$, and $\eta_{ik}\in \mathbb{R}$ with $u_i+u_k>\mathrm{arcsinh}\, (-\eta_{ik})$.
Consequently, the admissible space of $u$ is defined as
\begin{equation*}
\begin{aligned}
\mathcal{U}^{(12)}_{ijk}(\eta)
=\{(u_i,u_j,u_k)\in \mathbb{R}_{>0}\times\mathbb{R}^2_{<0}
\mid &\ u_j+u_k>-\mathrm{arccosh}\, \eta_{jk},\
u_i+u_j>\mathrm{arcsinh}\, (-\eta_{ij}),\ \\
&\ u_i+u_k>\mathrm{arcsinh}\, (-\eta_{ik})\},
\end{aligned}
\end{equation*}
where $\eta_{jk}>1$, $\eta_{ij}\in \mathbb{R}$, and $\eta_{ik}\in \mathbb{R}$.

\item[(XIII)]
Assume $(\alpha_i,\alpha_j,\alpha_k)=(0,1,-1)$.
In Subsection \ref{section a4}, it follows that $u_j+u_k>\mathrm{arcsinh}\ (-\eta_{jk})$ and $\eta_{jk}>0$.
Moreover, $(u_i+u_j)\in \mathbb{R}$ for $\eta_{ij}\geq0$, and $u_i+u_j>\log(-\eta_{ij})$ for $\eta_{ij}<0$.
Similarly, $\eta_{ik}>0$ and $u_i+u_k<\log \eta_{ik}$.
When $\eta_{ij}<0$, the inequalities
$\log(-\eta_{ij})+\mathrm{arcsinh}\, (-\eta_{jk})<u_i+2u_j+u_k<2u_j+\log \eta_{ik}<\log \eta_{ik}$ implies $(\eta_{ik}-\eta_{ij})(\eta_{ik}+\eta_{ij})>2\eta_{jk}\eta_{ij}\eta_{ik}$.
For simplicity, we assume $\eta_{ik}+\eta_{ij}>0$.
Consequently, the admissible space of $u$ is defined as
\begin{equation*}
\begin{aligned}
\mathcal{U}^{(13)}_{ijk}(\eta)
=\{(u_i,u_j,u_k)\in \mathbb{R}\times\mathbb{R}^2_{<0}
\mid &\ u_j+u_k>\mathrm{arcsinh}\, (-\eta_{jk}),\
e^{u_i+u_j}>-\eta_{ij},\ \\
&\ u_i+u_k<\log \eta_{ik} \},
\end{aligned}
\end{equation*}
where $\eta_{jk}>0$, $\eta_{ik}>0$, and $\eta_{ij}$ satisfies either $\eta_{ij}\geq0$ or $\eta_{ij}<0$ with $\eta_{ij}+\eta_{ik}>0$.

The admissible space $\mathcal{U}^{(13)}_{ijk}(\eta)$ can be explicitly characterized as follows:
\begin{description}
\item[(i)]For $\eta_{jk}>0$, $\eta_{ik}>0$, and $\eta_{ij}\geq0$, the admissible space is given by
\begin{equation*}
\mathcal{U}_{ijk}^{(13)}(\eta)
=\{(u_i,u_j,u_k)\in \mathbb{R}\times\mathbb{R}^2_{<0}
\mid u_j+u_k>\mathrm{arcsinh}\, (-\eta_{jk}),
u_i+u_k<\log \eta_{ik} \}.
\end{equation*}
\item[(ii)] For $\eta_{jk}>0$, $\eta_{ik}>0$, $\eta_{ij}<0$, and $\eta_{ij}+\eta_{ik}>0$, the admissible space is given by
\begin{align*}
\mathcal{U}_{ijk}^{(13)}(\eta)
=\{(u_i,u_j,u_k)\in \mathbb{R}\times\mathbb{R}^2_{<0}
\mid &\ u_j+u_k>\mathrm{arcsinh}\, (-\eta_{jk}),\
u_i+u_k<\log \eta_{ik},\ \\
&\ u_i+u_j>\log(-\eta_{ij}) \}.
\end{align*}
\end{description}

\item[(XIV)]
Assume $(\alpha_i,\alpha_j,\alpha_k)=(1,1,-1)$.
It is known that $u_j+u_k>\mathrm{arcsinh}\, (-\eta_{jk})$ and $\eta_{jk}>0$.
Moreover, $\eta_{ij}> -1$ and $(u_i+u_j)\in \mathbb{R}$.
Similarly, $\eta_{ik}\in \mathbb{R}$ and $u_i+u_k<\mathrm{arcsinh}\, \eta_{ik}$.
Since $u_i>0$ and $u_j<0$, it follows that $\mathrm{arcsinh}\, (-\eta_{jk})<u_j+u_k<u_i+u_k<\mathrm{arcsinh}\, \eta_{ik}$, and hence $\eta_{ik}+\eta_{jk}>0$.
Consequently, the admissible space of $u$ is defined as
\begin{equation*}
\mathcal{U}^{(14)}_{ijk}(\eta)
=\{(u_i,u_j,u_k)\in \mathbb{R}_{>0}\times\mathbb{R}^2_{<0}
\mid u_j+u_k>\mathrm{arcsinh}\, (-\eta_{jk}),\
u_i+u_k<\mathrm{arcsinh}\, \eta_{ik} \},
\end{equation*}
where $\eta_{jk}>0$, $\eta_{ij}>-1$, $\eta_{ik}\in \mathbb{R}$, and $\eta_{ik}+\eta_{jk}>0$.

\item[(XV)]
Assume $(\alpha_i,\alpha_j,\alpha_k)=(-1,1,-1)$.
It is known that $u_j+u_k>\mathrm{arcsinh}\, (-\eta_{jk})$ and $\eta_{jk}>0$.
Moreover, $\eta_{ij}\in \mathbb{R}$ and $u_i+u_j>\mathrm{arcsinh}\, (-\eta_{ij})$.
Similarly, $\eta_{ik}>1$ and $-\mathrm{arccosh}\, \eta_{ik}<u_i+u_k<\mathrm{arccosh}\, \eta_{ik}$.
Combining these results, we derive the inequality chain: $\mathrm{arcsinh}\, (-\eta_{jk})+\mathrm{arcsinh}\, (-\eta_{ij})<u_i+2u_j+u_k<2u_j+\mathrm{arccosh}\, \eta_{ik}<\mathrm{arccosh}\, \eta_{ik}$.
For simplicity, we assume $\eta_{ij}>0$.
Consequently, the admissible space of $u$ is defined as
\begin{equation*}
\begin{aligned}
\mathcal{U}^{(15)}_{ijk}(\eta)
=\{(u_i,u_j,u_k)\in \mathbb{R}_{>0}\times\mathbb{R}^2_{<0}
\mid &\ u_j+u_k>\mathrm{arcsinh}\, (-\eta_{jk}),\
u_i+u_j>\mathrm{arcsinh}\, (-\eta_{ij}),\ \\
&\ -\mathrm{arccosh}\, \eta_{ik}<u_i+u_k<\mathrm{arccosh}\, \eta_{ik}\},
\end{aligned}
\end{equation*}
where $\eta_{jk}>0$, $\eta_{ij}>0$, and $\eta_{ik}>1$.

\item[(XVI)]
Assume $(\alpha_i,\alpha_j,\alpha_k)=(0,-1,-1)$.
In Subsection \ref{section a4}, it follows that $u_j+u_k>-\mathrm{arccosh}\, \eta_{jk}$ and $\eta_{jk}>1$.
Moreover, $\eta_{ij}>0$ with $u_i+u_j<\log \eta_{ij}$, and $\eta_{ik}>0$ with $u_i+u_k<\log \eta_{ik}$.
Consequently, the admissible space of $u$ is defined as
\begin{equation*}
\begin{aligned}
\mathcal{U}^{(16)}_{ijk}(\eta)
=\{(u_i,u_j,u_k)\in \mathbb{R}\times\mathbb{R}^2_{<0}
\mid &\ u_j+u_k>-\mathrm{arccosh}\, \eta_{jk},\
u_i+u_j<\log \eta_{ij},\ \\
&\ u_i+u_k<\log \eta_{ik}\},
\end{aligned}
\end{equation*}
where $\eta_{jk}>1$, $\eta_{ij}>0$, and $\eta_{ik}>0$.

\item[(XVII)]
Assume $(\alpha_i,\alpha_j,\alpha_k)=(1,-1,-1)$.
It is known that $u_j+u_k>-\mathrm{arccosh}\, \eta_{jk}$ and $\eta_{jk}>1$.
Moreover, $\eta_{ij}\in \mathbb{R}$ and $u_i+u_j<\mathrm{arcsinh}\, \eta_{ij}$.
Since $u_i>0$ and $u_k<0$, it follows that  $-\mathrm{arccosh}\, \eta_{jk}<u_j+u_k<u_i+u_j<\mathrm{arcsinh}\, \eta_{ij}$.
For simplicity, we assume $\eta_{ij}>0$.
Similarly, $\eta_{ik}\in \mathbb{R}$ and $u_i+u_k<\mathrm{arcsinh}\, \eta_{ik}$, thus we assume $\eta_{ik}>0$.
Consequently, the admissible space of $u$ is defined as
\begin{equation*}
\begin{aligned}
\mathcal{U}^{(17)}_{ijk}(\eta)
=\{&(u_i,u_j,u_k)\in \mathbb{R}_{>0}\times\mathbb{R}^2_{<0}
\mid u_j+u_k>-\mathrm{arccosh}\, \eta_{jk}, \\
&\ u_i+u_j<\mathrm{arcsinh}\, \eta_{ij},\
u_i+u_k<\mathrm{arcsinh}\, \eta_{ik}\},
\end{aligned}
\end{equation*}
where $\eta_{jk}>1$, $\eta_{ij}>0$, and $\eta_{ik}>0$.

\item[(XVIII)]
Assume $(\alpha_i,\alpha_j,\alpha_k)=(-1,-1,-1)$.
It is known that $u_j+u_k>-\mathrm{arccosh}\, \eta_{jk}$ and $\eta_{jk}>1$.
Moreover, $\eta_{ij}>1$ and $-\mathrm{arccosh}\, \eta_{ij}<u_i+u_j<\mathrm{arccosh}\, \eta_{ij}$.
Similarly, $\eta_{ik}>1$ and $-\mathrm{arccosh}\, \eta_{ik}<u_i+u_k<\mathrm{arccosh}\, \eta_{ik}$.
Consequently, the admissible space of $u$ is defined as
\begin{equation*}
\begin{aligned}
\mathcal{U}^{(18)}_{ijk}(\eta)
=\{&(u_i,u_j,u_k)\in \mathbb{R}_{>0}\times\mathbb{R}^2_{<0}
\mid u_j+u_k>-\mathrm{arccosh}\, \eta_{jk}, \\
&-\mathrm{arccosh}\, \eta_{ij}<u_i+u_j<\mathrm{arccosh}\, \eta_{ij},\ \\
&-\mathrm{arccosh}\, \eta_{ik}<u_i+u_k<\mathrm{arccosh}\, \eta_{ik}\},
\end{aligned}
\end{equation*}
where $\eta_{jk}>1$, $\eta_{ij}>1$, and $\eta_{ik}>1$.
\end{description}

Across the various admissible spaces $\mathcal{U}^{(n)}_{ijk}(\eta)$ discussed above, 
the weight $\eta$ takes on different value ranges. 
These ranges cannot be expressed through a single consistent formula. 
For simplicity, we introduce the notation $\Xi$ to denote the set of all such weight ranges associated with the entire family of admissible spaces $\mathcal{U}^{(n)}_{ijk}(\eta)$, where $n = 1, 2, \dots, 18$.

\begin{theorem}\label{Thm: ASC 3}
Let $\{ijk\}\in F$ be a right-angled hyperbolic  hexagon with edge lengths given by (\ref{Eq: F38}).
Then the admissible space $\mathcal{U}^{(n)}_{ijk}(\eta)$ of $u$ admits 18 distinct types, each of which is a convex polytope.
\end{theorem}

We rewrite the admissible space in (\ref{Eq: admissible space a2}) as $\mathcal{U}^{(0)}_{ijk}(\eta)$.
Combining Theorem \ref{Thm: ASC 1} and Theorem \ref{Thm: ASC 3} gives the following result.

\begin{corollary}\label{Cor: ASC 3}
Suppose $(\Sigma,\mathcal{T},\alpha,\eta)$ is a weighted triangulated surface with boundary, where $\alpha: B\rightarrow \{-1,0,1\}$.
If one of the following conditions holds: 
\begin{description}
\item[(1)]
All right-angled hyperbolic hexagons are Type-I (i.e., the edge lengths of any right-angled hyperbolic hexagons are given by (\ref{Eq: F38})),
and the weight $\eta$ belongs to $\Xi$;
\item[(2)]
Not all right-angled hyperbolic hexagons are Type-I (i.e., there exist some right-angled hyperbolic hexagons with edge lengths given by (\ref{Eq: DCS3})).
In this case, for those right-angled hyperbolic hexagons with edge lengths given by (\ref{Eq: F38}),
the weight $\eta$ belongs to $\Xi$,
while for those with edge lengths given by (\ref{Eq: DCS3}), $\eta$ satisfies the condition in Theorem \ref{Thm: ASC 1}.
\end{description}
Then the admissible space
$\mathcal{U}(\eta)=\bigcap_{\{ijk\}\in F}\mathcal{U}_{ijk}(\eta)$ is a convex polytope on $(\Sigma,\mathcal{T},\eta)$.
Here $\mathcal{U}_{ijk}(\eta)\in \{\mathcal{U}^{(n)}_{ijk}(\eta),\ n=0,1,2,...,18\}$.
\end{corollary}

%\begin{remark}
%The admissible space $\mathcal{U}(\eta)$ is a finite intersection of 18 types of admissible spaces $\mathcal{U}^{(n)}_{ijk}(\eta)$ for $n=1,2,...,18$.
%But this finite intersection is not arbitrary because of the connectedness of the admissible space $\mathcal{U}(\eta)$.
%For example, if there exist the admissible space $\mathcal{U}^{(1)}_{ijk}(\eta)$ and $\mathcal{U}^{(18)}_{ijk}(\eta)$ for some right-angled hyperbolic  hexagons $\{ijk\}\in F$ with edge lengths given by (\ref{Eq: F38}), then there must be the third type of the admissible space $\mathcal{U}^{(n)}_{ijk}(\eta)$ for $n\neq1,18$.
%\end{remark}

\subsection{Negative definiteness of the Jacobian}

\begin{theorem}\label{Thm: matrix negative 4}
Under the same assumptions as those in Theorem \ref{Thm: ASC 3},
the Jacobian $\Lambda_{ijk}=\frac{\partial(\theta^{jk}_i,\theta^{ik}_j,\theta^{ij}_k)}
{\partial(u_i,u_j,u_k)}$ is symmetric and negative definite on $\mathcal{U}_{ijk}(\eta)$.
\end{theorem}
\proof
For simplicity, we adopt the notations introduction in Subsection \ref{subsec: matrix}.
Formula (\ref{Eq: F27}) implies
\begin{equation*}
\frac{\partial(f_i,f_j,f_k)}
{\partial(u_i,u_j,u_k)}
=\left(
   \begin{array}{ccc}
     -C_i & 0 & 0 \\
     0 & C_j & 0 \\
     0 & 0 & C_k \\
   \end{array}
 \right)
 :=Q_8,
\end{equation*}
where $C_r=\sqrt{1+\alpha_re^{2f_r}}$.
Then $\Lambda_{ijk}$ can be expressed as 
\begin{equation*}
\Lambda_{ijk}=
\frac{\partial(\theta_i,\theta_j,\theta_k)}
{\partial(u_i,u_j,u_k)}
=-\frac{1}{A}
 \left(
   \begin{array}{ccc}
     \sinh l_i & 0 & 0 \\
     0 & \sinh l_j & 0 \\
     0 & 0 & \sinh l_k \\
   \end{array}
 \right)Q_1Q_2Q_8,
\end{equation*}
where $Q_1$ and $Q_2$ are defined by (\ref{Eq: Q_1})  and (\ref{Eq: Q_2}), respectively.
Furthermore, $\det Q_1>0$.
For any right-angled hyperbolic hexagon $\{ijk\}\in F$ with edge lengths given by (\ref{Eq: F38}),
since $d_{ij}<0$, $d_{ik}<0$, $d_{ji}>0$, $d_{ki}>0$, $d_{jk}>0$, and $d_{kj}>0$,
it follows from (\ref{Eq: F23}) that $\det Q_2<0$.
The remaining part of the proof is analogously to that of Theorem \ref{Thm: matrix negative 3}, and is thus omitted here.

\qed

\begin{corollary}
Under the same assumptions as those in Corollary \ref{Cor: ASC 3},
the Jacobian $\Lambda=\frac{\partial (K_i,..., K_N)}{\partial(u_i,...,u_N)}$ is symmetric and negative definite on $\mathcal{U}(\eta)$.
\end{corollary}

\subsection{Rigidity of the mixed discrete conformal structure I}
The following theorem establishes the rigidity of the mixed discrete conformal structure I,
which generalizes the rigidity part of Theorem \ref{Thm: rigidity and image} (iv).
As its proof is similar to that of Theorem \ref{Thm: rigidity 1}, it is omitted here.

\begin{theorem}\label{Thm: rigidity 4}
Under the same assumptions as those in Corollary \ref{Cor: ASC 3},
the discrete conformal factor $f$ is uniquely determined by its generalized combinatorial curvature $K\in\mathbb{R}^N_{>0}$.
%In particular, the map from $f$ to $K$ is a smooth embedding.
\end{theorem}

\subsection{Existence of the mixed discrete conformal structure I}\label{subsection: image 6}

For Type-I right-angled hyperbolic hexagons, we restrict our consideration to four types of admissible spaces, specifically corresponding to $(0,0,0)$, $(1,0,0)$, $(0,0,1)$, and $(1,0,1)$. 
These restrictions are motivated by two considerations: first, all cases where $(\alpha_i, \alpha_j, \alpha_k)$ includes $-1$ are excluded, as the inclusion of $-1$ introduces significant technical complexities in the proof. 
Second, the case $\alpha_j = \alpha_k = 1$ is also excluded, with a detailed explanation provided in Remark \ref{Rmk: 4}.
For simplicity, we introduce a new notation $\widetilde{\Xi}$ to denote the set consisting of the ranges of weights $\eta$ across all admissible spaces $\mathcal{U}^{(n)}_{ijk}$ with $n = 1, 2, 4, 5$.

The following theorem characterizes the image of $K$ for the mixed discrete conformal structure I,
which generalizes the existence part of Theorem \ref{Thm: rigidity and image} (iv).

\begin{theorem}\label{Thm: image 4}
Suppose $(\Sigma,\mathcal{T},\alpha,\eta)$ is a weighted triangulated surface with boundary, where $\alpha: B\rightarrow \{0,1\}$.
If one of the following conditions holds: 
\begin{description}
\item[(1)]
All right-angled hyperbolic hexagons are Type-I (i.e., the edge lengths of any right-angled hyperbolic hexagon are given by (\ref{Eq: F38})), $\alpha_j$ and $\alpha_k$ cannot both be 1 simultaneously, and $\eta\in \widetilde{\Xi}$;
\item[(2)]
Not all right-angled hyperbolic hexagons are Type-I (i.e., there exist some right-angled hyperbolic hexagons with edge lengths given by (\ref{Eq: DCS3})).
In this case, for those right-angled hyperbolic hexagons with edge lengths given by (\ref{Eq: F38}), $\alpha_j$ and $\alpha_k$ cannot both be 1 simultaneously, and $\eta\in \widetilde{\Xi}$,
while for those with edge lengths given by (\ref{Eq: DCS3}), $\eta$ satisfies the condition in Theorem \ref{Thm: ASC 1}.
\end{description}
Then the image of $K$ is $\mathbb{R}^N_{>0}$.
\end{theorem}

To prove Theorem \ref{Thm: image 4}, we need the following four lemmas.

\begin{lemma}\label{Lem: limit d1}
Suppose $(\Sigma,\mathcal{T},\alpha,\eta)$ is a weighted triangulated surface with boundary, where the weights are given by $\alpha: B\rightarrow \{0,1\}$ and $\eta\in \widetilde{\Xi}$.
Let $\{ijk\}\in F$ be a right-angled hyperbolic  hexagon with edge lengths given by (\ref{Eq: F38}), and assume $\alpha_j$ and $\alpha_k$ cannot both be 1 simultaneously.
If one of the following conditions is satisfied
\begin{description}
\item[(i)] $\lim f_i=+\infty,\ \lim f_j=+\infty,\ \lim f_k=+\infty$;
\item[(ii)] $\lim f_i=+\infty,\ \lim f_j=+\infty,\ \lim f_k=c_1$;
\item[(iii)] $\lim f_i=+\infty,\ \lim f_j=c_2,\ \lim f_k=c_3$,
\end{description}
where$c_1,c_2,c_3\in \mathbb{R}$ are constants.
Then $\lim\theta^{jk}_i(f_i,f_j,f_k)=0$.
\end{lemma}
\proof
Formula (\ref{Eq: F38}) implies
\begin{equation}\label{Eq: F3}
\begin{aligned}
&\cosh l_{ij}=\bigg[\eta_{ij}+\sqrt{(e^{-2f_i}+\alpha_i)
(e^{-2f_j}+\alpha_j)}\bigg]e^{f_i+f_j},\\
&\cosh l_{ik}=\bigg[\eta_{ik}+\sqrt{(e^{-2f_i}+\alpha_i)
(e^{-2f_k}+\alpha_k)}\bigg]e^{f_i+f_k},\\
&\cosh l_{jk}=\bigg[\eta_{jk}-\sqrt{(e^{-2f_j}+\alpha_j)
(e^{-2f_k}+\alpha_k)}\bigg]e^{f_j+f_k}.
\end{aligned}
\end{equation}
For the case (i), it follows that $\lim\cosh l_{ij}:=\lim (c_ke^{f_i+f_j})$, $\lim\cosh l_{ik}:=\lim (c_je^{f_i+f_k})$, and $\lim\cosh l_{jk}:=\lim (c_ie^{f_j+f_k})$, 
where $c_i,c_j,c_k$ are positive constants.
By the cosine law, we have
\begin{equation*}
\cosh \theta^{jk}_i=\frac{\cosh l_{jk}+\cosh l_{ij}\cosh l_{ik}}{\sinh l_{ij}\sinh l_{ik}}\rightarrow
\frac{c_ie^{f_j+f_k}}{c_ke^{f_i+f_j}\cdot c_je^{f_i+f_k}}+1
\rightarrow 1,
\end{equation*}
which implies $\theta^{jk}_i\rightarrow 0$.
The proofs for the cases (ii) and (iii) follow analogously and are omitted here.
\qed

\begin{lemma}\label{Lem: limit d2}
Under the same assumptions as those in Lemma \ref{Lem: limit d1},
if $\lim f_k=-\infty$, then $\lim f_j=+\infty$.
Moreover, if one of the following conditions is satisfied
\begin{description}
\item[(i)] $\lim f_i=+\infty,\ \lim f_j=+\infty,\ \lim f_k=-\infty$,
\item[(ii)] $\lim f_i=c,\ \lim f_j=+\infty,\ \lim f_k=-\infty$,
\end{description}
then $\lim\theta^{ij}_k(f_i,f_j,f_k)=+\infty$,
where $c\in \mathbb{R}$ is a constant.
\end{lemma}
\proof
Combining (\ref{Eq: F27}) and Remark \ref{Rmk: 2},
if $\lim f_k=-\infty$, then $\lim u_k=-\infty$.
As established in Subsection \ref{subsection: AS6}, the condition $l_{jk}>0$ implies $u_j+u_k>C(\eta_{jk})$, 
where $C(\eta_{jk})$ is a constant depending on $\eta_{jk}$.
This, in turn, forces $\lim u_j=+\infty$.
By combining (\ref{Eq: F27}) and Remark \ref{Rmk: 2} again, we have $\alpha_j=0$ and $\lim f_j=+\infty$.

Set $\widetilde{L}_1=\frac{\cosh l_{ij}}{\sinh l_{ik}\sinh l_{jk}}$ and $\widetilde{L}_2=\frac{\cosh l_{ik}\cosh l_{jk}}{\sinh l_{ik}\sinh l_{jk}}$,
then $\cosh \theta_k^{ij}=\widetilde{L}_1+\widetilde{L}_2$.
Moreover, (\ref{Eq: F38}) reduces to
\begin{equation*}
\begin{aligned}
\cosh l_{ij}
&=\sqrt{1+\alpha_ie^{2f_i}}
+\eta_{ij}e^{f_i+f_j},\\
\cosh l_{ik}
&\rightarrow\sqrt{1+\alpha_ie^{2f_i}}
+\eta_{ik}e^{f_i+f_k},\\
\cosh l_{jk}
&\rightarrow -1+\eta_{jk}e^{f_j+f_k}.
\end{aligned}
\end{equation*}

For the case (i), it follows that $\lim\cosh l_{ij}:=\lim c_ke^{f_i+f_j}$, $\lim\cosh l_{jk}=\lim (-1+\eta_{jk}e^{f_j+f_k})$, and
\begin{eqnarray*}
\lim\cosh l_{ik}=
\begin{cases}
\lim (1+\eta_{ik}e^{f_i+f_k}),\ & \text{if}\ \alpha_i=0,  \\
\lim(c_je^{f_i}),\ & \text{if}\ \alpha_i=1,
\end{cases}
\end{eqnarray*}
where $c_j,c_k$ are positive constants.
If $\alpha_i=0$, we have
\begin{equation*}
\begin{aligned}
\widetilde{L}_1\rightarrow&\frac{c_ke^{f_i+f_j}}
{\sqrt{(-1+\eta_{jk}e^{f_j+f_k})^2-1}\cdot
\sqrt{(1+\eta_{ik}e^{f_i+f_k})^2-1}}\\
=&\frac{c_ke^{f_i+f_j}}
{\sqrt{\eta_{jk}\eta_{ik}}e^{f_i+f_j}}\cdot
\frac{1}{\sqrt{\eta_{jk}e^{2f_k}-2e^{f_k-f_j}}}\cdot
\frac{1}{\sqrt{\eta_{ik}e^{2f_k}+2e^{f_k-f_i}}}
\rightarrow+\infty.
\end{aligned}
\end{equation*}
Since $\widetilde{L}_2>0$, it follows that $\cosh \theta^{ij}_k=\widetilde{L}_1+\widetilde{L}_2\rightarrow +\infty$,
which implies $\theta^{ij}_k\rightarrow +\infty$.
The proof of the case that $\alpha_i=1$ is similar, so we omit it here.

For the case (ii), it follows that $\lim\cosh l_{ij}:=\lim c^\prime_ke^{f_j}$, $\lim\cosh l_{jk}=\lim (-1+\eta_{jk}e^{f_j+f_k})$, and
\begin{eqnarray*}
\lim\cosh l_{ik}=
\begin{cases}
1,\ & \text{if}\ \alpha_i=0,  \\
c^\prime_j,\ & \text{if}\ \alpha_i=1,
\end{cases}
\end{eqnarray*}
where $c^\prime_j,c^\prime_k$ are positive constants.
If $\alpha_i=0$, then $\cosh l_{ik}\rightarrow 1$, and hence $l_{ik}\rightarrow 0$.
This implies $\cosh \theta^{ij}_k>\widetilde{L}_2>\frac{\cosh l_{ik}}{\sinh l_{ik}}\rightarrow +\infty$, and hence
$\theta^{ij}_k\rightarrow +\infty$.
If $\alpha_i=1$, similar to the proof of the case (i), we derive $\widetilde{L}_1\rightarrow +\infty$.
Consequently, $\theta^{ij}_k\rightarrow +\infty$.
This completes the proof.

\qed

\begin{lemma}\label{Lem: limit d3}
Under the same assumptions as those in Lemma \ref{Lem: limit d1},
if one of the following conditions is satisfied
\begin{description}
\item[(i)] $\lim f_i=-\infty,\ \lim f_j=+\infty,\ \lim f_k=+\infty$,
\item[(ii)] $\lim f_i=-\infty,\ \lim f_j=+\infty,\ \lim f_k=c_1$,
\item[(iii)] $\lim f_i=-\infty,\ \lim f_j=c_2,\ \lim f_k=c_3$,
\item[(iv)] $\lim f_i=-\infty,\ \lim f_j=+\infty,\ \lim f_k=-\infty$,
\end{description}
where $c_1,c_2,c_3\in \mathbb{R}$ are constants.
Then $\lim\theta^{jk}_i(f_i,f_j,f_k)=+\infty$.
\end{lemma}
\proof
Set $L_1=\frac{\cosh l_{jk}}{\sinh l_{ij}\sinh l_{ik}}$ and $L_2=\frac{\cosh l_{ij}\cosh l_{ik}}{\sinh l_{ij}\sinh l_{ik}}$, then $\cosh \theta^{jk}_i=L_1+L_2$.
Given that $\lim f_i=-\infty$, the formula (\ref{Eq: F38}) reduces to
\begin{equation}\label{Eq: F2}
\begin{aligned}
\cosh l_{ij}
&\rightarrow\sqrt{1+\alpha_je^{2f_j}}
+\eta_{ij}e^{f_i+f_j},\\
\cosh l_{ik}
&\rightarrow\sqrt{1+\alpha_ke^{2f_k}}
+\eta_{ik}e^{f_i+f_k},\\
\cosh l_{jk}
&=-\sqrt{(1+\alpha_je^{2f_j})(1+\alpha_ke^{2f_k})}
+\eta_{jk}e^{f_j+f_k}.
\end{aligned}
\end{equation}

For the case (i), it follows from (\ref{Eq: F2}) that $\lim\cosh l_{jk}:=\lim (c_ie^{f_j+f_k})$,
\begin{eqnarray*}
\lim\cosh l_{ij}=
\begin{cases}
\lim(1+\eta_{ij}e^{f_i+f_j}),\ & \text{if}\ \alpha_j=0,  \\
\lim(c_ke^{f_j}),\ & \text{if}\ \alpha_j=1,
\end{cases}
\end{eqnarray*}
and
\begin{eqnarray*}
\lim\cosh l_{ik}=
\begin{cases}
\lim(1+\eta_{ik}e^{f_i+f_k}),\ & \text{if}\ \alpha_k=0,  \\
\lim(c_je^{f_k}),\ & \text{if}\ \alpha_k=1,
\end{cases}
\end{eqnarray*}
where $c_i,c_j,c_k$ are positive constants.
If $\alpha_j=\alpha_k=0$, then
\begin{equation*}
\begin{aligned}
L_1\rightarrow &\frac{c_ie^{f_j+f_k}}
{\sqrt{(1+\eta_{ij}e^{f_i+f_j})^2-1}\cdot
\sqrt{(1+\eta_{ik}e^{f_i+f_k})^2-1}}\\
=&\frac{c_ie^{f_j+f_k}}
{\sqrt{\eta_{ij}\eta_{ik}}e^{f_j+f_k}}\cdot
\frac{1}{\sqrt{\eta_{ij}e^{2f_i}+2e^{f_i-f_j}}}\cdot
\frac{1}{\sqrt{\eta_{ik}e^{2f_i}+2e^{f_i-f_k}}}
\rightarrow+\infty.
\end{aligned}
\end{equation*}
Since $L_2>0$, it follows that $\cosh \theta^{jk}_i=L_1+L_2\rightarrow +\infty$, 
which implies $\theta^{jk}_i\rightarrow +\infty$.
For the cases that $\alpha_j=0,\ \alpha_k=1$ and $\alpha_j=1,\ \alpha_k=0$, analogous reasoning yields   $\theta^{jk}_i\rightarrow +\infty$.

For the case (ii), it follows from (\ref{Eq: F2}) that  $\lim\cosh l_{jk}:=\lim(c^\prime_ie^{f_j})$,
and
\begin{eqnarray*}
\lim\cosh l_{ik}=
\begin{cases}
1,\ & \text{if}\ \alpha_k=0,  \\
c^\prime_j,\ & \text{if}\ \alpha_k=1,
\end{cases}
\end{eqnarray*}
where $c^\prime_i,c^\prime_j$ are positive constants.
If $\alpha_k=0$, then $l_{ik}\rightarrow 0$, and hence $L_2>\frac{\cosh l_{ik}}{\sinh l_{ik}}\rightarrow +\infty$.
This implies $\cosh\theta^{jk}_i\rightarrow +\infty$.
If $\alpha_k=1$, then by the assumption, we have $\alpha_j=0$.
Hence, $\lim\cosh l_{ij}=\lim(1+\eta_{ij}e^{f_i+f_j})$,
$\lim\cosh l_{jk}=\lim(c^\prime_ie^{f_j})$, and
$\lim\cosh l_{ik}=c^\prime_j$.
Consequently, 
\begin{equation*}
\begin{aligned}
L_1\rightarrow &\frac{c^\prime_ie^{f_j}}
{\sqrt{(1+\eta_{ij}e^{f_i+f_j})^2-1}\cdot
\sqrt{(c^\prime_j)^2-1}}\\
=&\frac{c^\prime_ie^{f_j}}
{\sqrt{\eta_{ij}}\sqrt{(c^\prime_j)^2-1}e^{f_j}}\cdot
\frac{1}{\sqrt{\eta_{ij}e^{2f_i}+2e^{f_i-f_j}}}
\rightarrow+\infty.
\end{aligned}
\end{equation*}
Since $L_2>0$, it follows that $\cosh \theta^{jk}_i=L_1+L_2\rightarrow +\infty$, 
which implies $\theta^{jk}_i\rightarrow +\infty$.

For the case (iii), it follows from (\ref{Eq: F2}) that $\lim\cosh l_{jk}:=c^{\prime\prime}_i$, and
\begin{eqnarray*}
\lim\cosh l_{ij}=
\begin{cases}
1,\ & \text{if}\ \alpha_j=0,  \\
c^{\prime\prime}_k,\ & \text{if}\ \alpha_j=1,
\end{cases}
\end{eqnarray*}
where $c^{\prime\prime}_i,c^{\prime\prime}_k$ are positive constants.
If $\alpha_j=0$, then $l_{ij}\rightarrow 0$, and hence $L_2\rightarrow +\infty$.
This implies $\theta^{jk}_i\rightarrow +\infty$.
If $\alpha_j=1$, then by the assumption, $\alpha_k=0$.
According to the arguments in the case (ii),
we still have $\theta^{jk}_i\rightarrow +\infty$.

For the case (iv), it follows from (\ref{Eq: F2}) that $\lim\cosh l_{ik}=1$, and hence $l_{ik}\rightarrow 0$.
Consequently, $\cosh\theta^{jk}_i>L_2\geq\frac{\cosh l_{ik}}{\sinh l_{ik}}\rightarrow +\infty$, which implies $\theta^{jk}_i\rightarrow +\infty$.
This completes the proof.
\qed

\begin{remark}\label{Rmk: 4}
As shown in the proof of case (iii) in Lemma \ref{Lem: limit d3}, if $\alpha_j = \alpha_k = 1$, then $\lim\cosh l_{ij} = c''_k$, $\lim\cosh l_{jk} = c''_i$, and $\lim\cosh l_{ik} = c'_j$. 
In this case, $\cosh \theta^{jk}_i$, $\cosh \theta^{ik}_j$ and $\cosh \theta^{ij}_k$ all tend to finite constants, 
which neither approach 1 nor $+\infty$. 
Consequently, no contradiction can be derived from this case in the subsequent analysis, 
thereby justifying the exclusion of the case that $\alpha_j = \alpha_k = 1$.
\end{remark}

\begin{lemma}\label{Lem: limit d4}
Under the same assumptions as those in Lemma \ref{Lem: limit d1},
if one of the following conditions is satisfied
\begin{description}
\item[(i)] $\lim f_j=+\infty,\ \lim f_i=c_1,\ \lim f_k=+\infty$,
\item[(ii)] $\lim f_j=+\infty,\ \lim f_i=c_2,\ \lim f_k=c_3$,
\end{description}
where $c_1,c_2,c_3\in \mathbb{R}$ are constants.
Then $\lim\theta^{ik}_j(f_i,f_j,f_k)=0$.
\end{lemma}
\proof
For the case (i), by (\ref{Eq: F3}), we have
$\lim\cosh l_{ij}:=\lim(c_ke^{f_j})$,
$\lim\cosh l_{ik}:=c_je^{f_k}$, and $\lim\cosh l_{jk}:=\lim(c_ie^{f_j+f_k})$,
where $c_i,c_j,c_k$ are positive constants.
Then
\begin{equation*}
\cosh \theta^{ik}_j=\frac{\cosh l_{ik}+\cosh l_{ij}\cosh l_{jk}}{\sinh l_{ij}\sinh l_{jk}}\rightarrow
\frac{c_je^{f_k}}{c_ke^{f_j}\cdot c_ie^{f_j+f_k}}+1
\rightarrow 1,
\end{equation*}
which implies $\theta^{ik}_j\rightarrow 0$.
The proof of the case (ii) is similar, so we omit it here.
\qed

\noindent\textbf{Proof of Theorem \ref{Thm: image 4}:}
Let $X$ be the image of $K$.
The openness of $X$ in $\mathbb{R}^N_{>0}$ follows from the smooth injectivity of the map.
To prove the closedness of $X$ in $\mathbb{R}^N_{>0}$, we prove that if a sequence $\{f^{(m)}\}\in\Omega(\eta)=\bigcap_{\{ijk\}\in F}\Omega_{ijk}(\eta)$ satisfies $\lim_{m\rightarrow +\infty}K^{(m)}\in \mathbb{R}^N_{>0}$, then there exists a subsequence,
still denoted by $\{f^{(m)}\}$, such that $\lim_{m\rightarrow +\infty}f^{(m)}\in\Omega(\eta)$.

Suppose otherwise, there exists a subsequence $\{f^{(m)}\}$, such that $\lim_{m\rightarrow +\infty}f^{(m)}$ lies on the boundary of $\Omega(\eta)$.
The boundary of the admissible space $\Omega(\eta)$ in $[-\infty,+\infty]^N$ consists of the following five parts:
\begin{description}
\item[(i)]
There exists an ideal edge $\{rs\}\in E$ such that $l_{rs}=0$.
\item[(ii)]
There exists $i\in B$ in a right-angled hyperbolic  hexagon with edge lengths given by (\ref{Eq: F38}), such that $f_i^{(m)}\rightarrow \pm\infty$.
\item[(iii)]
There exists $j\in B$ in a right-angled hyperbolic  hexagon with edge lengths given by (\ref{Eq: F38}), such that $f_j^{(m)}\rightarrow \pm\infty$.
\item[(iv)]
There exists $k\in B$ in a right-angled hyperbolic  hexagon with edge lengths given by (\ref{Eq: F38}), such that $f_k^{(m)}\rightarrow \pm\infty$.
\item[(v)]
There exists $j\in B$ in a right-angled hyperbolic  hexagon with edge lengths given by (\ref{Eq: DCS3}), such that $f_j^{(m)}\rightarrow \pm\infty$.
\end{description}
The proof for case (i) is analogous to that of Theorem \ref{Thm: image 1}, and is thus omitted here.

We first establish that the aforementioned cases (ii), (iii), and (iv) cannot hold on $\mathcal{U}^{(n)}(\eta)=\bigcap_{\{ijk\}\in F}\mathcal{U}^{(n)}_{ijk}(\eta)$, where $n=1,2,4,5$.
The corresponding results are stated as follows.
\begin{itemize}
\item[(1)]
For $(0,0,0)$-type admissible space $\mathcal{U}^{(1)}_{ijk}(\eta)$ in (\ref{Eq: F4}), we have $u_j+u_k>\log \frac{2}{\eta_{jk}}$.
By Remark \ref{Rmk: 2}, $u_j=f_j$ and $u_k=f_k$, then $f_j+f_k$ has a lower bound.
Hence, both $f_j$ and $f_k$ cannot tend to $-\infty$ simultaneously.

If $f_i^{(m)}\rightarrow +\infty$, $f_j^{(m)}\nrightarrow -\infty$, and $f_k^{(m)}\nrightarrow -\infty$,
then $(\theta^{jk}_i)^{(m)}\rightarrow 0$ by Lemma \ref{Lem: limit d1}.
If $f_i^{(m)}\rightarrow +\infty$ and $f_k^{(m)}\rightarrow -\infty$, then $f_j^{(m)}\rightarrow +\infty$.
It follows from Lemma \ref{Lem: limit d2} (i) that $(\theta^{ij}_k)^{(m)}\rightarrow +\infty$.
By swapping the positions of $j$ and $k$ in Lemma \ref{Lem: limit d2} (i), we conclude that as $f_i^{(m)}\rightarrow +\infty$, $f_j^{(m)}\rightarrow -\infty$, and $f_k^{(m)}\rightarrow +\infty$,
it follows that $(\theta^{ik}_j)^{(m)}\rightarrow +\infty$.
Hence, as $f_i^{(m)}\rightarrow +\infty$, either $K^{(m)}_i\rightarrow 0$ or $K^{(m)}_s\rightarrow +\infty$ for $s\in\{j,k\}$.
Both scenarios contradict the assumption that $\lim_{m\rightarrow +\infty}K^{(m)}\in \mathbb{R}^N_{>0}$.
If $f_i^{(m)}\rightarrow -\infty$, by Lemma \ref{Lem: limit d3}, then $(\theta^{jk}_i)^{(m)}\rightarrow +\infty$ and hence $K^{(m)}_i\rightarrow +\infty$.
Consequently, the subsequence $f^{(m)}_i$ cannot tend to $\pm\infty$ in the admissible space $\mathcal{U}^{(1)}_{ijk}(\eta)$.

If $f_j^{(m)}\rightarrow +\infty$ and  $f_k^{(m)}\rightarrow -\infty$, by Lemma \ref{Lem: limit d2} (ii), then $(\theta^{ij}_k)^{(m)}\rightarrow +\infty$ and hence $K^{(m)}_k\rightarrow +\infty$.
If $f_j^{(m)}\rightarrow +\infty$ and  $f_k^{(m)}\nrightarrow -\infty$, by Lemma \ref{Lem: limit d4},
then $(\theta^{ik}_j)^{(m)}\rightarrow 0$ and hence $K^{(m)}_j\rightarrow 0$.
If $f_j^{(m)}\rightarrow -\infty$, then $f_k^{(m)}\rightarrow +\infty$.
Swapping the positions of $j$ and $k$ in Lemma \ref{Lem: limit d2} (ii) gives $(\theta^{ij}_k)^{(m)}\rightarrow +\infty$, and
hence $K^{(m)}_k\rightarrow +\infty$.
Consequently, the subsequence $f^{(m)}_j$ cannot tend to $\pm\infty$ in the admissible space $\mathcal{U}^{(1)}_{ijk}(\eta)$.

By swapping the positions of $j$ and $k$, the subsequence $f^{(m)}_k$ cannot tend to $\pm\infty$ in the admissible space $\mathcal{U}^{(1)}_{ijk}(\eta)$.

\item[(2)]
For $(1,0,0)$-type admissible space $\mathcal{U}^{(2)}_{ijk}(\eta)$ in (\ref{Eq: F5}), we have $u_j+u_k>\log \frac{2}{\eta_{jk}}$.
The proof follows identically to that of the $(0,0,0)$-type admissible space $\mathcal{U}^{(1)}_{ijk}(\eta)$, and thus is omitted here.

\item[(3)]
For $(0,0,1)$-type admissible space $\mathcal{U}^{(4)}_{ijk}(\eta)$ in (\ref{AS4}), the positions of $j$ and $k$ are not interchangeable.
\begin{itemize}
\item[(a)]
In the admissible space (\ref{Eq: F48}),
we have $u_j+u_k>\log \frac{1}{\eta_{jk}}$.
Combining (\ref{Eq: F27}) and Remark \ref{Rmk: 2} yields $u_j=f_j\in \mathbb{R}$ and $u_k=-\mathrm{arcsinh}\, (e^{-f_k})\in \mathbb{R}_{<0}$.
It follows that $f_j$ cannot tend to $-\infty$.

If $f_i^{(m)}\rightarrow +\infty$, $f_j^{(m)}\nrightarrow -\infty$, and $f_k^{(m)}\nrightarrow -\infty$,
by Lemma \ref{Lem: limit d1},
then $(\theta^{jk}_i)^{(m)}\rightarrow 0$ and hence $K^{(m)}_i\rightarrow 0$.
If $f_i^{(m)}\rightarrow +\infty$ and $f_k^{(m)}\rightarrow -\infty$, then $f_j^{(m)}\rightarrow +\infty$.
By Lemma \ref{Lem: limit d2} (i),  $(\theta^{ij}_k)^{(m)}\rightarrow +\infty$ and hence $K^{(m)}_k\rightarrow +\infty$.
If $f_i^{(m)}\rightarrow -\infty$, by Lemma \ref{Lem: limit d3}, then $(\theta^{jk}_i)^{(m)}\rightarrow +\infty$ and hence $K^{(m)}_i\rightarrow +\infty$.
Consequently, the subsequence $f^{(m)}_i$ cannot tend to $\pm\infty$ in the admissible space (\ref{Eq: F48}).

If $f_j^{(m)}\rightarrow +\infty$ and  $f_k^{(m)}\rightarrow -\infty$, by Lemma \ref{Lem: limit d2} (ii), then $(\theta^{ij}_k)^{(m)}\rightarrow +\infty$ and hence $K^{(m)}_k\rightarrow +\infty$.
If $f_j^{(m)}\rightarrow +\infty$ and  $f_k^{(m)}\nrightarrow -\infty$, by Lemma \ref{Lem: limit d4}, then $(\theta^{ik}_j)^{(m)}\rightarrow 0$ and hence $K^{(m)}_j\rightarrow 0$.
Consequently, the subsequence $f^{(m)}_j$ cannot tend to $+\infty$ in the admissible space (\ref{Eq: F48}).

If $f_k^{(m)}\rightarrow +\infty$, by the above arguments, then $f_i^{(m)}\rightarrow c_1$ and $f_j^{(m)}\rightarrow c_2$, where $c_1,c_2$ are constants.
By swapping the positions of $j$ and $k$ in Lemma \ref{Lem: limit d4} (ii), we have $(\theta^{ij}_k)^{(m)}\rightarrow 0$ and hence $K^{(m)}_k\rightarrow 0$.
If $f_k^{(m)}\rightarrow -\infty$, then $u_k^{(m)}\rightarrow -\infty$ and $u_j^{(m)}\rightarrow +\infty$.
This implies $f_j^{(m)}\rightarrow +\infty$, the case of which has been proved.
Consequently, the subsequence $f^{(m)}_k$ cannot tend to $\pm\infty$ in the admissible space (\ref{Eq: F48}).

\item[(b)]
In the admissible space (\ref{Eq: F42}), we have $u_j+u_k>\log \frac{1}{\eta_{jk}},\  u_i+u_k>\log(-\eta_{ik})$.
The case (a) shows that $f_j$ cannot tend to $-\infty$.
Similarly, combining (\ref{Eq: F27}) and Remark \ref{Rmk: 2} gives $u_i=-f_i$.
Thus $f_i$ cannot tend to $+\infty$.

If $f_i^{(m)}\rightarrow -\infty$, by Lemma \ref{Lem: limit d3}, then $(\theta^{jk}_i)^{(m)}\rightarrow +\infty$ and hence $K^{(m)}_i\rightarrow +\infty$.
Consequently, the subsequence $f^{(m)}_i$ cannot tend to $-\infty$ in the admissible space (\ref{Eq: F42}).

If $f_j^{(m)}\rightarrow +\infty$ and  $f_k^{(m)}\rightarrow -\infty$,
then $f_i^{(m)}\rightarrow -\infty$ must hold to satisfy the inequality $u_i+u_k>\log(-\eta_{ik})$.
This case has been included in the above arguments.
If $f_j^{(m)}\rightarrow +\infty$ and  $f_k^{(m)}\nrightarrow -\infty$,
by Lemma \ref{Lem: limit d4}, then $(\theta^{ik}_j)^{(m)}\rightarrow 0$ and hence $K^{(m)}_j\rightarrow 0$.
Consequently, the subsequence $f^{(m)}_j$ cannot tend to $+\infty$ in the admissible space (\ref{Eq: F42}).

If $f_k^{(m)}\rightarrow +\infty$, by the above arguments, $f_i^{(m)}\rightarrow c_1$ and $f_j^{(m)}\rightarrow c_2$, where $c_1,c_2$ are constants.
By swapping the positions of $j$ and $k$ in Lemma \ref{Lem: limit d4} (ii), we have $(\theta^{ij}_k)^{(m)}\rightarrow 0$ and hence $K^{(m)}_k\rightarrow 0$.
If $f_k^{(m)}\rightarrow -\infty$, then $u_k^{(m)}\rightarrow -\infty$, and $u_j^{(m)}\rightarrow +\infty$.
This implies $f_j^{(m)}\rightarrow +\infty$, the case of which has been proved.
Consequently, the subsequence $f^{(m)}_k$ cannot tend to $\pm\infty$ in the admissible space (\ref{Eq: F42}).
\end{itemize}

\item[(4)]
For $(1,0,1)$-type admissible space $\mathcal{U}^{(5)}_{ijk}(\eta)$ in (\ref{AS5}), we have $u_j+u_k>\log \frac{1}{\eta_{jk}}$.
The proof follows identically to that of the $(0,0,1)$-type admissible space (\ref{Eq: F48}), and thus is omitted here.
\end{itemize}

From the foregoing proofs, it follows that the generalized angle tends to either $0$ or $+\infty$ on the admissible space $\mathcal{U}^{(n)}(\eta)$, 
which is composed of the single admissible subspace $\mathcal{U}^{(n)}_{ijk}(\eta)$ for $n = 1, 2, 4, 5$.
Furthermore, for the case (v), the proof of Theorem \ref{Thm: image 1} demonstrates that the generalized angle $\theta_j^{ik}$ tends to either $0$ or $+\infty$ on $\mathcal{U}^{(0)}_{ijk}(\eta)$.
Through a lengthy yet analogous analysis, one can establish that the generalized angle tends to either $0$ or $+\infty$ on $\widetilde{\mathcal{U}}(\eta)$, which is formed by mixed admissible subspaces $\mathcal{U}^{(n)}_{ijk}(\eta)$ with $n = 0, 1, 2, 4, 5$.
Accordingly, the aforementioned cases (ii), (iii), (iv), and (v) cannot hold on $\widetilde{\mathcal{U}}(\eta)$. 
This completes the proof.

\qed

\appendix
\section{Proofs of Lemmas}\label{appendix a}

For simplicity, we set $h_t=h_{rs}$ and $q_t=q_{rs}$, and assume $\Vert v_i\Vert =\Vert v_j\Vert =\Vert v_k\Vert =1$.
Please refer to Figure \ref{Figure_3}.

\begin{lemma}
If $c_{ijk}$ is time-like, then
\begin{align*}
&\sinh q_r=\cosh h_s\cdot\sinh d_{rt}=\cosh h_t\cdot\sinh d_{rs},\\
&\sinh h_r=\cosh q_s\cdot\sinh \theta_{rt}=\cosh q_t\cdot\sinh \theta_{rs}.
\end{align*}
\end{lemma}
\proof
Without loss of generality, we assume $r=i, s=j, t=k$.

In the generalized right-angled hyperbolic triangle $c_{ijk}c_{ik}v_i$,
by Proposition \ref{Prop: Right angle},
we have $-v_i\ast c_{ijk}=(c_{ik}\ast c_{ijk})(v_i\ast c_{ik})$.
By Proposition \ref{Prop: Ratcliffe 1} (i), $c_{ik}\ast c_{ijk}=-\cosh h_j$.
Moreover, by Proposition \ref{Prop: Ratcliffe 1} (ii),
if $c_{ik}$ and $v_i$ lie on the same side of the hyperplane $v_i^\bot$,
then $$v_i\ast c_{ik}=+\sinh(d_\mathbb{H}(c_{ik},v_i^\bot))=+\sinh (-d_{ik})=-\sinh d_{ik};$$
if $c_{ik}$ and $v_i$ lie on opposite sides of the hyperplane $v_i^\bot$, then
$$v_i\ast c_{ik}=-\sinh(d_\mathbb{H}(c_{ik},v_i^\bot))=-\sinh (+d_{ik})=-\sinh d_{ik}.$$
In either case, $v_i\ast c_{ik}=-\sinh d_{ik}$.
Similarly,
if $c_{ijk}$ and $v_i$ lie on the same side of the hyperplane $v_i^\bot$,
then $$v_i\ast c_{ijk}=+\sinh(d_\mathbb{H}(c_{ijk},v_i^\bot))=+\sinh (-q_i)=-\sinh q_i;$$
if $c_{ijk}$ and $v_i$ lie on opposite sides of the hyperplane $v_i^\bot$, then $$v_i\ast c_{ijk}=-\sinh(d_\mathbb{H}(c_{ijk},v_i^\bot))=-\sinh (+q_i)=-\sinh q_i.$$
In either case, $v_i\ast c_{ijk}=-\sinh q_i$.
Consequently, $\sinh q_i=\cosh h_j\cdot\sinh d_{ik}$.

In the generalized right-angled hyperbolic triangle $c_{ijk}c_{ij}v_i$,
we have $-v_i\ast c_{ijk}=(c_{ij}\ast c_{ijk})(v_i\ast c_{ij})=-\cosh h_k\cdot(v_i\ast c_{ij})$.
Similar arguments can imply $v_i\ast c_{ij}=-\sinh d_{ij}$.
Thus $\sinh q_i=\cosh h_k\cdot\sinh d_{ij}$.

In the generalized right-angled hyperbolic triangle $c_{ijk}c^\prime_{ik}v^\prime_i$,
we have $-v^\prime_i\ast c_{ijk}=(c^\prime_{ik}\ast c_{ijk})(v^\prime_i\ast c^\prime_{ik})=-\cosh q_j\cdot(v^\prime_i\ast c^\prime_{ik})$.
Similarly, $v^\prime_i\ast c_{ijk}=-\sinh h_i$ and
$v^\prime_i\ast c^\prime_{ik}=-\sinh \theta_{ik}$.
Thus $\sinh h_i=\cosh q_j\cdot\sinh \theta_{ik}$.

In the generalized right-angled hyperbolic triangle $c_{ijk}c^\prime_{ij}v^\prime_i$,
we have $-v^\prime_i\ast c_{ijk}=(c^\prime_{ij}\ast c_{ijk})(v^\prime_i\ast c^\prime_{ij})=-\cosh q_k\cdot(v^\prime_i\ast c^\prime_{ij})$.
Similarly,
$v^\prime_i\ast c_{ijk}=-\sinh h_i$ and
$v^\prime_i\ast c^\prime_{ij}=-\sinh \theta_{ij}$.
Thus $\sinh h_i=\cosh q_k\cdot\sinh \theta_{ij}$.
\qed

\begin{proposition}
If $c_{ijk}\in \mathbb{R}^3$ lies in the open domain $D_\mathrm{V}$, then
\begin{equation*}
\begin{aligned}
&\cosh q_i=-\sinh h_j\cdot\sinh d_{ik}=-\sinh h_k\cdot\sinh d_{ij},\\
&\cosh q_j=\sinh h_i\cdot\sinh d_{jk}=\sinh h_k\cdot\sinh d_{ji},\\
&\cosh q_k=\sinh h_i\cdot\sinh d_{kj}=\sinh h_j\cdot\sinh d_{ki},\\
&\cosh h_i=\sinh q_j\cdot\sinh \theta_{ik}=\sinh q_k\cdot\sinh \theta_{ij},\\
&\cosh h_j=-\sinh q_i\cdot\sinh \theta_{jk}=\sinh q_k\cdot\sinh \theta_{ji},\\
&\cosh h_k=-\sinh q_i\cdot\sinh \theta_{kj}=\sinh q_j\cdot\sinh \theta_{ki}.
\end{aligned}
\end{equation*}
If $c_{ijk}\in \mathbb{R}^3$ lies in the open domain $D_\mathrm{II}$, then we swap $h$ with $q$, and $\theta$ with $d$.
\end{proposition}
\proof
\begin{description}
\item[(1)] Assume $c_{ijk}\in D_\mathrm{V}$.
In the generalized right-angled hyperbolic triangle $c_{ijk}c_{ik}v_i$,
by Proposition \ref{Prop: Right angle},
we have $-v_i\ast c_{ijk}=(c_{ik}\ast c_{ijk})(v_i\ast c_{ik})$.
By Proposition \ref{Prop: Ratcliffe 1} (iii),
$v_i\ast c_{ijk}
=\cosh(d_\mathbb{H}(v_i^\bot,c_{ijk}^\bot))
=\cosh(-q_i)=\cosh q_i$.
Moreover, by Proposition \ref{Prop: Ratcliffe 1} (ii),
$v_i\ast c_{ik}=+\sinh d_\mathbb{H}(v_i^\bot,c_{ik})=+\sinh(-d_{ik})=-\sinh d_{ik}$.
If $c_{ik}$ and $c_{ijk}$ lie on the same side of the hyperplane $c_{ijk}^\bot$,
then
\begin{equation*}
c_{ik}\ast c_{ijk}=+\sinh(d_\mathbb{H}(c_{ik},c_{ijk}^\bot))
=+\sinh(-h_j)=-\sinh h_j.
\end{equation*}
If $c_{ik}$ and $c_{ijk}$ lie on opposite sides of the hyperplane $c_{ijk}^\bot$,
then
\begin{equation*}
c_{ik}\ast c_{ijk}=-\sinh(d_\mathbb{H}(c_{ik},c_{ijk}^\bot))
=-\sinh h_j.
\end{equation*}
In either case, $c_{ik}\ast c_{ijk}=-\sinh h_j$.
Consequently, $\cosh q_i=-\sinh h_j\cdot\sinh d_{ik}$.
By Lemma \ref{Lem: coplaner}, it is easy to check that $d_{ik}<0$.
It follows that $h_j>0$ and hence $d_\mathbb{H}(c_{ik},c_{ijk}^\bot)=h_j$.
This implies that $c_{ik}$ and $c_{ijk}$ lie on opposite sides of $c_{ijk}^\bot$.

In the generalized right-angled hyperbolic triangle $c_{ijk}c_{ij}v_i$,
we have $-v_i\ast c_{ijk}=(c_{ij}\ast c_{ijk})(v_i\ast c_{ij})=-\cosh q_i$.
Similarly, $v_i\ast c_{ij}=-\sinh d_{ij}$.
If $c_{ij}$ and $c_{ijk}$ lie on the same side of the hyperplane $c_{ijk}^\bot$,
then $c_{ij}\ast c_{ijk}=+\sinh(d_\mathbb{H}(c_{ij},c_{ijk}^\bot))
=+\sinh(-h_k)=-\sinh h_k$.
If $c_{ij}$ and $c_{ijk}$ lie on opposite sides of the hyperplane $c_{ijk}^\bot$, then $c_{ij}\ast c_{ijk}=-\sinh(d_\mathbb{H}(c_{ij},c_{ijk}^\bot))
=-\sinh h_k$.
In either case, $c_{ij}\ast c_{ijk}=-\sinh h_k$.
Consequently, $\cosh q_i=-\sinh h_k\cdot\sinh d_{ij}$.
Similarly, $d_{ij}<0$ by Lemma \ref{Lem: coplaner} and hence $h_k>0$.
So $d_\mathbb{H}(c_{ij},c_{ijk}^\bot)=h_k$,
which implies that $c_{ij}$ and $c_{ijk}$ lie on opposite sides of $c_{ijk}^\bot$.

In the generalized right-angled hyperbolic triangle $c_{ijk}c_{jk}v_j$,
we have $-v_j\ast c_{ijk}=(c_{jk}\ast c_{ijk})(v_j\ast c_{jk})$.
By Proposition \ref{Prop: Ratcliffe 1} (iii),
$v_j\ast c_{ijk}
=-\cosh(d_\mathbb{H}(v_j^\bot,c_{ijk}^\bot))=-\cosh q_j$.
Moreover, by Proposition \ref{Prop: Ratcliffe 1} (ii),
$c_{jk}\ast c_{ijk}=-\sinh (d_\mathbb{H}(c_{jk},c_{ijk}^\bot))=-\sinh h_i$.
Similarly, $v_j\ast c_{jk}=-\sinh d_{jk}$.
Thus $\cosh q_j=\sinh h_i\cdot \sinh d_{jk}$.

In the generalized right-angled hyperbolic triangle $c_{ijk}c_{ij}v_j$,
we have $-v_j\ast c_{ijk}=(c_{ij}\ast c_{ijk})(v_j\ast c_{ij})$.
Since $v_j\ast c_{ijk}=-\cosh q_j$,
$v_j\ast c_{ij}=-\sinh d_{ji}$,
and $c_{ij}\ast c_{ijk}=-\sinh h_k$,
it follows that $\cosh q_j=\sinh h_k\cdot \sinh d_{ji}$.

In the generalized right-angled hyperbolic triangle $c_{ijk}c_{jk}v_k$,
we have $-v_k\ast c_{ijk}=(c_{jk}\ast c_{ijk})(v_k\ast c_{jk})$.
Since $v_k\ast c_{ijk}=-\cosh q_k$,
$v_k\ast c_{jk}=-\sinh d_{kj}$,
and $c_{jk}\ast c_{ijk}=-\sinh h_i$,
it follows that $\cosh q_k=\sinh h_i\cdot \sinh d_{kj}$.

In the generalized right-angled hyperbolic triangle $c_{ijk}c_{ik}v_k$,
we have $-v_k\ast c_{ijk}=(c_{ik}\ast c_{ijk})(v_k\ast c_{ik})$.
Since $v_k\ast c_{ijk}=-\cosh q_k$,
$v_k\ast c_{ik}=-\sinh d_{ki}$, and
$c_{ik}\ast c_{ijk}=-\sinh h_j$,
it follows that $\cosh q_k=\sinh h_j\cdot\sinh d_{ki}$.

In the generalized right-angled hyperbolic triangle $c_{ijk}c^\prime_{ik}v^\prime_i$,
we have $-v^\prime_i\ast c_{ijk}=(c^\prime_{ik}\ast c_{ijk})(v^\prime_i\ast c^\prime_{ik})$.
Since $v^\prime_i\ast c_{ijk}=-\cosh h_i$,
$v^\prime_i\ast c^\prime_{ik}=-\sinh \theta_{ik}$, and
$c^\prime_{ik}\ast c_{ijk}=-\sinh q_j$,
it follows that $\cosh h_i=\sinh q_j\cdot\sinh \theta_{ik}$.

In the generalized right-angled hyperbolic triangle $c_{ijk}c^\prime_{ij}v^\prime_i$,
we have $-v^\prime_i\ast c_{ijk}=(c^\prime_{ij}\ast c_{ijk})(v^\prime_i\ast c^\prime_{ij})$.
Since $v^\prime_i\ast c_{ijk}=-\cosh h_i$,
$v^\prime_i\ast c^\prime_{ij}=-\sinh \theta_{ij}$, and
$c^\prime_{ij}\ast c_{ijk}=-\sinh q_k$,
it follows that $\cosh h_i=\sinh q_k\cdot\sinh \theta_{ij}$.

In the generalized right-angled hyperbolic triangle $c_{ijk}c^\prime_{jk}v^\prime_j$,
we have $-v^\prime_j\ast c_{ijk}=(c^\prime_{jk}\ast c_{ijk})(v^\prime_j\ast c^\prime_{jk})$.
By Proposition \ref{Prop: Ratcliffe 1} (iii),
$v^\prime_j\ast c_{ijk}
=-\cosh(d_\mathbb{H}((v^\prime_j)^\bot,c_{ijk}^\bot))
=-\cosh h_j$.
Since $c^\prime_{jk}$ and $c_{ijk}$ lie on opposite sides of the hyperplane $c_{ijk}^\bot$, 
it follows that $c^\prime_{jk}\ast c_{ijk}=-\sinh(d_\mathbb{H}(c^\prime_{jk},c_{ijk}^\bot))=-\sinh (-q_i)=+\sinh q_i$.
Due to $v^\prime_j\ast c^\prime_{jk}=-\sinh \theta_{jk}$,
we conclude that $\cosh h_j=-\sinh q_i\cdot\sinh \theta_{jk}$.

In the generalized right-angled hyperbolic triangle $c_{ijk}c^\prime_{ij}v^\prime_j$,
we have $-v^\prime_j\ast c_{ijk}=(c^\prime_{ij}\ast c_{ijk})(v^\prime_j\ast c^\prime_{ij})$.
Since $v^\prime_j\ast c_{ijk}=-\cosh h_j$,
$v^\prime_j\ast c^\prime_{ij}=-\sinh \theta_{ji}$, and
$c^\prime_{ij}\ast c_{ijk}=-\sinh q_k$,
it follows that $\cosh h_j=\sinh q_k\cdot\sinh \theta_{ji}$.

In the generalized right-angled hyperbolic triangle $c_{ijk}c^\prime_{jk}v^\prime_k$,
we have $-v^\prime_k\ast c_{ijk}=(c^\prime_{jk}\ast c_{ijk})(v^\prime_k\ast c^\prime_{jk})$.
Since $v^\prime_k\ast c_{ijk}=-\cosh h_k$,
$v^\prime_k\ast c^\prime_{jk}=-\sinh \theta_{kj}$, and
$c^\prime_{jk}\ast c_{ijk}=\sinh q_i$,
it follows that $\cosh h_k=-\sinh q_i\cdot\sinh \theta_{kj}$.

In the generalized right-angled hyperbolic triangle $c_{ijk}c^\prime_{ik}v^\prime_k$,
we have $-v^\prime_k\ast c_{ijk}=(c^\prime_{ik}\ast c_{ijk})(v^\prime_k\ast c^\prime_{ik})$.
Since $v^\prime_k\ast c_{ijk}=-\cosh h_k$,
$v^\prime_k\ast c^\prime_{ik}=-\sinh \theta_{ki}$, and
$c^\prime_{ik}\ast c_{ijk}=-\sinh q_j$,
it follows that$\cosh h_k=\sinh q_j\cdot\sinh \theta_{ki}$.

\item[(2)] Assume $c_{ijk}\in D_\mathrm{II}$.
In the generalized right-angled hyperbolic triangle $c_{ijk}c_{ik}v_i$,
we have $-v_i\ast c_{ijk}=(c_{ik}\ast c_{ijk})(v_i\ast c_{ik})$.
By Proposition \ref{Prop: Ratcliffe 1} (iii), $v_i\ast c_{ijk}=-\cosh q_i$ and $c_{ik}\ast c_{ijk}=-\sinh h_j$.
Since $v_i\ast c_{ik}=-\sinh d_\mathbb{H}(v_i^\bot,c_{ik})=-\sinh(+d_{ik})=-\sinh d_{ik}$,
it follows that $\cosh q_i=\sinh h_j\cdot\sinh d_{ik}$.

In the generalized right-angled hyperbolic triangle $c_{ijk}c_{ij}v_i$,
we have $-v_i\ast c_{ijk}=(c_{ij}\ast c_{ijk})(v_i\ast c_{ij})=\cosh q_i$.
Since $v_i\ast c_{ij}=-\sinh d_{ij}$ and
$c_{ij}\ast c_{ijk}=-\sinh h_k$,
it follows that $\cosh q_i=\sinh h_k\cdot\sinh d_{ij}$.

In the generalized right-angled hyperbolic triangle $c_{ijk}c_{jk}v_j$,
we have $-v_j\ast c_{ijk}=(c_{jk}\ast c_{ijk})(v_j\ast c_{jk})$.
By Proposition \ref{Prop: Ratcliffe 1} (iii),
$v_j\ast c_{ijk}
=-\cosh(d_\mathbb{H}(v_j^\bot,c_{ijk}^\bot))
=-\cosh q_j$.
Since $c_{jk}$ and $c_{ijk}$ lie on opposite sides of the hyperplane $c_{ijk}^\bot$, 
it follows that $c_{jk}\ast c_{ijk}=-\sinh(d_\mathbb{H}(c_{jk},c_{ijk}^\bot))
=-\sinh (-h_i)=+\sinh h_i$.
Due to $v_j\ast c_{jk}=-\sinh d_{jk}$,
we conclude that $\cosh q_j=-\sinh h_i\cdot \sinh d_{jk}$.

In the generalized right-angled hyperbolic triangle $c_{ijk}c_{ij}v_j$,
we have $-v_j\ast c_{ijk}=(c_{ij}\ast c_{ijk})(v_j\ast c_{ij})$.
Since $v_j\ast c_{ijk}=-\cosh q_j$,
$v_j\ast c_{ij}=-\sinh d_{ji}$,
and $c_{ij}\ast c_{ijk}=-\sinh h_k$,
it follows that $\cosh q_j=\sinh h_k\cdot \sinh d_{ji}$.

In the generalized right-angled hyperbolic triangle $c_{ijk}c_{jk}v_k$,
we have $-v_k\ast c_{ijk}=(c_{jk}\ast c_{ijk})(v_k\ast c_{jk})$.
Since $v_k\ast c_{ijk}=-\cosh q_k$,
$v_k\ast c_{jk}=-\sinh d_{kj}$,
and $c_{jk}\ast c_{ijk}=\sinh h_i$,
it follows that $\cosh q_k=-\sinh h_i\cdot \sinh d_{kj}$.

In the generalized right-angled hyperbolic triangle $c_{ijk}c_{ik}v_k$,
we have $-v_k\ast c_{ijk}=(c_{ik}\ast c_{ijk})(v_k\ast c_{ik})$.
Since $v_k\ast c_{ijk}=-\cosh q_k$,
$v_k\ast c_{ik}=-\sinh d_{ki}$, and
$c_{ik}\ast c_{ijk}=-\sinh h_j$,
it follows that $\cosh q_k=\sinh h_j\cdot\sinh d_{ki}$.

In the generalized right-angled hyperbolic triangle $c_{ijk}c^\prime_{ik}v^\prime_i$,
we have $-v^\prime_i\ast c_{ijk}=(c^\prime_{ik}\ast c_{ijk})(v^\prime_i\ast c^\prime_{ik})$.
Similarly, $v^\prime_i\ast c_{ijk}=\cosh h_i$ and
$v^\prime_i\ast c^\prime_{ik}=-\sinh \theta_{ik}$.
If $c^\prime_{ik}$ and $c_{ijk}$ lie on the same side of the hyperplane $c_{ijk}^\bot$,
then $c^\prime_{ik}\ast c_{ijk}=+\sinh(d_\mathbb{H}(c^\prime_{ik},c_{ijk}^\bot))
=+\sinh (-q_j)=-\sinh q_j$.
If $c^\prime_{ik}$ and $c_{ijk}$ lie on opposite sides of the hyperplane $c_{ijk}^\bot$, then $c^\prime_{ij}\ast c_{ijk}=-\sinh(d_\mathbb{H}(c^\prime_{ik},c_{ijk}^\bot))
=-\sinh (+q_j)=-\sinh q_j$.
In either case, $c^\prime_{ik}\ast c_{ijk}=-\sinh q_j$.
Consequently, $\cosh h_i=-\sinh q_j\cdot\sinh \theta_{ik}$.
By Lemma \ref{Lem: coplaner 2}, it is easy to check that $\theta_{ik}<0$.
It follows that $q_j>0$ and hence $d_\mathbb{H}(c^\prime_{ik},c_{ijk}^\bot)=q_j$.
This implies that $c^\prime_{ik}$ and $c_{ijk}$ lie on opposite sides of the hyperplane $c_{ijk}^\bot$.

In the generalized right-angled hyperbolic triangle $c_{ijk}c^\prime_{ij}v^\prime_i$,
we have $-v^\prime_i\ast c_{ijk}=(c^\prime_{ij}\ast c_{ijk})(v^\prime_i\ast c^\prime_{ij})$.
Similarly, $v^\prime_i\ast c_{ijk}=\cosh h_i$ and $v^\prime_i\ast c^\prime_{ij}=-\sinh \theta_{ij}$.
If $c^\prime_{ij}$ and $c_{ijk}$ lie on the same side of the hyperplane $c_{ijk}^\bot$,
then $c^\prime_{ij}\ast c_{ijk}=+\sinh(d_\mathbb{H}(c^\prime_{ij},c_{ijk}^\bot))
=+\sinh (-q_k)=-\sinh q_k$.
If $c^\prime_{ij}$ and $c_{ijk}$ lie on opposite sides of the hyperplane $c_{ijk}^\bot$, then $c^\prime_{ij}\ast c_{ijk}=-\sinh(d_\mathbb{H}(c^\prime_{ij},c_{ijk}^\bot))
=-\sinh (+q_k)=-\sinh q_k$.
In either case, $c^\prime_{ij}\ast c_{ijk}=-\sinh q_k$.
Consequently, $\cosh h_i=-\sinh q_k\cdot\sinh \theta_{ij}$.
Similarly, $\theta_{ij}<0$ by Lemma \ref{Lem: coplaner 2} and hence $q_k>0$.
Then $d_\mathbb{H}(c^\prime_{ij},c_{ijk}^\bot)=q_k$, which implies that $c^\prime_{ij}$ and $c_{ijk}$ are on opposite sides of the hyperplane $c_{ijk}^\bot$.

In the generalized right-angled hyperbolic triangle $c_{ijk}c^\prime_{jk}v^\prime_j$,
we have $-v^\prime_j\ast c_{ijk}=(c^\prime_{jk}\ast c_{ijk})(v^\prime_j\ast c^\prime_{jk})$.
Since $v^\prime_j\ast c_{ijk}=-\cosh h_j$,
$v^\prime_j\ast c^\prime_{jk}=-\sinh \theta_{jk}$,
and $c^\prime_{jk}\ast c_{ijk}=-\sinh q_i$,
it follows that $\cosh h_j=\sinh q_i\cdot\sinh \theta_{jk}$.

In the generalized right-angled hyperbolic triangle $c_{ijk}c^\prime_{ij}v^\prime_j$,
we have $-v^\prime_j\ast c_{ijk}=(c^\prime_{ij}\ast c_{ijk})(v^\prime_j\ast c^\prime_{ij})$.
Since $v^\prime_j\ast c_{ijk}=-\cosh h_j$,
$v^\prime_j\ast c^\prime_{ij}=-\sinh \theta_{ji}$, and
$c^\prime_{ij}\ast c_{ijk}=-\sinh q_k$,
it follows that $\cosh h_j=\sinh q_k\cdot\sinh \theta_{ji}$.

In the generalized right-angled hyperbolic triangle $c_{ijk}c^\prime_{jk}v^\prime_k$,
we have $-v^\prime_k\ast c_{ijk}=(c^\prime_{jk}\ast c_{ijk})(v^\prime_k\ast c^\prime_{jk})$.
Since $v^\prime_k\ast c_{ijk}=-\cosh h_k$,
$v^\prime_k\ast c^\prime_{jk}=-\sinh \theta_{kj}$, and
$c^\prime_{jk}\ast c_{ijk}=-\sinh q_i$,
it follows that $\cosh h_k=\sinh q_i\cdot\sinh \theta_{kj}$.

In the generalized right-angled hyperbolic triangle $c_{ijk}c^\prime_{ik}v^\prime_k$,
we have $-v^\prime_k\ast c_{ijk}=(c^\prime_{ik}\ast c_{ijk})(v^\prime_k\ast c^\prime_{ik})$.
Since $v^\prime_k\ast c_{ijk}=-\cosh h_k$,
$v^\prime_k\ast c^\prime_{ik}=-\sinh \theta_{ki}$, and
$c^\prime_{ik}\ast c_{ijk}=-\sinh q_j$,
it follows that $\cosh h_k=\sinh q_j\cdot\sinh \theta_{ki}$.
\end{description}
\qed

\begin{proposition}
Let $c_{ijk}\in \mathbb{R}^3$ lie in the open domain $D_\mathrm{vi}$, then
\begin{equation*}
\begin{aligned}
&\cosh q_i=\sinh h_j\cdot\sinh d_{ik}
=-\sinh h_k\cdot\sinh d_{ij},\\
&\cosh q_j=\sinh h_i\cdot\sinh d_{jk}
=\sinh h_k\cdot\sinh d_{ji},\\
&\cosh q_k=\sinh h_i\cdot\sinh d_{kj}
=-\sinh h_j\cdot\sinh d_{ki},\\
&\cosh h_i=\sinh q_j\cdot\sinh \theta_{ik}
=\sinh q_k\cdot\sinh \theta_{ij},\\
&\cosh h_j=\sinh q_i\cdot\sinh \theta_{jk}
=-\sinh q_k\cdot\sinh \theta_{ji},\\
&\cosh h_k=-\sinh q_i\cdot\sinh \theta_{kj}
=\sinh q_j\cdot\sinh \theta_{ki}.
\end{aligned}
\end{equation*}
\end{proposition}
\proof
In the generalized right-angled hyperbolic triangle $c_{ijk}c_{ik}v_i$,
we have $-v_i\ast c_{ijk}=(c_{ik}\ast c_{ijk})(v_i\ast c_{ik})$.
By Proposition \ref{Prop: Ratcliffe 1} (iii),
$v_i\ast c_{ijk}=\cosh d_\mathbb{H}(v_i^\bot,c_{ijk}^\bot)=\cosh q_i$.
By Proposition \ref{Prop: Ratcliffe 1} (ii),
$v_i\ast c_{ik}=+\sinh d_\mathbb{H}(v_i^\bot,c_{ik})=+\sinh(-d_{ik})=-\sinh d_{ik}$.
If $c_{ik}$ and $c_{ijk}$ lie on the same side of the hyperplane $c_{ijk}^\bot$,
then $c_{ik}\ast c_{ijk}=+\sinh(d_\mathbb{H}(c_{ik},c_{ijk}^\bot))
=+\sinh h_j$;
if $c_{ik}$ and $c_{ijk}$ lie on opposite sides of the hyperplane $c_{ijk}^\bot$, then $c_{ik}\ast c_{ijk}=-\sinh(d_\mathbb{H}(c_{ik},c_{ijk}^\bot))
=-\sinh (-h_j)=+\sinh h_j$.
In either case, $c_{ik}\ast c_{ijk}=\sinh h_j$.
Consequently, $\cosh q_i=\sinh h_j\cdot\sinh d_{ik}$.
By Lemma \ref{Lem: coplaner}, it is easy to check that $d_{ik}<0$.
It follows that $h_j<0$ and hence
$d_\mathbb{H}(c_{ik},c_{ijk}^\bot)=-h_j$.
This implies that $c_{ik}$ and $c_{ijk}$ lie on opposite sides of $c_{ijk}^\bot$.

In the generalized right-angled hyperbolic triangle $c_{ijk}c_{ij}v_i$,
we have $-v_i\ast c_{ijk}=(c_{ij}\ast c_{ijk})(v_i\ast c_{ij})$.
Since $v_i\ast c_{ij}=-\sinh d_{ij}$ and
$c_{ij}\ast c_{ijk}=-\sinh h_k$,
it follows that $\cosh q_i=-\sinh h_k\cdot\sinh d_{ij}$.

In the generalized right-angled hyperbolic triangle $c_{ijk}c_{jk}v_j$,
we have $-v_j\ast c_{ijk}=(c_{jk}\ast c_{ijk})(v_j\ast c_{jk})$.
Since $v_j\ast c_{ijk}=-\cosh q_j$,
$v_j\ast c_{jk}=-\sinh d_{jk}$,
and $c_{jk}\ast c_{ijk}=-\sinh h_i$,
it follows that $\cosh q_j=\sinh h_i\cdot \sinh d_{jk}$.

In the generalized right-angled hyperbolic triangle $c_{ijk}c_{ij}v_j$,
we have $-v_j\ast c_{ijk}=(c_{ij}\ast c_{ijk})(v_j\ast c_{ij})$.
Since $v_j\ast c_{ijk}=-\cosh q_j$,
$v_j\ast c_{ij}=-\sinh d_{ji}$,
and $c_{ij}\ast c_{ijk}=-\sinh h_k$,
it follows that $\cosh q_j=\sinh h_k\cdot \sinh d_{ji}$.

In the generalized right-angled hyperbolic triangle $c_{ijk}c_{jk}v_k$,
we have $-v_k\ast c_{ijk}=(c_{jk}\ast c_{ijk})(v_k\ast c_{jk})$.
Since $v_k\ast c_{ijk}=-\cosh q_k$,
$v_k\ast c_{jk}=-\sinh d_{kj}$,
and $c_{jk}\ast c_{ijk}=-\sinh h_i$,
it follows that $\cosh q_k=\sinh h_i\cdot \sinh d_{kj}$.

In the generalized right-angled hyperbolic triangle $c_{ijk}c_{ik}v_k$,
we have $-v_k\ast c_{ijk}=(c_{ik}\ast c_{ijk})(v_k\ast c_{ik})$.
Since $v_k\ast c_{ijk}=-\cosh q_k$,
$v_k\ast c_{ik}=-\sinh d_{ki}$, and
$c_{ik}\ast c_{ijk}=\sinh h_j$,
it follows that $\cosh q_k=-\sinh h_j\cdot\sinh d_{ki}$.

In the generalized right-angled hyperbolic triangle $c_{ijk}c^\prime_{ik}v^\prime_i$,
we have $-v^\prime_i\ast c_{ijk}=(c^\prime_{ik}\ast c_{ijk})(v^\prime_i\ast c^\prime_{ik})$.
Since $v^\prime_i\ast c_{ijk}=-\cosh h_i$,
$v^\prime_i\ast c^\prime_{ik}=-\sinh \theta_{ik}$, and
$c^\prime_{ik}\ast c_{ijk}=-\sinh q_j$,
it follows that $\cosh h_i=\sinh q_j\cdot\sinh \theta_{ik}$.

In the generalized right-angled hyperbolic triangle $c_{ijk}c^\prime_{ij}v^\prime_i$,
we have $-v^\prime_i\ast c_{ijk}=(c^\prime_{ij}\ast c_{ijk})(v^\prime_i\ast c^\prime_{ij})$.
Since $v^\prime_i\ast c_{ijk}=-\cosh h_i$,
$v^\prime_i\ast c^\prime_{ij}=-\sinh \theta_{ij}$, and
$c^\prime_{ij}\ast c_{ijk}=-\sinh q_k$,
it follows that $\cosh h_i=\sinh q_k\cdot\sinh \theta_{ij}$.

In the generalized right-angled hyperbolic triangle $c_{ijk}c^\prime_{jk}v^\prime_j$,
we have $-v^\prime_j\ast c_{ijk}=(c^\prime_{jk}\ast c_{ijk})(v^\prime_j\ast c^\prime_{jk})$.
Moreover, $v^\prime_j\ast c_{ijk}=\cosh h_j$ and
$v^\prime_j\ast c^\prime_{jk}=-\sinh \theta_{jk}$.
If $c^\prime_{jk}$ and $c_{ijk}$ lie on the same side of the hyperplane $c_{ijk}^\bot$,
then $c^\prime_{jk}\ast c_{ijk}=+\sinh(d_\mathbb{H}(c^\prime_{jk},c_{ijk}^\bot))
=+\sinh q_i$;
if $c^\prime_{jk}$ and $c_{ijk}$ lie on opposite sides of the hyperplane $c_{ijk}^\bot$, then $c^\prime_{jk}\ast c_{ijk}=-\sinh(d_\mathbb{H}(c^\prime_{jk},c_{ijk}^\bot))
=-\sinh (-q_i)=+\sinh q_i$.
In either case, $c^\prime_{jk}\ast c_{ijk}=\sinh q_i$.
Consequently, $\cosh h_j=\sinh q_i\cdot\sinh \theta_{jk}$.
By Lemma \ref{Lem: coplaner 2}, it is easy to check that $\theta_{jk}<0$.
It follows that $q_i<0$ and hence
$d_\mathbb{H}(c^\prime_{jk},c_{ijk}^\bot)=-q_i$.
This implies that $c^\prime_{jk}$ and $c_{ijk}$ lie on opposite sides of the hyperplane $c_{ijk}^\bot$.

In the generalized right-angled hyperbolic triangle $c_{ijk}c^\prime_{ij}v^\prime_j$,
we have $-v^\prime_j\ast c_{ijk}=(c^\prime_{ij}\ast c_{ijk})(v^\prime_j\ast c^\prime_{ij})$.
Since $v^\prime_j\ast c_{ijk}=\cosh h_j$,
$v^\prime_j\ast c^\prime_{ij}=-\sinh \theta_{ji}$, and
$c^\prime_{ij}\ast c_{ijk}=-\sinh q_k$,
it follows that $\cosh h_j=-\sinh q_k\cdot\sinh \theta_{ji}$.

In the generalized right-angled hyperbolic triangle $c_{ijk}c^\prime_{jk}v^\prime_k$,
we have $-v^\prime_k\ast c_{ijk}=(c^\prime_{jk}\ast c_{ijk})(v^\prime_k\ast c^\prime_{jk})$.
Since $v^\prime_k\ast c_{ijk}=-\cosh h_k$,
$v^\prime_k\ast c^\prime_{jk}=-\sinh \theta_{kj}$, and
$c^\prime_{jk}\ast c_{ijk}=\sinh q_i$,
it follows that $\cosh h_k=-\sinh q_i\cdot\sinh \theta_{kj}$.

In the generalized right-angled hyperbolic triangle $c_{ijk}c^\prime_{ik}v^\prime_k$,
we have $-v^\prime_k\ast c_{ijk}=(c^\prime_{ik}\ast c_{ijk})(v^\prime_k\ast c^\prime_{ik})$.
Since $v^\prime_k\ast c_{ijk}=-\cosh h_k$,
$v^\prime_k\ast c^\prime_{ik}=-\sinh \theta_{ki}$, and
$c^\prime_{ik}\ast c_{ijk}=-\sinh q_j$,
it follows that $\cosh h_k=\sinh q_j\cdot\sinh \theta_{ki}$.
\qed

\end{document}